\documentclass[12pt]{article}
\usepackage{geometry}
\usepackage{setspace}			
\usepackage{authblk}			
\usepackage{appendix}			
\usepackage{natbib}			    
\usepackage[hidelinks]{hyperref}	    

\usepackage{amsmath,amsfonts,amssymb}	
\usepackage{bm}				            
\usepackage{bbm}			
\usepackage[bb=stix]{mathalpha}	

\usepackage{graphicx}			
\usepackage{caption}			
\usepackage{subcaption}			
\usepackage{float}				
\usepackage[section]{placeins}	

\usepackage{longtable}			
\usepackage{booktabs}			
\usepackage{multirow}			
\usepackage{multicol}			

\usepackage{color}				
\usepackage{ulem}				
\usepackage{cancel}				

\begin{document}

\title{Performance and Risk Analytics of Asian Exchange-Traded Funds}

\author[1]{Bhathiya Divelgama\thanks{Corresponding author, bdivelga@ttu.edu}}
\author[1]{Nancy Asare Nyarko}
\author[1]{Naa Sackley Dromo Aryee}
\author[2]{Abootaleb Shirvani}
\author[1]{Svetlozar T. Rachev}
\affil[1]{Department of Mathematics \& Statistics, Texas Tech University, Lubbock, TX, USA}
\affil[2]{Department of Mathematical Science, Kean University, Union, NJ, USA}

\date{}
\maketitle
\begin{abstract}
Investing in Asian markets through exchange-traded funds (ETFs) provides investors with access to rapidly expanding economies and valuable diversification opportunities. 
This study examines the advantages and challenges of investing in Asian ETFs by conducting comprehensive risk assessments, 
portfolio analyses, and performance comparisons. 
The dataset comprises 29 ETFs offering exposure across a wide spectrum of Asian markets, 
including broad regional funds, 
country-specific ETFs, as well as sector-focused funds, 
dividend-oriented ETFs, small-cap portfolios, 
and emerging market bond ETFs.

To evaluate risk and return dynamics, 
the study employs Markowitz’s efficient frontier to identify optimal portfolios for given levels of risk, 
and conditional value-at-risk (CVaR) to capture potential extreme losses for a more comprehensive risk assessment. 
Multiple portfolio configurations are analyzed under long-only and long–short investment strategies to assess adaptability across varying market conditions. Furthermore, 
key performance risk measures including the Sharpe ratio, Rachev ratio, 
and stable tail-adjusted return ratio (STARR) are calculated to provide an 
in-depth evaluation of reward-to-risk efficiency, 
with particular emphasis on the role of tail behavior in portfolio performance.

This research aims to deliver deeper insights into the risk–return trade-offs, 
tail-risk behavior, and efficiency of Asian ETFs, 
offering investors a practical foundation for constructing robust and 
well-diversified portfolios across both emerging and developed Asian markets.
\end{abstract}

\section{Introduction}

Exchange-traded funds (ETFs) are open-end index funds that provide daily transparency of their holdings and trade on exchanges like stocks, 
offering investors real-time pricing. 
The first ETF was launched in 1990 with the Toronto exchange traded index participation to track the TSE-35, 
and early research focused on their structure, 
tax efficiency, 
and unique creation-rede ion mechanism \citep{mittal2017}. 
Over the past quarter-century, 
ETFs have emerged as a highly popular passive investment vehicle due to increasing demand for passive strategies, 
low transaction costs, and high liquidity, 
which distinguish them from traditional index funds \citep{liebi2020}. 
By 2016, ETFs accounted for over 10\% of U.S. market capitalization and more than 30\% of overall trading volume \citep{Ben-David2016}.

Since their inception, 
ETFs have attracted substantial academic attention, 
resulting in a growing body of literature examining their performance, 
market dynamics, and international expansion. 
Research has investigated ETFs in diverse regions, 
including Japan \citep{rompotis2024b}, 
India \citep{malhotra2023}, China \citep{wu2021, ning2024}, 
Russia \citep{tarassov2016}, 
and the Asia-Pacific region \citep{marszk2019, rompotis2024a}, 
highlighting differences in market development, trading efficiency, 
and risk characteristics. 
Studies have documented that ETFs in emerging markets often exhibit higher trading costs, 
lower liquidity, 
and limited arbitrage efficiency compared to their developed-market counterparts \citep{hilliard2022, cheng2008, jares2004}. 
Moreover, 
empirical analyses show that ETFs can influence the tail dependence and risk contagion of underlying securities \citep{ning2024}, 
and provide resilience during market shocks \citep{malhotra2023}. 
Recent work on Pakistan-exposed ETFs further underscores the relevance of frontier markets, 
revealing that careful risk assessment and portfolio optimization are critical for capturing growth opportunities while managing exposure \citep{jaffri2025}.

Beyond geography, 
research has also examined ETF performance and efficiency across asset classes. 
Evidence suggests that, 
while ETFs are generally considered low-cost, 
passive instruments, they do not always outperform the market \citep{blitz2021etf}. 
Deviations from net asset value, 
tracking errors, and structural frictions, 
especially in international ETFs, 
can influence returns and challenge the replication of underlying indices \citep{petajisto2017, bahadar2020etf, zawadzki2020etf}. 
Systematic reviews further highlight that ETF research has expanded to cover volatility, liquidity, 
risk-return trade-offs, 
and emerging trends such as ESG, AI-driven strategies, 
and new asset-class ETFs \citep{joshi2024etf}. 
Collectively, 
these findings demonstrate that ETFs’ global reach and cross-asset presence create both opportunities and challenges for investors, 
emphasizing the importance of careful analysis when constructing diversified portfolios.

In addition to equity ETFs, 
commodity-focused ETFs have also attracted significant attention as liquid and accessible instruments for diversification and risk management. 
Building on the broader research into ETFs across markets and asset classes, 
studies on gold ETFs specifically have examined their role as safe-haven assets during periods of market stress, 
while also assessing their impact on prices and overall market stability \citep{Yenileg2025, wang2010gold}. 
This dual function underscores their importance not only for portfolio management but also for broader financial dynamics, 
highlighting how ETFs across different asset classes can influence both individual investment strategies and systemic market behavior.

Building on this literature, 
the present study focuses on Asian ETFs, 
aiming to extend existing research through a comprehensive analysis of 29 ETFs across the region. 
We assess the advantages and limitations of investing in these ETFs by conducting detailed risk assessments, 
portfolio analyses, and performance comparisons. 
The risk–return structure is evaluated using Markowitz’s efficient frontier, which identifies optimal portfolios that maximize returns for a given level of variance-based risk, 
alongside conditional value at risk (CVaR), 
which captures potential extreme losses. 
Prior studies indicate that ETF portfolios may under-perform when cost and liquidity factors are overlooked, 
however, 
accounting for structural advantages allows ETFs to outperform both broader markets and their underlying assets \citep{vu2021etfperformance}. 
Our analysis employs time-series price data, 
efficient frontier estimation, robust regression, 
and fat-tail risk assessment. 
This comprehensive approach provides new insights into the opportunities and challenges of Asian ETFs for portfolio diversification and risk management, building on earlier research in Asian ETF portfolio optimization \citep{young2019meanvariance}.

\section{Method}
\subsection{Data}
The dataset for this study was obtained from Bloomberg Professional Services, 
covering the period from December 10, 2014, to January 6, 2025. 
A total of 29 exchange-traded funds (ETFs) were selected to provide comprehensive coverage of Asian markets, 
with the inclusion of a small number of broad emerging market funds to enable comparison with the wider emerging market universe. 
Many of the Asian economies represented in this sample, 
such as China, India, South Korea, Taiwan, Malaysia, Thailand, Indonesia, the Philippines, and Vietnam, 
are classified as emerging markets by MSCI and FTSE. 
Consequently, 
the ETFs in our dataset capture both Asian regional exposures and emerging market characteristics, 
reflecting Asia’s dual role as one of the largest and fastest growing components of global emerging market indices.

The sample incorporates a diverse mix of broad regional ETFs, country specific ETFs, 
sector focused ETFs, dividend oriented funds, small cap ETFs, and fixed income ETFs. 
This composition ensures that the dataset reflects both the structural growth opportunities within Asia and the diversification benefits that arise from cross country and cross sector allocation. 
Table 1 in Appendix A presents the complete list of ETFs included in the analysis, 
along with their Bloomberg tickers, full fund names, and investment focus.

For benchmarking purposes, 
we construct an equally weighted portfolio (EQW), 
assuming a \$100 long only investment in each ETF as of December 10, 2014, 
with daily re-balancing to maintain equal weights across all 29 funds. 
In addition, the Dow Jones Industrial Average (DJIA 30) is employed as the reference market index, 
while the 3-month and 1-year U.S. Treasury yields, 
obtained from Bloomberg, 
are used as proxies for the risk-free rate over the same period.
Historical optimization is performed using a rolling window of approximately four trading years (4 × 252 trading days).

\subsection{ETF and Benchmark Performance}

\begin{figure}[H]
\centering
    \includegraphics[width=1\textwidth]{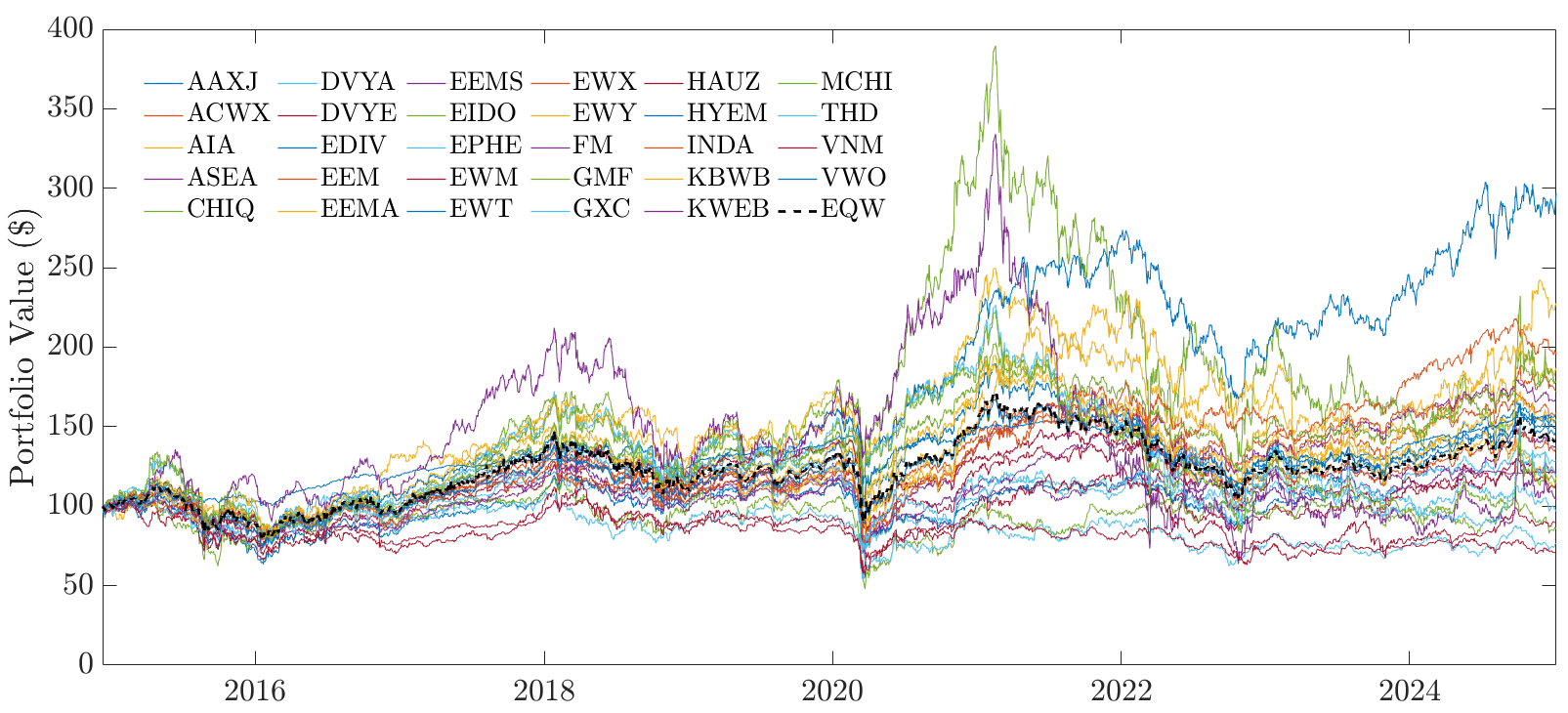}
    \caption{Cumulative prices for individual Asian-traded ETFs and EQW, assuming a \$100 investment on December 10, 2014.}
    \label{fig:cumpric}
\end{figure}

The EQW portfolio provides a robust benchmark in our analysis. 
By definition (see \cite{markowitz1952portfolio}, \cite{sharpe1999investments}), 
an EQW is a non-optimized portfolio that allocates an equal proportion of investment capital to each included asset, 
providing a neutral, diversified baseline. 
As shown in Figure~\ref{fig:cumpric}, 
the EQW demonstrates greater stability than most individual ETFs. 
While some country-specific funds such as EWT (Taiwan) and INDA (India) outperformed over the long run, others, 
including VNM (Vietnam) and HYEM (Emerging Markets High Yield Bonds), 
lagged behind. 
The COVID-19 market crash of early 2020 produced a sharp draw-down across all ETFs, 
but both the EQW and most constituents recovered, 
underscoring resilience in Asian and emerging markets.
This empirical evidence aligns with prior studies. 
\cite{jaffri2025} show that equal weighting enhances diversification and resilience in emerging market-focused ETFs. 
\cite{malladi2017} provide long-term evidence that EQW outperforms 
value-weighted strategies, 
largely due to systematic re-balancing. 
\cite{Taljaard2021} highlight that EQW can under-perform in the short run when market concentration rises but remains superior over longer horizons. 

Taken together, 
these findings support our use of the EQW portfolio as the central benchmark. 
Its simplicity, empirical robustness, 
and consistent performance advantages make it a suitable standard for evaluation.  

\section{Portfolio Optimization}

Portfolio optimization is a fundamental tool in investment management, 
aiming to construct portfolios that balance risk and return while accommodating investor preferences, 
constraints, and uncertainties. 
The concept was first formalized by \citep{markowitz1952portfolio} through the 
mean--variance framework, 
which seeks to maximize expected returns for a given level of risk, 
where risk is quantified by the portfolio’s return variance.
Over time, 
research has significantly evolved beyond the original mean variance formulation. 
\cite{salo2024fifty} provide a comprehensive historical overview of portfolio optimization, 
highlighting the development of integrated frameworks that address decision makers’ preferences and multiple objectives across planning horizons. 
Building on this, 
\cite{liu2023portfolio} introduce dynamic risk preference models and optimization under uncertain environments, 
applying average value at risk (VaR) frameworks with genetic algorithms to enhance wealth creation and efficiency.
Furthermore, 
comparative analyses such as \cite{sen2023portfolio} demonstrate that different optimization approaches, 
ranging from traditional mean variance portfolios to machine learning based designs, 
yield varying performance outcomes depending on market conditions. 
Complementing these insights, 
\cite{ou2023portfolio} emphasize the continued relevance of modern portfolio theory (MPT) in understanding how regulatory, 
industry specific, 
and client driven constraints shape portfolio performance.
 
In \cite{lindquist2021advanced}, 
both historical and dynamic optimization frameworks are discussed.
Historical optimization estimates expected returns and risks using past data, 
assuming that market relationships remain stable within each estimation window. 
A rolling-window approach is typically used, providing simplicity and interpretability but limited responsiveness to sudden market shifts or crises.
Dynamic optimization, 
on the other hand, 
models time-varying behavior in returns and volatility using methods such as ARMA–GARCH processes, 
fat-tailed (e.g., Student’s t) distributions, 
and copula-based dependence structures. 
This approach better captures evolving correlations and tail risks, 
producing more adaptive portfolios at the cost of higher model complexity and computational demand.

In this paper, 
we adopt the historical optimization approach.
Using daily ETF return data from 2014–2025, 
we use a rolling window where each window spans approximately four years. 
The optimizer estimates expected returns and covariances for assets within each window and generates portfolio weights for the next trading day. 
This recursive process produces a continuous series of optimized portfolio weights and returns, 
reflecting how the strategy would have performed under real historical conditions.

In this section, 
we first analyze the Markowitz and CVaR efficient frontiers.
The CVaR optimizations are conducted at the 95\% and 99\% quantile levels to evaluate how varying degrees of tail-risk sensitivity influence portfolio outcomes.
Next, 
we consider six optimized portfolios: 
the minimum-variance portfolio (MVP) and the tangent portfolio (TVP) derived from the Markowitz framework, 
as well as the minimum CVaR portfolios (M95 and M99) and the tangent CVaR portfolios (T95 and T99) corresponding to the 95\% and 99\% quantile levels, respectively.
In a series of experiments, 
we assess the performance of these portfolios using the equally weighted (EQW) portfolio as a benchmark under the long-only and long–short investment strategies outlined in \citep{lindquist2021advanced}.

\subsection{Markowitz and CVaR efficient frontier}
We analyze the Markowitz efficient frontier, 
which represents the set of optimal portfolios that achieve the highest expected return for a given level of risk. 
In constructing the efficient frontier, 
we use four years of historical data, 
selecting the most recent four year window available in our dataset, 
ending on January 30, 2025. 
This sample length provides a balanced window that is long enough to obtain stable return and covariance estimates while still reflecting current market conditions. Following the framework of \cite{lindquist2021advanced}, 
we replicate the mean–variance optimization procedure to estimate daily portfolio weights that minimize portfolio variance subject to a specified expected return and full investment constraint. 
This approach traces the fundamental trade-off between risk and return, 
forming the classical efficient frontier.

Formally, 
the Markowitz mean–variance optimization problem can be expressed as minimizing portfolio variance subject to a desired expected return and full investment constraint. 
Let $\mathbf{r}$ denote the vector of expected returns, 
$\boldsymbol{\Sigma}$ the covariance matrix of asset returns, 
and $\mathbf{w}$ the portfolio weights. 
The optimization problem is given by:
\[
\begin{aligned}
\min_{\mathbf{w}} \quad & \mathbf{w}^\top \boldsymbol{\Sigma} \mathbf{w} \\
\text{subject to} \quad & \mathbf{r}^\top \mathbf{w} = \bar{r}_p, \\
& \mathbf{e}_n^\top \mathbf{w} = 1,
\end{aligned}
\]
where $\bar{r}_p$ is the target portfolio return and $\mathbf{e}_n$ is an $n$-dimensional vector of ones. 
This problem is solved using the standard Lagrange multiplier method, 
which introduces auxiliary variables \( q \) and \( \theta_0 \) corresponding to the constraints \( r^T w = \bar{r}_p \) and \( e_n^T w = 1 \). 
The Lagrangian form of the optimization problem is written as:

\[
\min_{\mathbf{w}, q, \theta_0} L(\mathbf{w}, q, \theta_0) 
= \min_{\mathbf{w}, q, \theta_0} 
\left( 
\frac{1}{2}\mathbf{w}^\top \boldsymbol{\Sigma} \mathbf{w} 
+ q(\bar{r}_p - \mathbf{r}^\top \mathbf{w}) 
+ \theta_0 (1 - \mathbf{e}_n^\top \mathbf{w})
\right).
\]

Minimizing this Lagrangian with respect to $\mathbf{w}$ yields the first-order conditions that define the optimal portfolio weights. 
The analytical solution can be expressed as:

\[
\mathbf{w}^* = \bar{r}_p \mathbf{w}_1 + \mathbf{w}_2,
\]

where

\[
\mathbf{w}_1 = \frac{1}{\Delta} (B \boldsymbol{\Sigma}^{-1}\bar{\mathbf{r}} - C \boldsymbol{\Sigma}^{-1}\mathbf{e}_n), 
\quad 
\mathbf{w}_2 = \frac{1}{\Delta} (A \boldsymbol{\Sigma}^{-1}\mathbf{e}_n - C \boldsymbol{\Sigma}^{-1}\bar{\mathbf{r}}),
\]

and the scalar constants are defined as

\[
A = \bar{\mathbf{r}}^\top \boldsymbol{\Sigma}^{-1}\bar{\mathbf{r}}, 
\quad 
B = \mathbf{e}_n^\top \boldsymbol{\Sigma}^{-1}{\mathbf{e}_n}, 
\quad 
C = \bar{\mathbf{r}}^\top \boldsymbol{\Sigma}^{-1}\mathbf{e}_n, 
\quad 
\Delta = AB - C^2.
\]

The efficient frontier is then obtained by plotting the portfolio's standard deviation 
$
\sigma_p = \sqrt{\mathbf{w}^{*\top} \boldsymbol{\Sigma} \mathbf{w}^*}
$
against its expected return $\bar{r}_p$, 
illustrating the trade-off between risk and return. 
This formulation provides the basis for constructing both the minimum-variance and tangent portfolios in subsequent analyses.

\begin{figure}[H]
\centering
\begin{subfigure}{0.45\textwidth}
    \includegraphics[width=\linewidth]{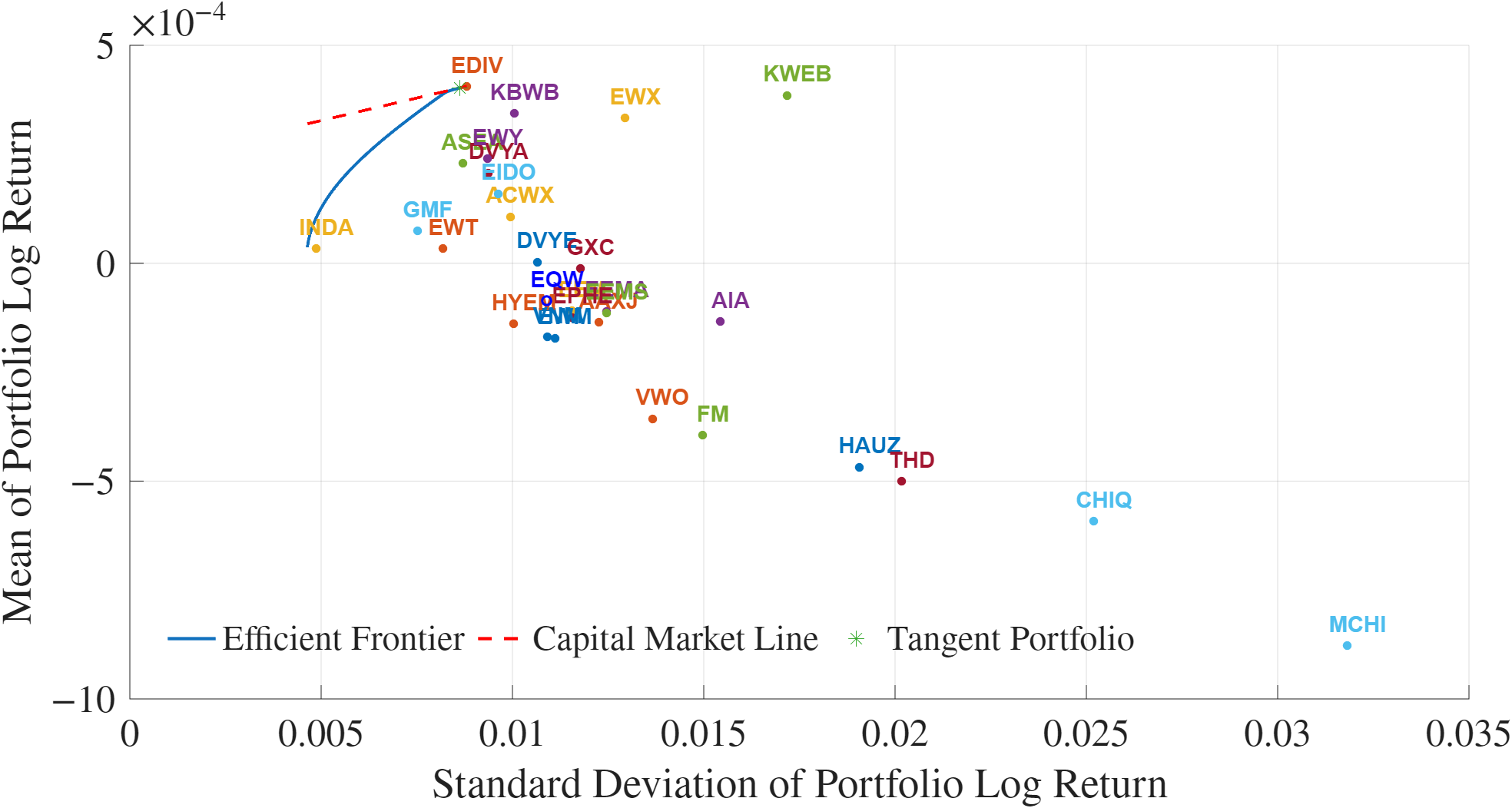}
\end{subfigure}
\hspace{0.5cm}
\begin{subfigure}{0.45\textwidth}
    \includegraphics[width=\linewidth]{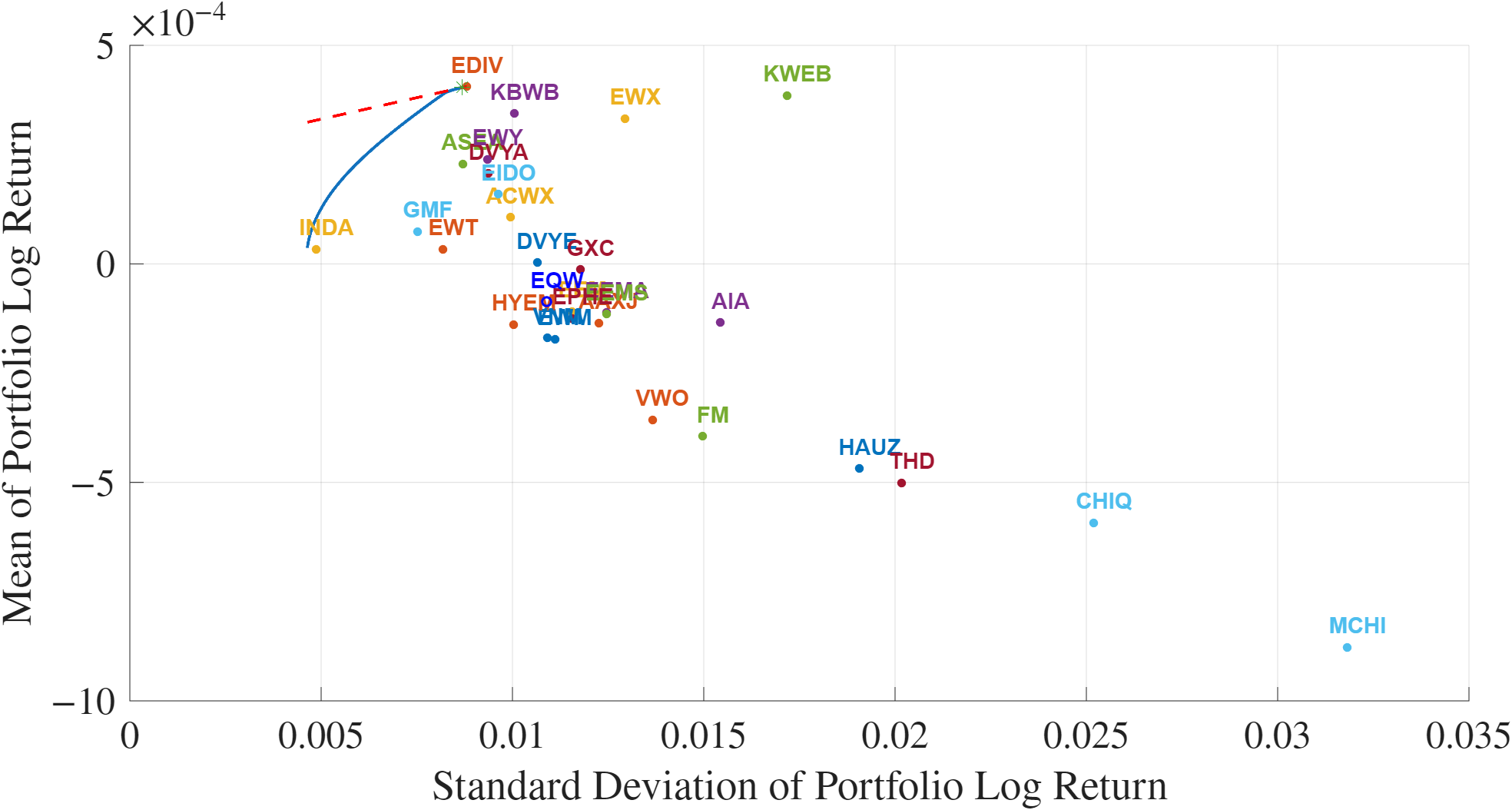}
\end{subfigure}
\caption{The Markowitz efficient frontiers using 1-year (left) and 3-month (right) U.S. Treasury yield}
\label{fig:mark_frontier}
\end{figure}

Figure \ref{fig:mark_frontier} shows the Markowitz efficient frontiers constructed using the one-year and the three-month Treasury yield as the risk-free rate. 
Each plot illustrates the trade-off between risk, 
measured by the standard deviation of portfolio log returns, 
and expected return, 
represented by the mean of portfolio log returns. 
The blue curve traces the efficient frontier, 
while the red dashed line represents the capital market line that is tangent to the frontier at the optimal, 
or tangent portfolio. 
Both graphs exhibit the typical upward-sloping convex shape, 
indicating that higher returns are associated with higher levels of portfolio risk. 
The similarity between the two frontiers suggests that the choice between the one-year and three-month Treasury yields has minimal effect on the resulting efficient frontier or the tangent portfolio position.

The EQW portfolio lies below the efficient frontier, 
exhibiting a modestly negative mean return that underscores the efficiency gains from optimization relative to a simple equal-weight allocation. 
The tangent portfolio, 
positioned near the upper-left edge of the ETF cluster, 
sits very close to EDIV, 
suggesting that this dividend-focused fund performs similarly to the optimal risk-adjusted portfolio. 
In practical terms, 
an investor holding EDIV would experience a performance profile that closely mirrors the tangent portfolio, 
indicating its strong standalone efficiency. 
Most ETFs cluster tightly around the efficient frontier, 
showing generally well-balanced risk–return profiles. 
However, 
funds such as KWEB deliver higher mean returns at the cost of greater volatility, 
while CHIQ and MCHI sit below the frontier, 
reflecting weaker risk-adjusted performance. 
Overall, 
the results suggest that diversified and income-oriented exposures dominate the tangent portfolio, 
reinforcing the advantage of combining stability and moderate growth within the Markowitz framework.

Building on the classical mean–variance framework, 
we next extend the analysis by incorporating downside tail risk through the CVAR measure, 
following the approach of \cite{rockafellar2000optimization} and \cite{krokhmal2002portfolio}. 
While the Markowitz model captures overall variance as a measure of risk, 
it assumes return distributions are symmetric. 
In contrast, 
CVaR focuses specifically on the tail of the loss distribution, 
measuring the expected loss beyond a specified quantile of returns. 
This makes it a more suitable framework for managing extreme downside risk, 
particularly in volatile or asymmetric markets.

Formally, 
the CVaR optimization problem seeks the portfolio weights $\mathbf{w}$ that minimize the expected loss in the worst $\alpha$ proportion of outcomes:
\[
\min_{\mathbf{w}} \; \mathrm{CVaR}_{\alpha}(\mathbf{w})
\quad \text{subject to} \quad 
\mathbf{r}^\top \mathbf{w} = \bar{r}_p, 
\quad 
\mathbf{e}_n^\top \mathbf{w} = 1.
\]
Here, 
$\alpha$ represents the tail probability (commonly 0.05 or 0.01), 
corresponding to 95\% or 99\% confidence levels. 
The CVAR can be expressed as
\[
\mathrm{CVaR}_{\alpha}(\mathbf{w}) 
= \mathrm{VaR}_{\alpha}(\mathbf{w}) 
+ \frac{1}{\alpha} \, 
\mathbb{E}\!\left[\, 
(-f(\mathbf{w}, \mathbf{r}) - \mathrm{VaR}_{\alpha}(\mathbf{w}))^{+} 
\right],
\]
where $(x)^{+} = \max(0, x)$ ensures that only returns exceeding the Value-at-Risk threshold contribute to the tail risk. 
This formulation transforms the problem into a convex optimization, 
allowing the CVaR frontier to be derived efficiently through linear programming techniques. 
Compared to the mean–variance frontier, 
the CVaR frontier explicitly accounts for asymmetry in return distributions and highlights the trade-off between expected return and potential extreme losses, 
offering a more realistic view of risk exposure under market stress.

\begin{figure}[H]
\centering
\begin{subfigure}{0.45\textwidth}
    \includegraphics[width=\linewidth]{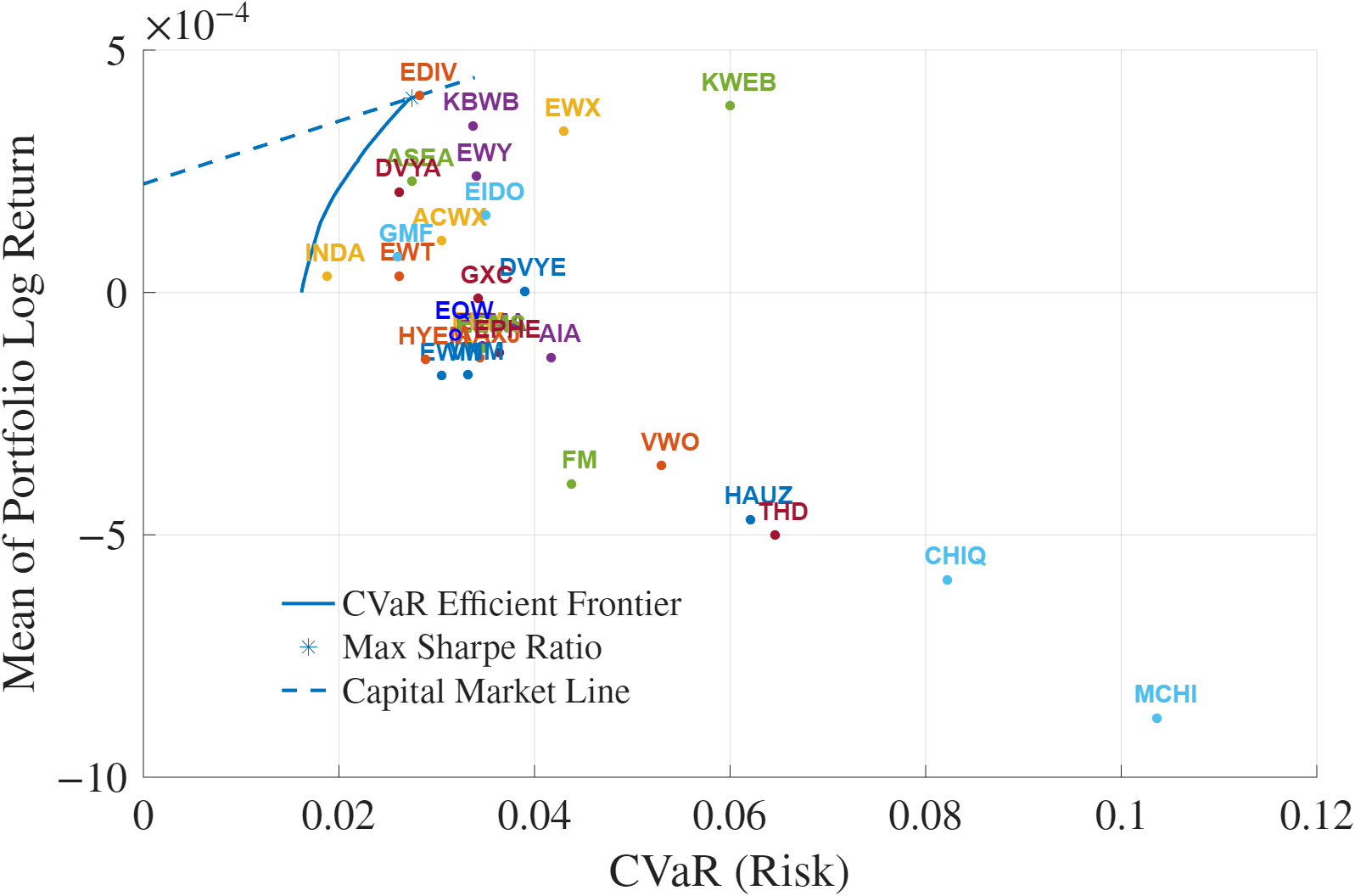}
\end{subfigure}
\hspace{0.5cm}
\begin{subfigure}{0.45\textwidth}
    \includegraphics[width=\linewidth]{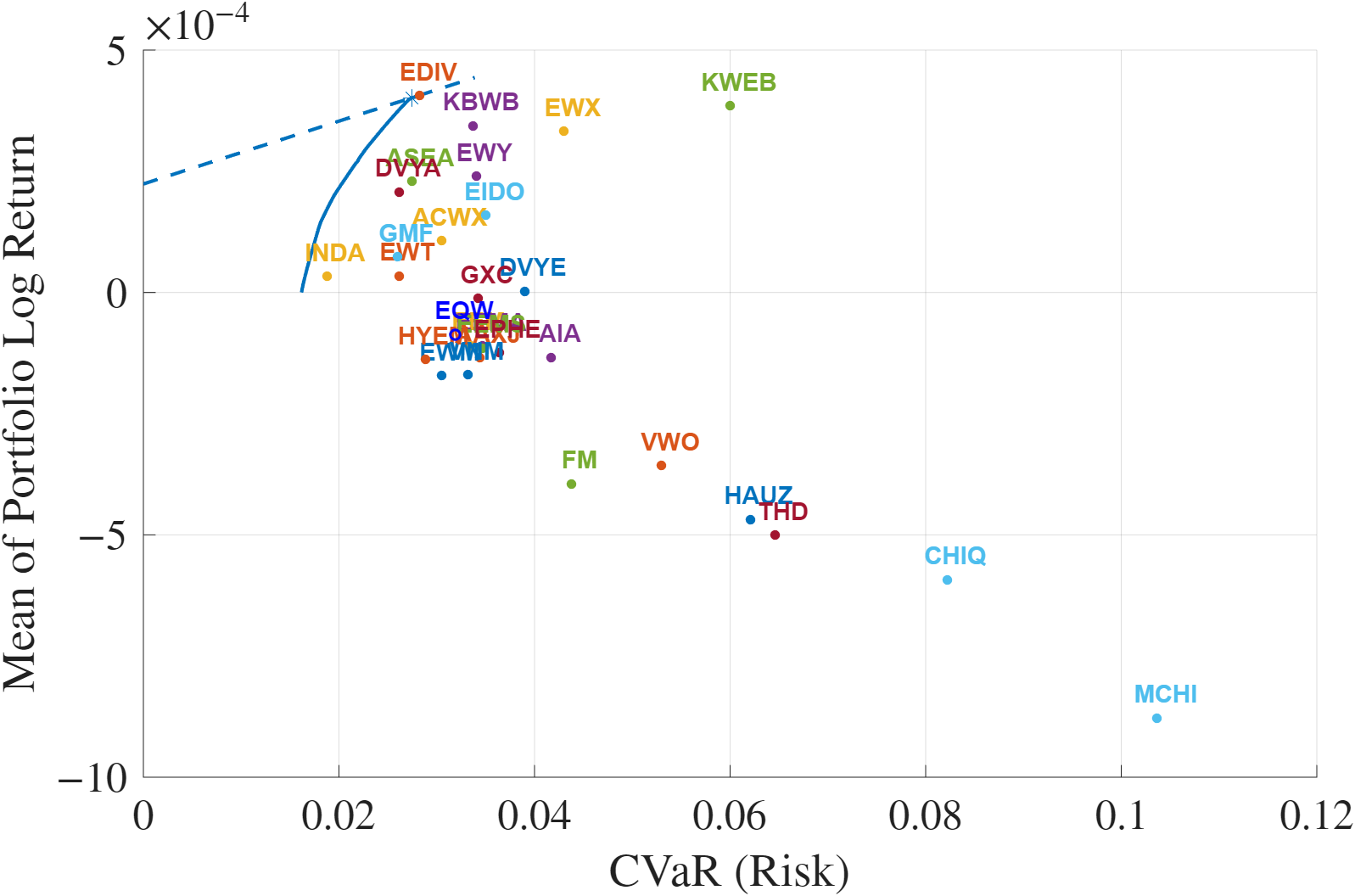}
\end{subfigure}

\vspace{0.1cm}

\begin{subfigure}{0.45\textwidth}
    \includegraphics[width=\linewidth]{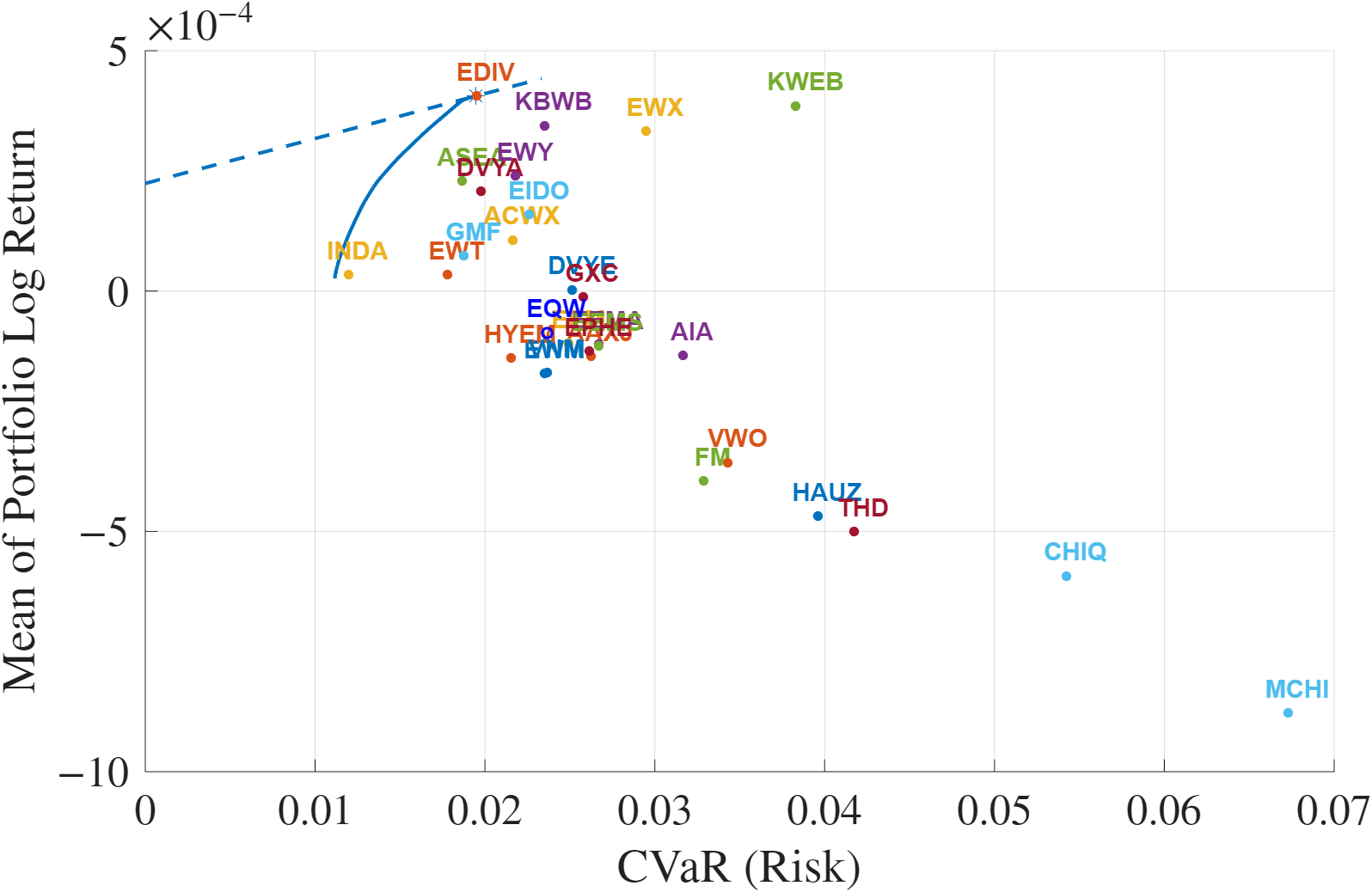}
\end{subfigure}
\hspace{0.5cm}
\begin{subfigure}{0.45\textwidth}
    \includegraphics[width=\linewidth]{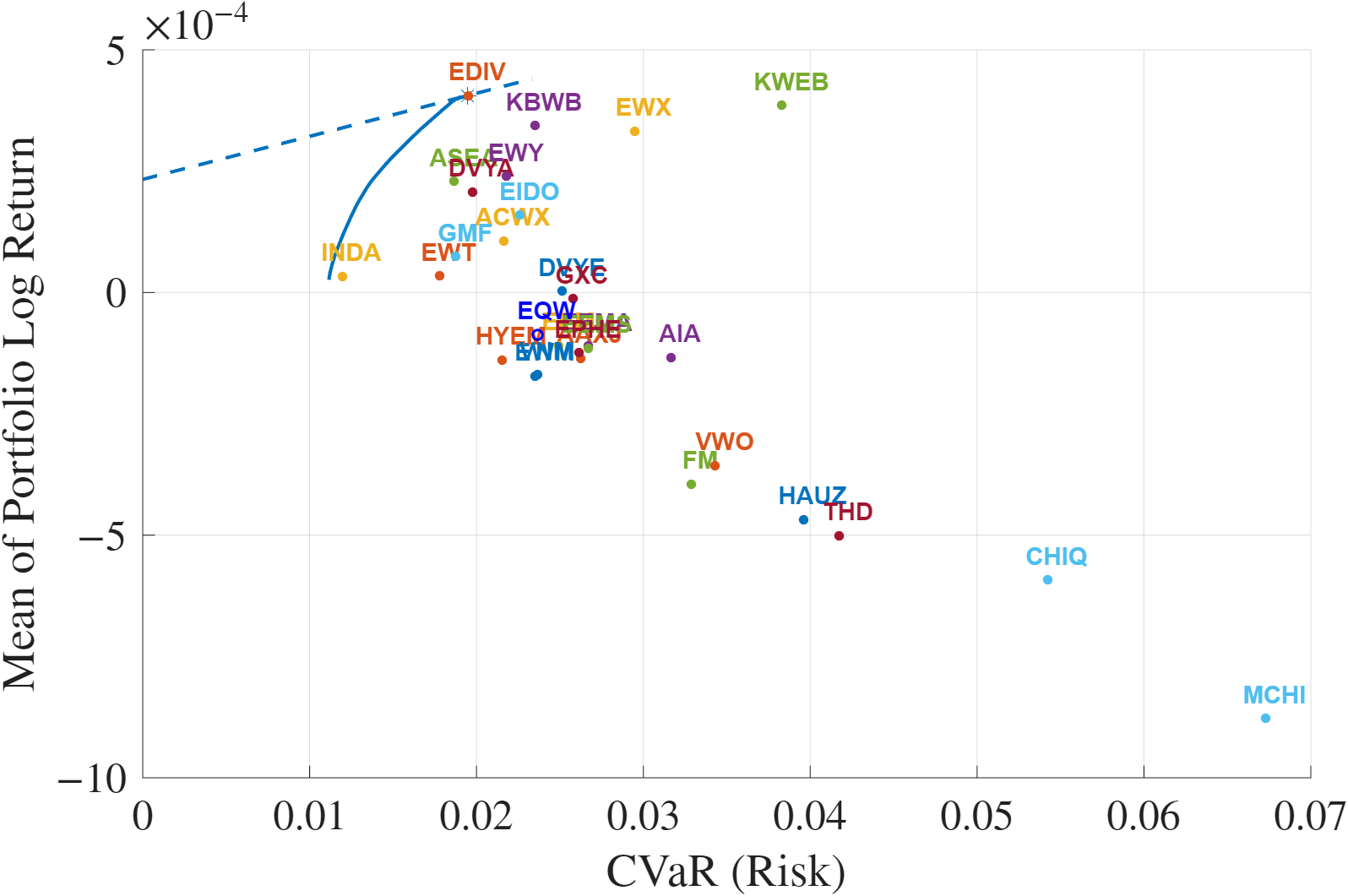}
\end{subfigure}

\caption{CVaR efficient frontiers at the 99\% and 95\% confidence levels using 1-year and 3-month U.S. Treasury yields as the risk-free rate. 
The top panels correspond to the 99\% level with the 1-year (left) and 3-month (right) yields, while the bottom panels show the 95\% level for the same yield horizons.}
\label{fig:cvar_frontier}
\end{figure}

Figure \ref{fig:cvar_frontier} presents the CVaR efficient frontiers at the 99\% and 95\% confidence levels using both the 1-year and 3-month U.S. 
Treasury yields as the risk-free rate. 
Each curve shows the trade-off between expected return and tail-risk exposure, 
measured by CVaR. 
As expected, 
the 99\% frontiers exhibit slightly higher CVaR values than the 95\%, 
reflecting the stricter loss quantile that captures rarer and more extreme downside events. 
This shift confirms that portfolios optimized under higher confidence levels face greater tail-risk for a given expected return, 
consistent with theory.

The maximum Sharpe ratio portfolios, 
which represent the best trade-off between expected return and tail-risk, 
remain relatively stable across yield horizons but differ slightly across confidence levels. 
At the 99\% level, 
the maximum Sharpe ratio portfolio lies close to EDIV, 
while at the 95\% level it aligns almost exactly with EDIV, 
suggesting that EDIV performs efficiently as a standalone ETF under both risk perspectives. 
The higher CVaR values observed at the 99\% level confirm that portfolios become more sensitive to extreme losses as the confidence threshold tightens. Overall, 
these results show that CVaR-based optimization preserves the general shape of the Markowitz frontier while emphasizing downside protection and resilience under severe market conditions.

When comparing the individual ETFs to the optimized frontiers, 
it is clear that most single funds lie below both the Markowitz and CVaR efficient frontiers, 
indicating lower risk–return efficiency. 
This suggests that while some ETFs perform reasonably well on their own, 
diversification across multiple ETFs produces portfolios that achieve higher returns for comparable or lower levels of risk.
\cite{lindquist2021advanced} noted that log returns can generally replace discrete returns in portfolio optimization when daily returns remain small. 
However, when daily returns become large, 
typically greater than about 1\%, 
the approximation can introduce significant errors. 
To confirm that this assumption held for our dataset, 
we compared the efficient frontiers computed using both log and arithmetic returns. 
The resulting frontiers were the same, indicating that our data were well within the range where the log-return approximation is valid. 
Accordingly, the analysis presented here uses log returns throughout.

\subsection{Evaluation of Portfolio Performance under Long-Only and Long–Short Strategies}
In this section, 
we evaluate the performance of optimized portfolios under both long-only and long–short investment strategies. 
The long-only approach restricts all asset weights to non-negative values, 
representing portfolios that invest exclusively in long positions. 
In contrast, 
long–short strategies allow both long and short exposures, 
enabling investors to profit from assets expected to appreciate while shorting those expected to decline. 
This flexibility can enhance return potential for a given level of volatility but introduces additional risks, 
such as margin requirements and potential losses from adverse price movements.

Following the framework of \cite{lindquist2021advanced}, 
we implement long–short models with leverage levels of 10\%, 20\%, and 30\%, 
reflecting varying degrees of short position exposure. 
The \cite{jacobs1999long} model extends the classical mean–variance optimization by allowing the long and short positions to be optimized simultaneously. 
In their formulation, 
asset weights are constrained within a range that permits shorting up to a specified proportion of the total portfolio value while maintaining overall balance in the portfolio. 
The approach assumes daily re-balancing and disregards transaction, 
margin, or borrowing costs.

\cite{lo2008new} propose an alternative long–short framework inspired by the 130/30 concept, 
in which 30\% of the capital is obtained through short selling and the proceeds are used to take 130\% long positions. 
Their model includes explicit leverage constraints to control both the extent of short selling and the overall exposure of the portfolio. 
This approach provides a more structured implementation of leverage while preserving the diversification and efficiency benefits of mean–variance optimization.

The performance of each portfolio is evaluated from December 2018 through January 2025, 
assuming a \$100 investment on December 11, 2018. 
The analysis considers six optimized portfolios: MVP, TVP, M95, T95, M99, and T99, 
representing both mean–variance and CVaR-based optimization frameworks across different confidence levels. 
The EQW portfolio serves as a benchmark for comparison. 
Although EQW is not an optimized portfolio, 
it provides a practical baseline that reflects the performance of a simple, 
fully diversified allocation without requiring return or risk estimates. 
Comparing the optimized portfolios to this benchmark highlights the performance gains achieved through systematic optimization under varying risk and leverage assumptions.

Extending the analysis to include long–short portfolios provides a clearer view of how leverage reshapes the risk–return dynamics of optimized strategies. Allowing short positions introduces flexibility to hedge downside exposures and capture relative-value opportunities, 
while leverage amplifies both gains and losses depending on market conditions. 
In this context, 
moderate leverage can improve returns by enhancing diversification and exploiting mispricings, 
whereas higher leverage increases portfolio sensitivity to volatility and draw-downs. 
Evaluating these effects across varying leverage levels (10\%, 20\%, and 30\%) highlights how the balance between risk control and return enhancement evolves once short selling is permitted.

\begin{figure}[H]
\centering

\begin{subfigure}{0.45\textwidth}
    \includegraphics[width=\linewidth]{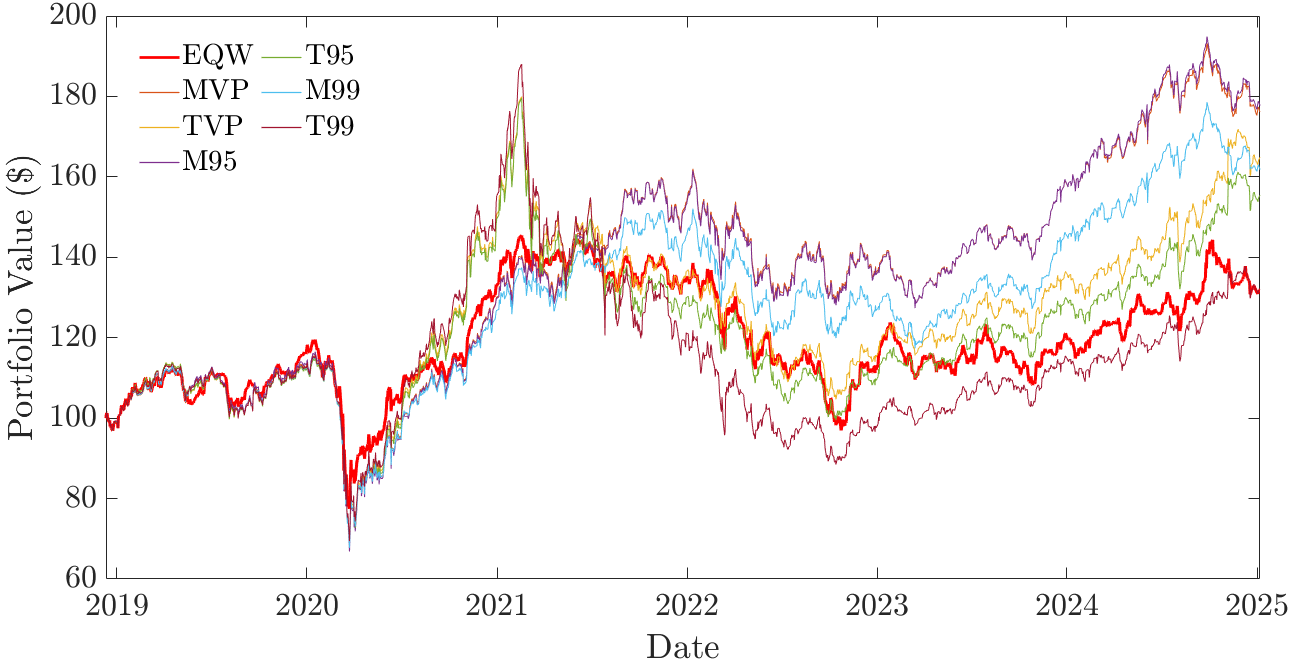}
\end{subfigure}
\hspace{0.5cm}
\begin{subfigure}{0.45\textwidth}
    \includegraphics[width=\linewidth]{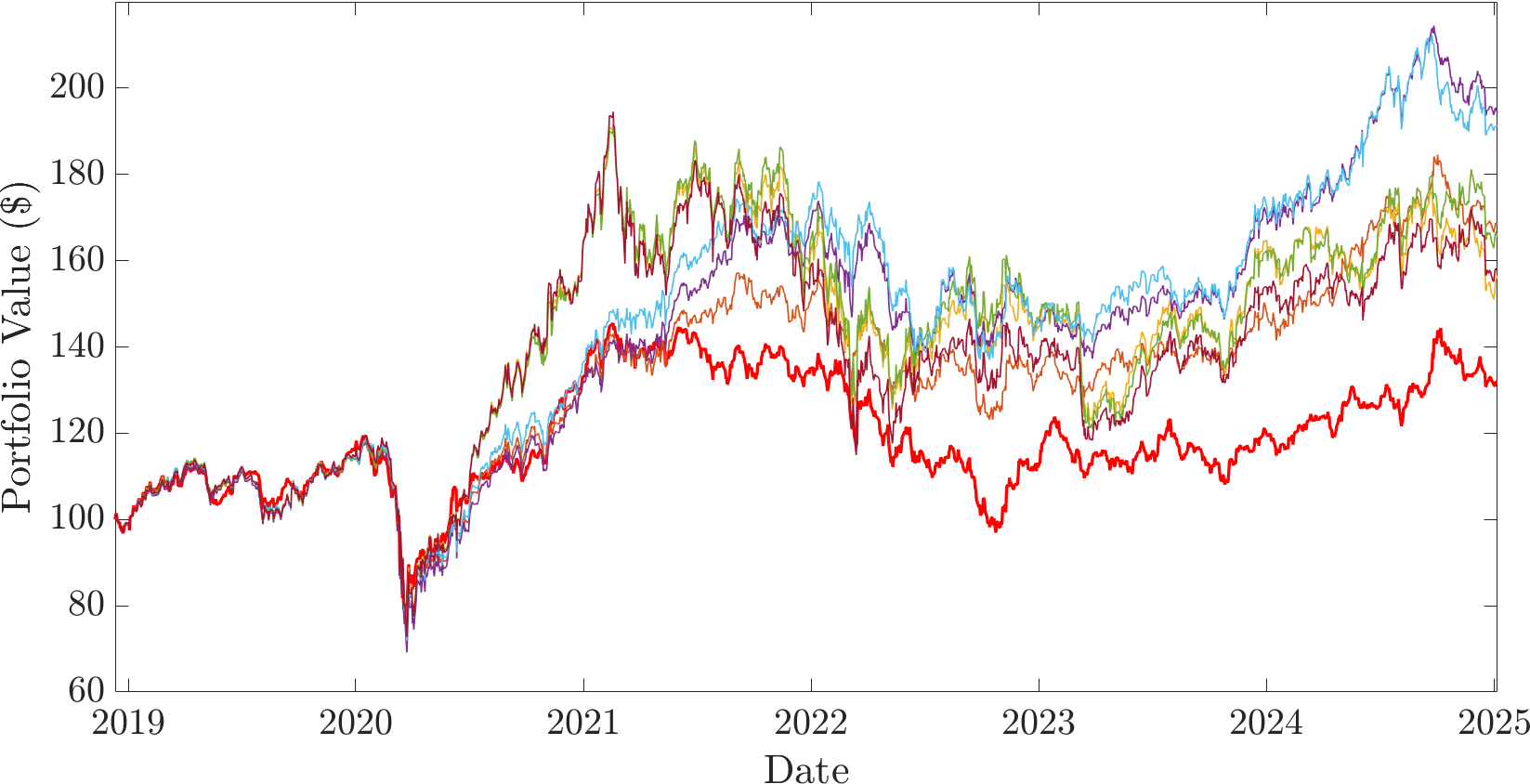}
\end{subfigure}

\vspace{0.1cm}

\begin{subfigure}{0.45\textwidth}
    \includegraphics[width=\linewidth]{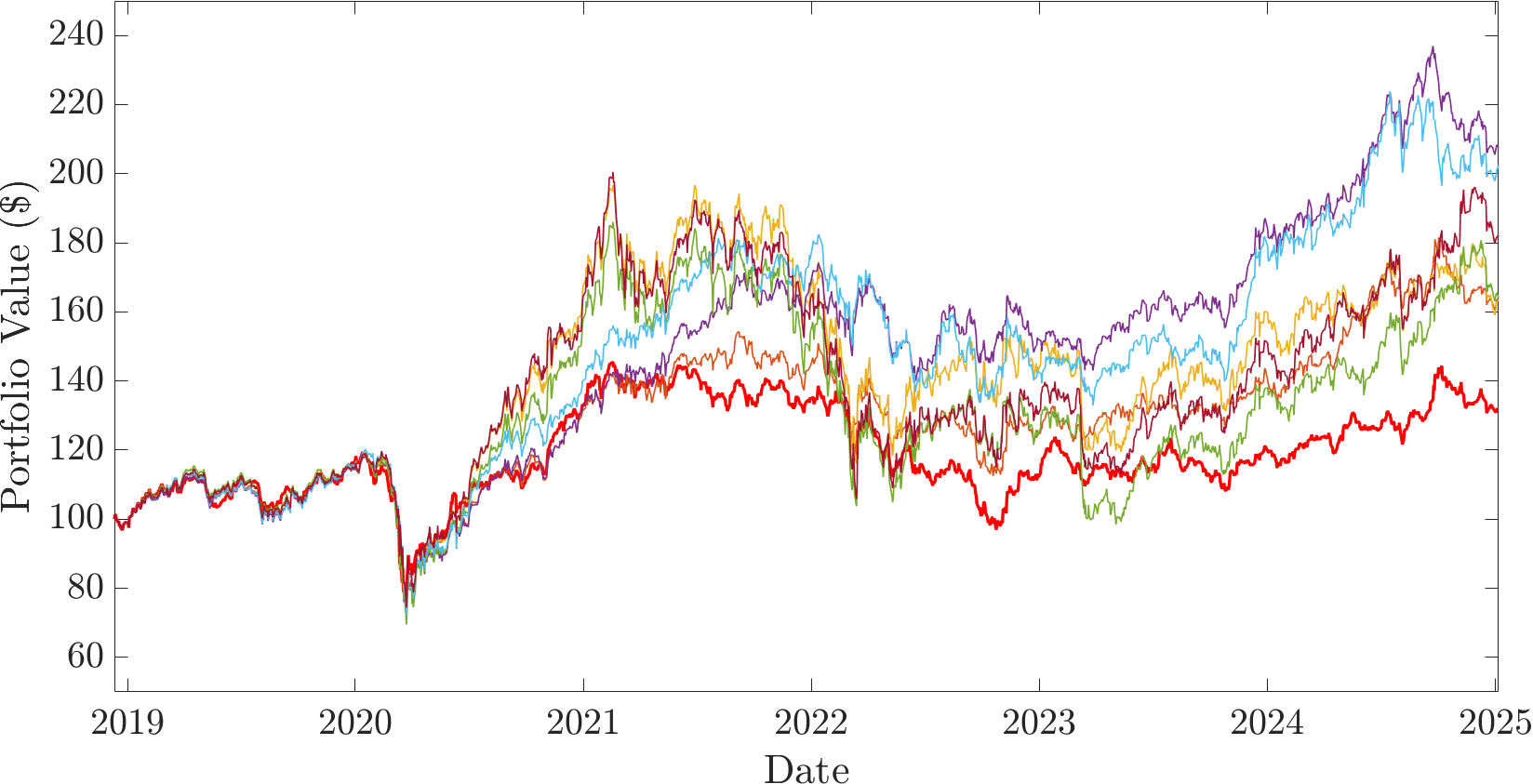}
\end{subfigure}
\hspace{0.5cm}
\begin{subfigure}{0.45\textwidth}
    \includegraphics[width=\linewidth]{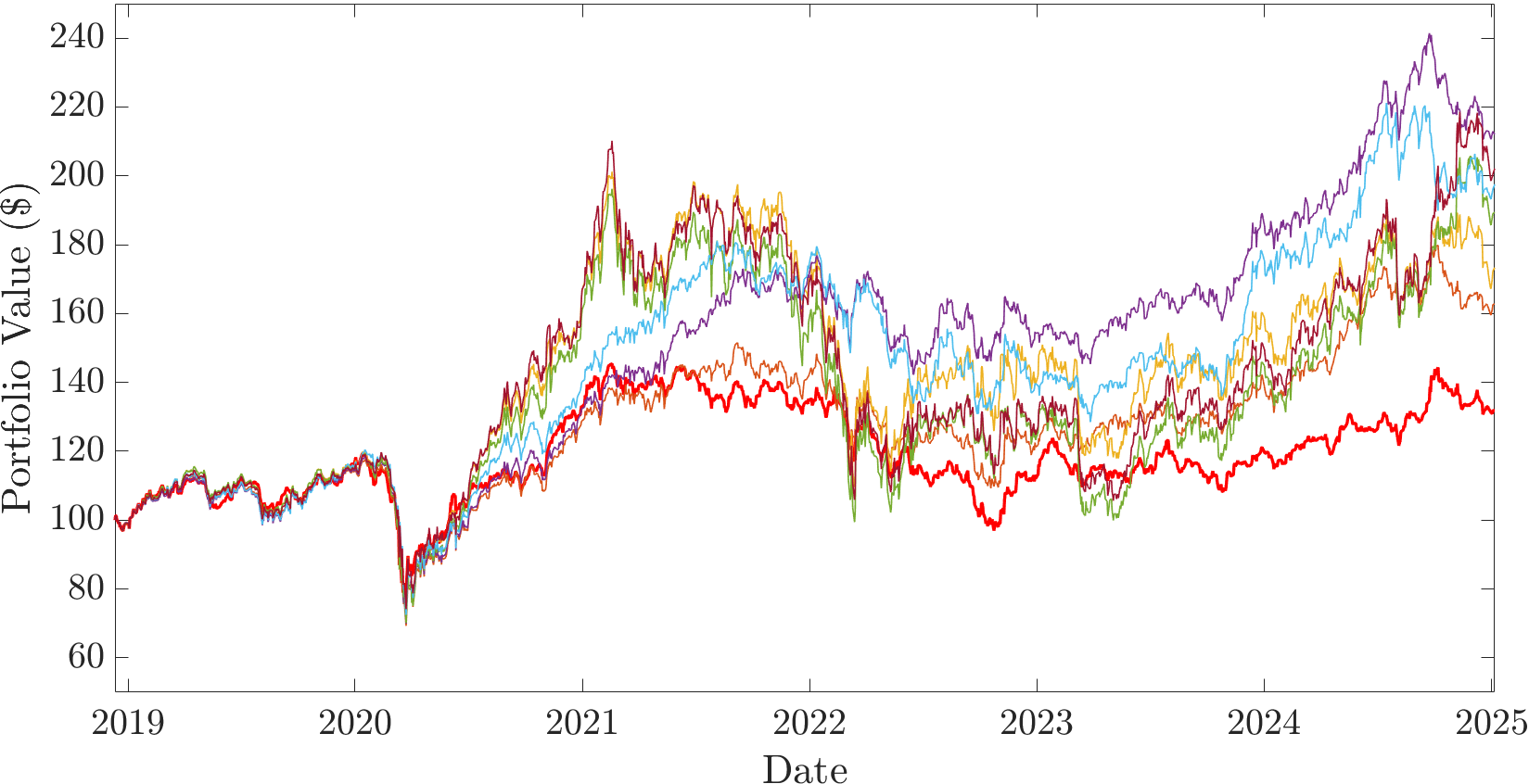}
\end{subfigure}

\caption{Cumulative price performance of the optimized portfolios under long-only and long–short strategies from December 2018 to January 2025. 
The top panels display the long-only strategy (left) and the long–short strategy with 10\% leverage (right). 
The bottom panels correspond to long–short strategies with leverage levels of 20\% (left) and 30\% (right). 
Each plot shows the growth of an initial \$100 investment on December 11, 2018 over the analysis period.}
\label{fig:lo_ls_cumprice}
\end{figure}

Across all strategies, 
the optimized portfolios consistently outperform the EQW, 
underscoring the value of systematic optimization in enhancing portfolio efficiency. 
This performance advantage is evident under both the long-only and long–short strategies, 
with the former delivering steady, 
risk-controlled growth and the latter achieving higher cumulative returns through the inclusion of short positions. 
\cite{lindquist2021advanced} note that the EQW portfolio is inherently long-only and typically excluded from long–short analyses, 
but it is retained here as a consistent benchmark to illustrate relative performance across all strategies. 
The long-only strategy exhibits smoother trajectories with moderate gains, 
whereas the introduction of leverage in the long–short strategy amplifies both returns and volatility. 
As leverage increases from 10\% to 30\%, 
portfolios experience progressively greater dispersion, 
demonstrating that while moderate leverage (10–20\%) enhances performance without destabilizing the portfolio, 
higher leverage (30\%) introduces pronounced fluctuations and sharper draw-downs during adverse market conditions.

In the long-only strategy, 
M95 achieves the strongest terminal value, 
with T95, T99, TVP, and M99 clustered closely behind. 
Introducing short exposure and moderate leverage (10–20\%) lifts the performance of return-seeking portfolios (T95, T99, TVP), 
though at the cost of greater volatility. 
At higher leverage (30\%), 
dispersion widens further as risk-focused portfolios (M95, M99) maintain smoother growth and better downside resilience. 
Overall, 
leverage amplifies both returns and risks, 
with moderate levels enhancing performance efficiency while excessive leverage magnifies fluctuations and draw-downs.

\section{Risk Return Performance Evaluation}

Performance evaluation provides insight into how efficiently a portfolio balances risk and return. 
While cumulative growth offers an overview of performance, 
risk adjusted measures allow for a more robust comparison across portfolios with varying risk exposures. 
The following ratios, Sharpe, stable tail adjusted return ratio (STARR), and Rachev, 
capture different aspects of portfolio efficiency and tail risk sensitivity. 
Together, 
they provide a comprehensive framework for understanding portfolio robustness under both normal and extreme market conditions 
\citep{goetzmann2002sharpening, lo2002statistics, gatfaoui2009sharpe, steiner2011sharpe, smetters2013sharper, dafonseca2020performance, choi2015reward, zhang2025risk}.

In this section, the performance of the optimized portfolios (MVP, TVP, M95, T95, M99, and T99) is evaluated using these three ratios under both the long only and long short strategies with leverage levels of 10\%, 20\%, and 30\%. 
The EQW portfolio is also included as a benchmark to assess the relative improvement achieved through optimization.

For a broader perspective, 
Appendix B presents box-plots of the Sharpe, STARR, and Rachev ratios for the 29 individual Asian ETFs included in our study. 
This highlights how each ETF performs on a risk adjusted basis, 
providing deeper insight into their individual behavior compared to their collective performance within the EQW and optimized portfolios. 
The results offer a more detailed view of asset level performance, 
complementing the portfolio level findings discussed in the main analysis.

\subsection{Sharpe Ratio}

The Sharpe ratio \citep{sharpe1999investments} measures the excess return earned per unit of total risk. 
It is one of the most widely used performance measures in finance, 
providing a straightforward way to evaluate the trade-off between risk and return.

There have been several studies on the Sharpe ratio. 
\cite{lo2002statistics} highlights that Sharpe ratios can be overstated due to serial correlation and time aggregation effects, 
meaning that monthly Sharpe ratios cannot simply be annualized by multiplying by $\sqrt{12}$. 
\cite{goetzmann2002sharpening} demonstrate that Sharpe ratios can be manipulated through option-like payoffs, 
raising caution when comparing managers who employ derivatives. 
\cite{gatfaoui2009sharpe} introduces a filtering method to extract ``fundamental'' Sharpe ratios that are robust to skewness and kurtosis biases, 
reflecting purer performance indicators.
The Sharpe ratio is formally defined as follows

\begin{equation}
SR(T) = \frac{\mathbb{E}[r_p(t) - r_f(t)]}{\sqrt{\mathrm{VAR}[r_p(t) - r_f(t)]}}
\end{equation}

A higher Sharpe ratio indicates that the portfolio achieves greater excess returns for a given level of volatility. 
A negative Sharpe ratio occurs when the expected excess return $\mathbb{E}(r_p - r_f)$ is negative, 
meaning the portfolio under-performs the risk free rate on average, 
while a positive Sharpe ratio reflects positive excess performance. 
However, 
because the Sharpe ratio relies on the assumption of normally distributed returns, 
it may not fully capture performance under extreme market conditions.

\begin{figure}[H]
\centering

\begin{subfigure}{0.45\textwidth}
    \includegraphics[width=\linewidth]{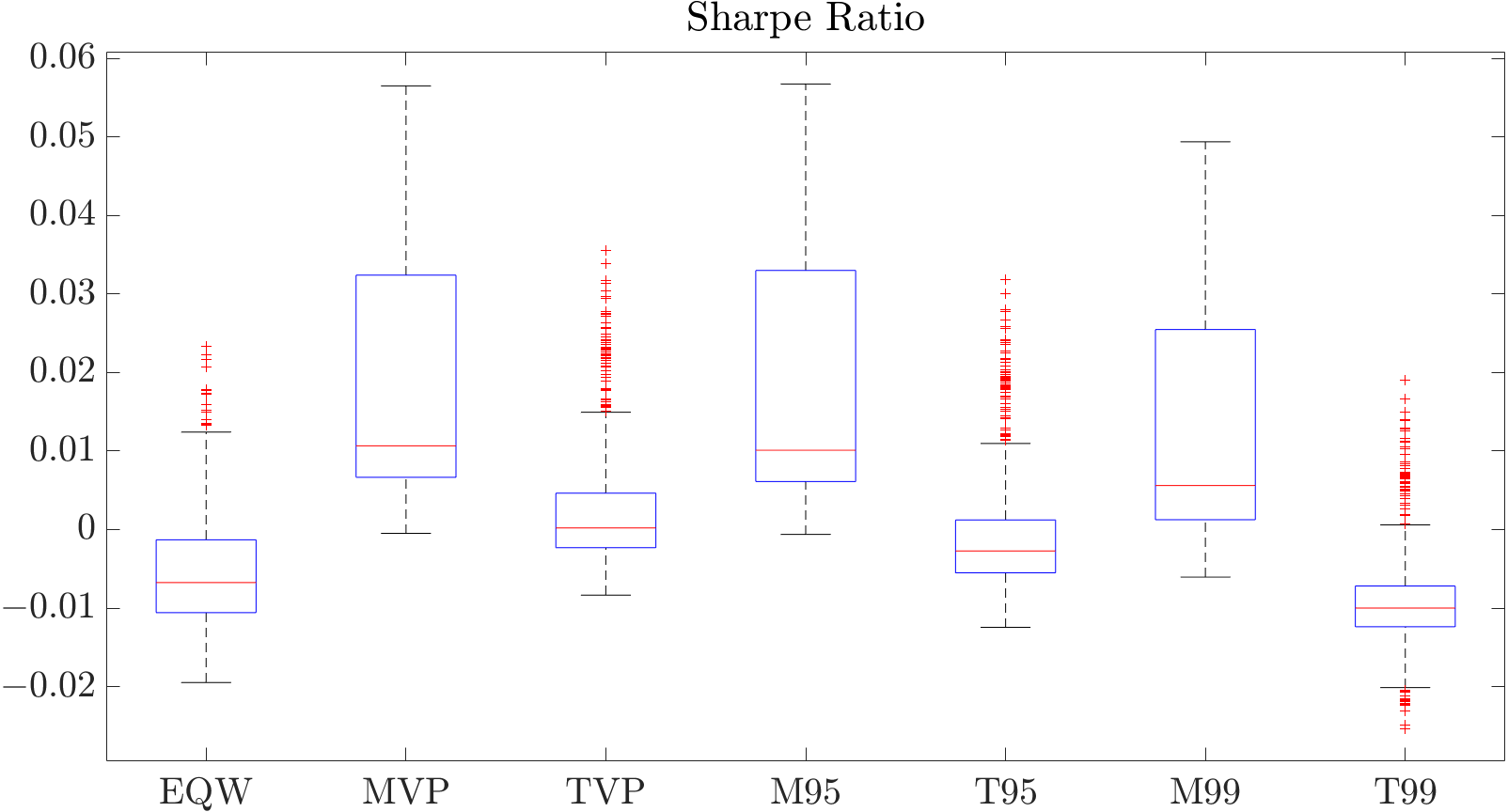}
\end{subfigure}
\hspace{0.5cm}
\begin{subfigure}{0.45\textwidth}
    \includegraphics[width=\linewidth]{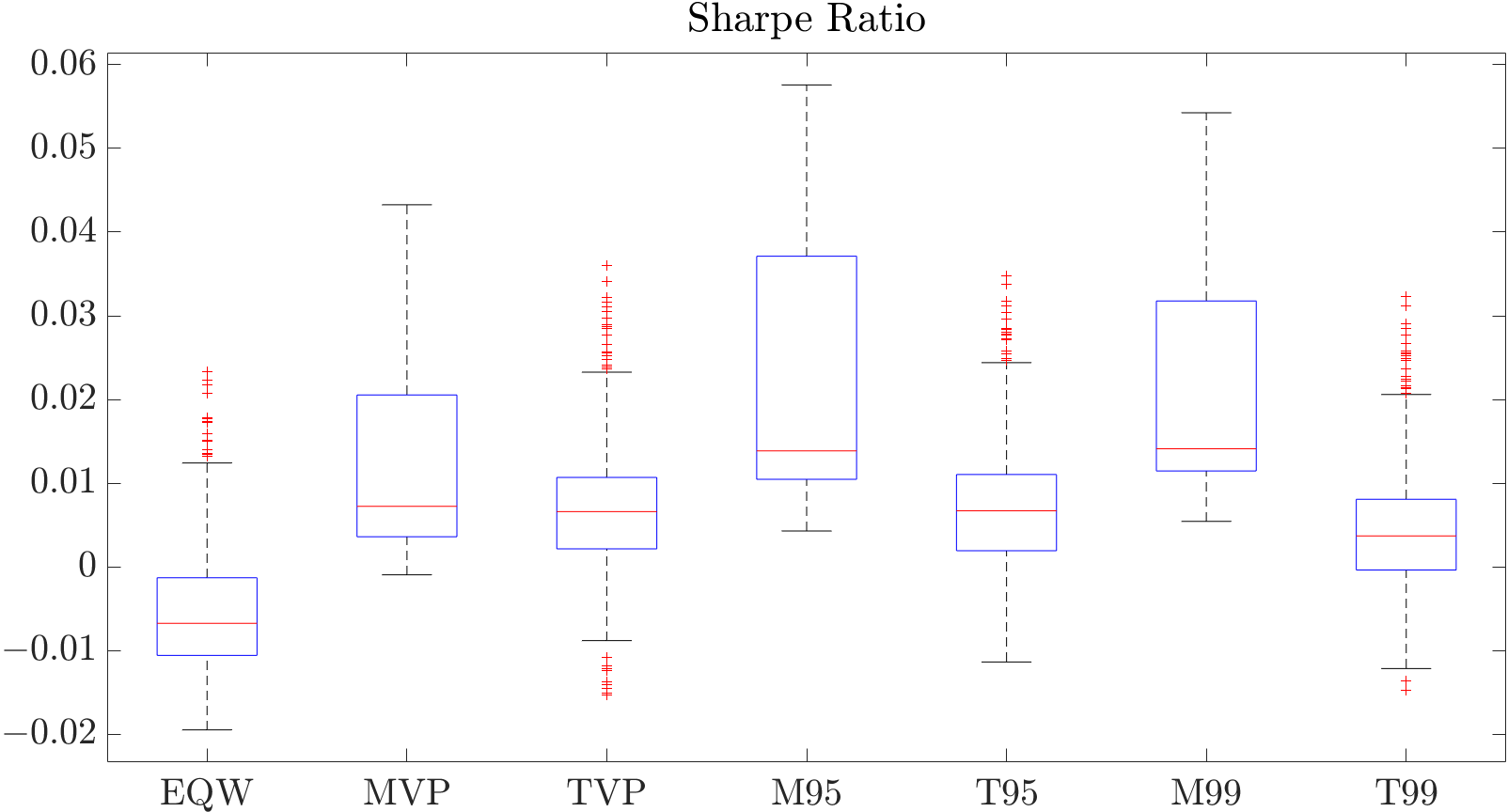}
\end{subfigure}

\vspace{0.1cm}

\begin{subfigure}{0.45\textwidth}
    \includegraphics[width=\linewidth]{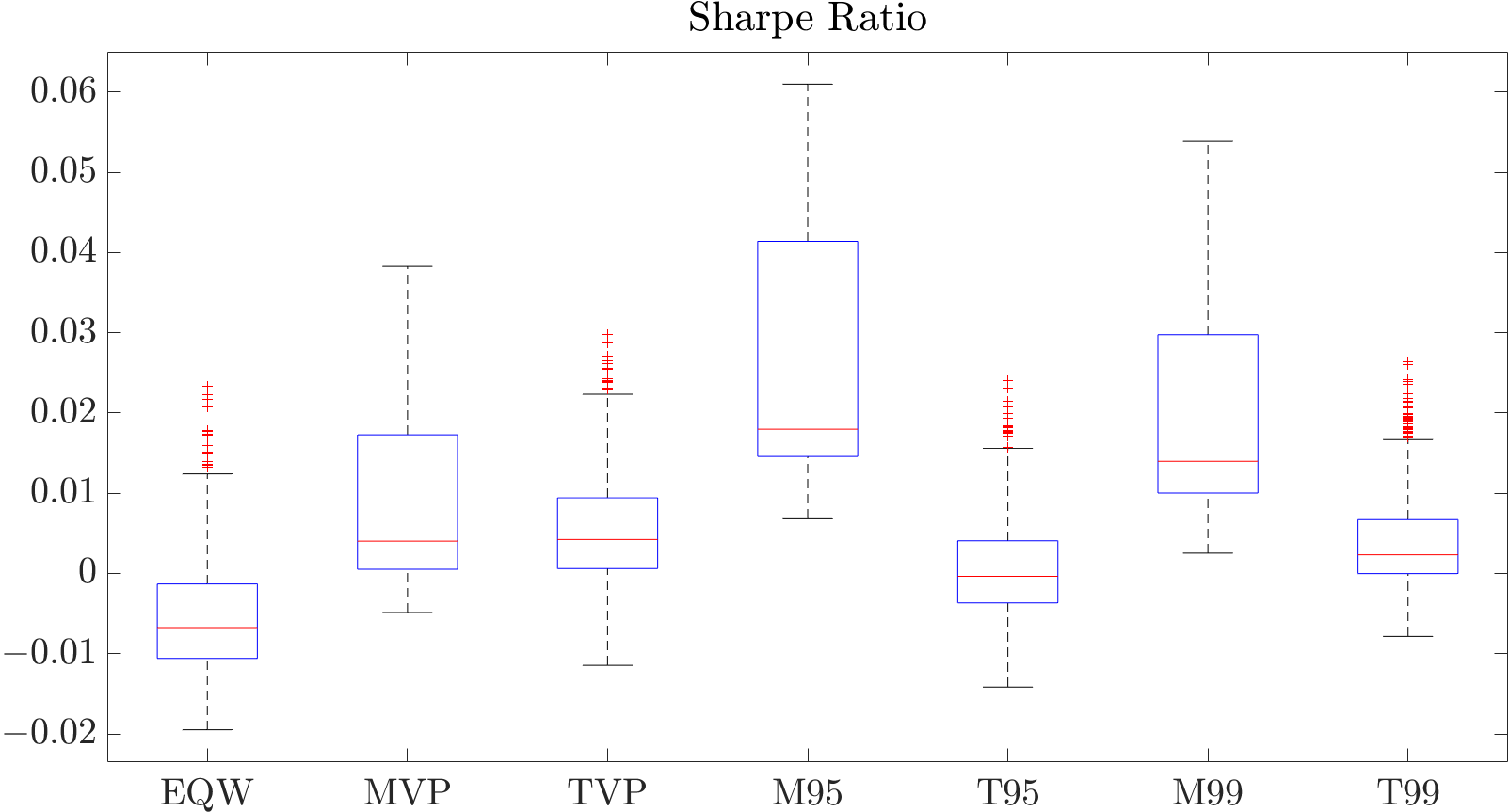}
\end{subfigure}
\hspace{0.5cm}
\begin{subfigure}{0.45\textwidth}
    \includegraphics[width=\linewidth]{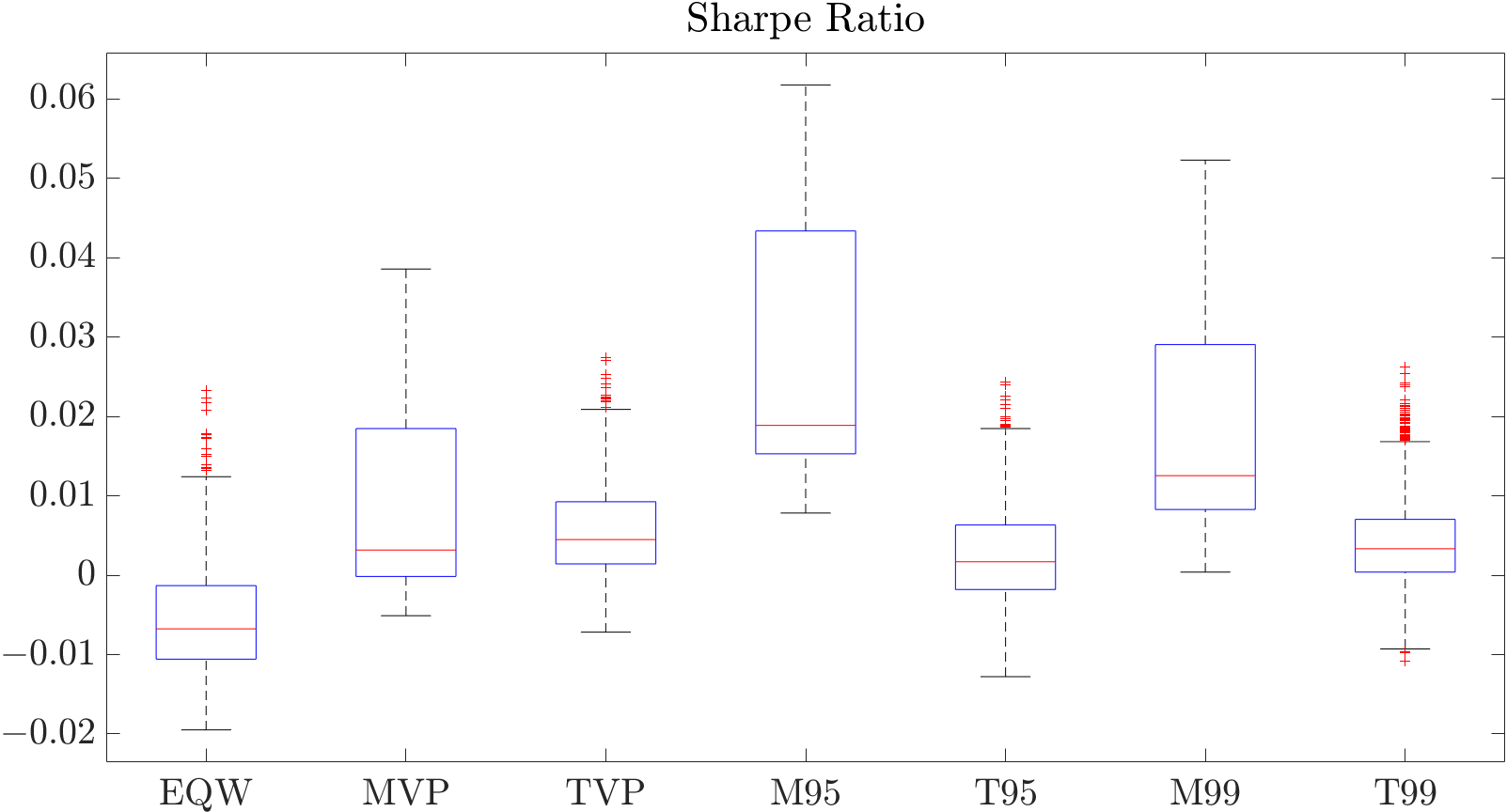}
\end{subfigure}
\caption{Sharpe ratios of the optimized portfolios and EQW under long-only and long–short strategies. 
Each figure displays results under four configurations: 
Long–Only (top-left), Long–Short 10\% (top-right), Long–Short 20\% (bottom-left), and Long–Short 30\% (bottom-right).}
\label{fig:sharpe_lo_ls}
\end{figure}

Figure \ref{fig:sharpe_lo_ls} illustrates the Sharpe ratios of the optimized portfolios under long-only and long-short strategies. 
Across all configurations, the MVP, M95, 
and M99 portfolios demonstrate superior performance relative to the EQW, TVP, T95, and T99 portfolios. 
Their higher median Sharpe ratios and narrower inter-quartile ranges indicate stronger and more consistent risk-adjusted returns. 
In contrast, 
the equally weighted and threshold-based portfolios show lower central tendencies and greater variability, 
suggesting less efficient reward-to-risk trade-offs. 
The performance gains of the optimized portfolios are especially pronounced under long-short strategies, 
although increased leverage (20\% and 30\%) amplifies return dispersion, 
reflecting the heightened sensitivity of portfolio outcomes to market fluctuations.

\subsection{Rachev Ratio}
The Rachev ratio (RR) \citep{rachev2008advanced} compares average gains in large positive market swings to average losses in large negative swings, 
making it particularly useful for understanding portfolio behavior under extreme market conditions. 
It focuses on asymmetry between the tails of the return distribution rather than overall volatility, 
providing insight into how well a strategy performs in favorable versus adverse market environments.

\citet{biglova2004different} introduced this ratio to emphasize the balance between expected tail gains and expected tail losses, 
and subsequent studies have applied it to stress testing and portfolio construction frameworks. 
\citet{dafonseca2020performance} also includes the RR among key performance measures when evaluating international portfolios, 
finding that it captures reward–risk efficiency better during periods of market stress.
The RR is mathematically defined as

\begin{equation}
RR_{\alpha,\beta}(T) = 
\frac{\mathrm{CVaR}_{\beta 
}(r_f(t) - r_p(t))}
{\mathrm{CVaR}_{\alpha}(r_p(t) - r_f(t))}
\end{equation}

In this paper, RR is evaluated at both the 95\% and 99\% confidence levels to analyze 
reward–risk efficiency under different degrees of tail sensitivity.
The lower this ratio is, 
the higher is the probability of extreme gains relative to extreme losses.

\begin{figure}[H]
\centering

\begin{subfigure}{0.45\textwidth}
    \includegraphics[width=\linewidth]{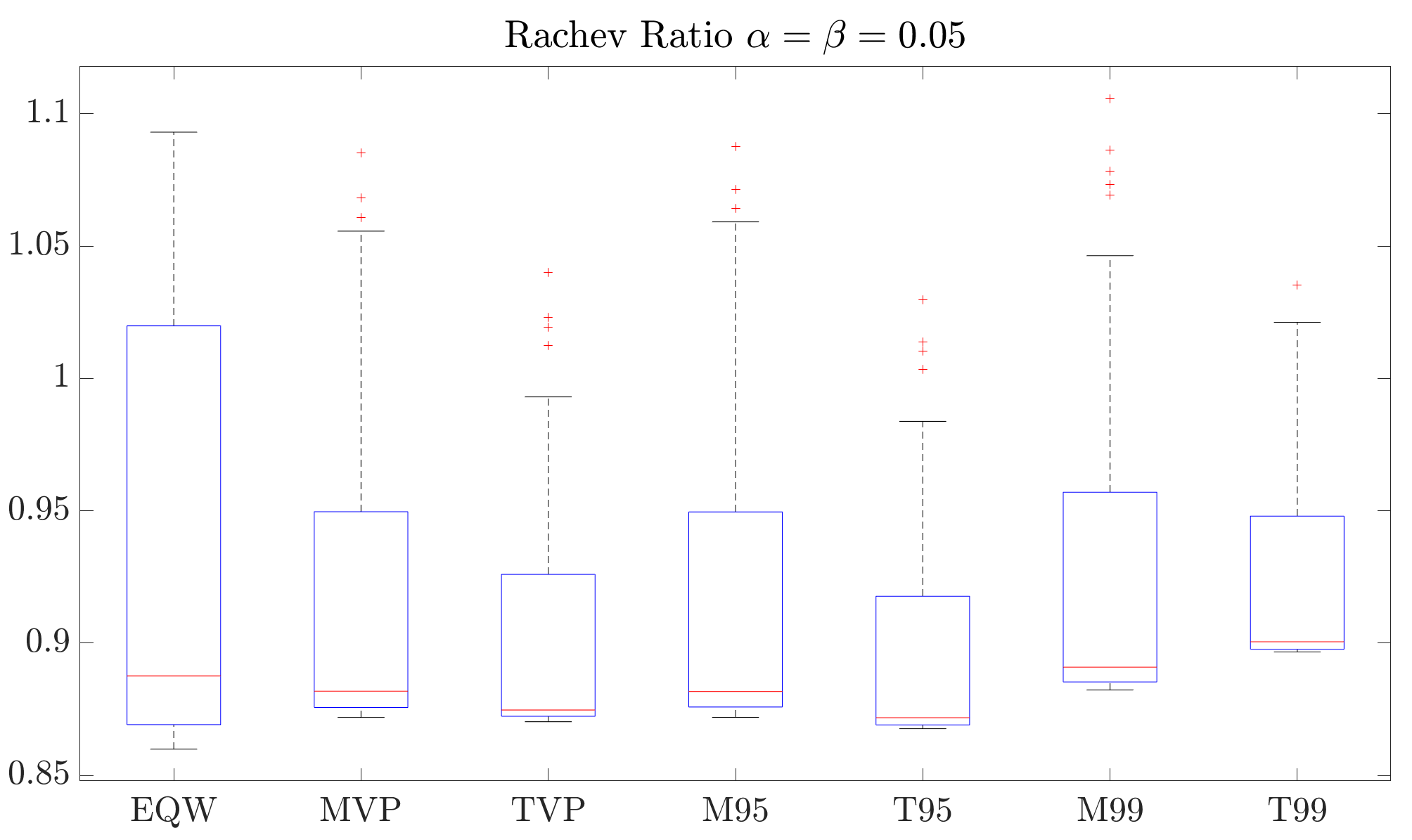}
\end{subfigure}
\hspace{0.5cm}
\begin{subfigure}{0.45\textwidth}
    \includegraphics[width=\linewidth]{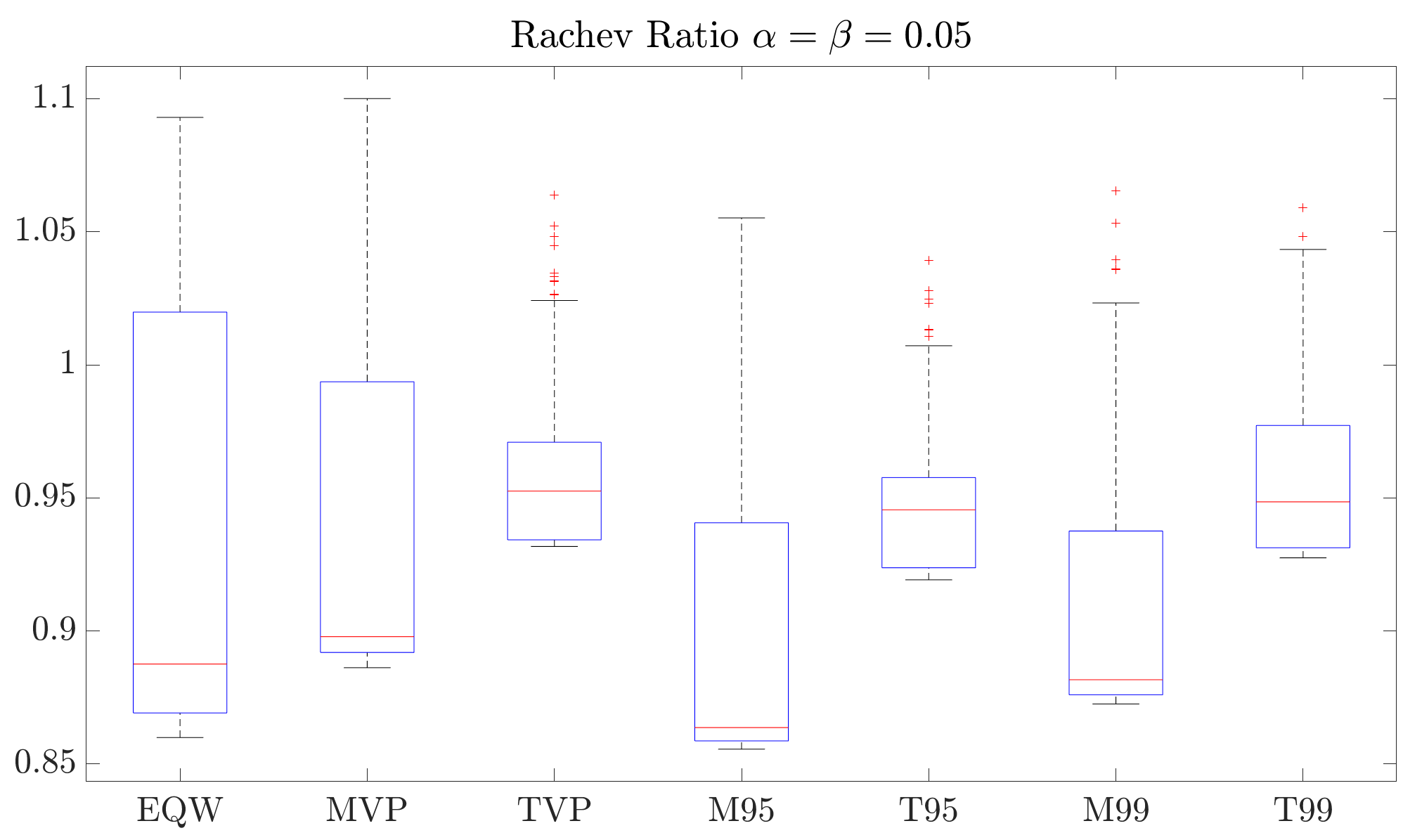}
\end{subfigure}

\vspace{0.1cm}

\begin{subfigure}{0.45\textwidth}
    \includegraphics[width=\linewidth]{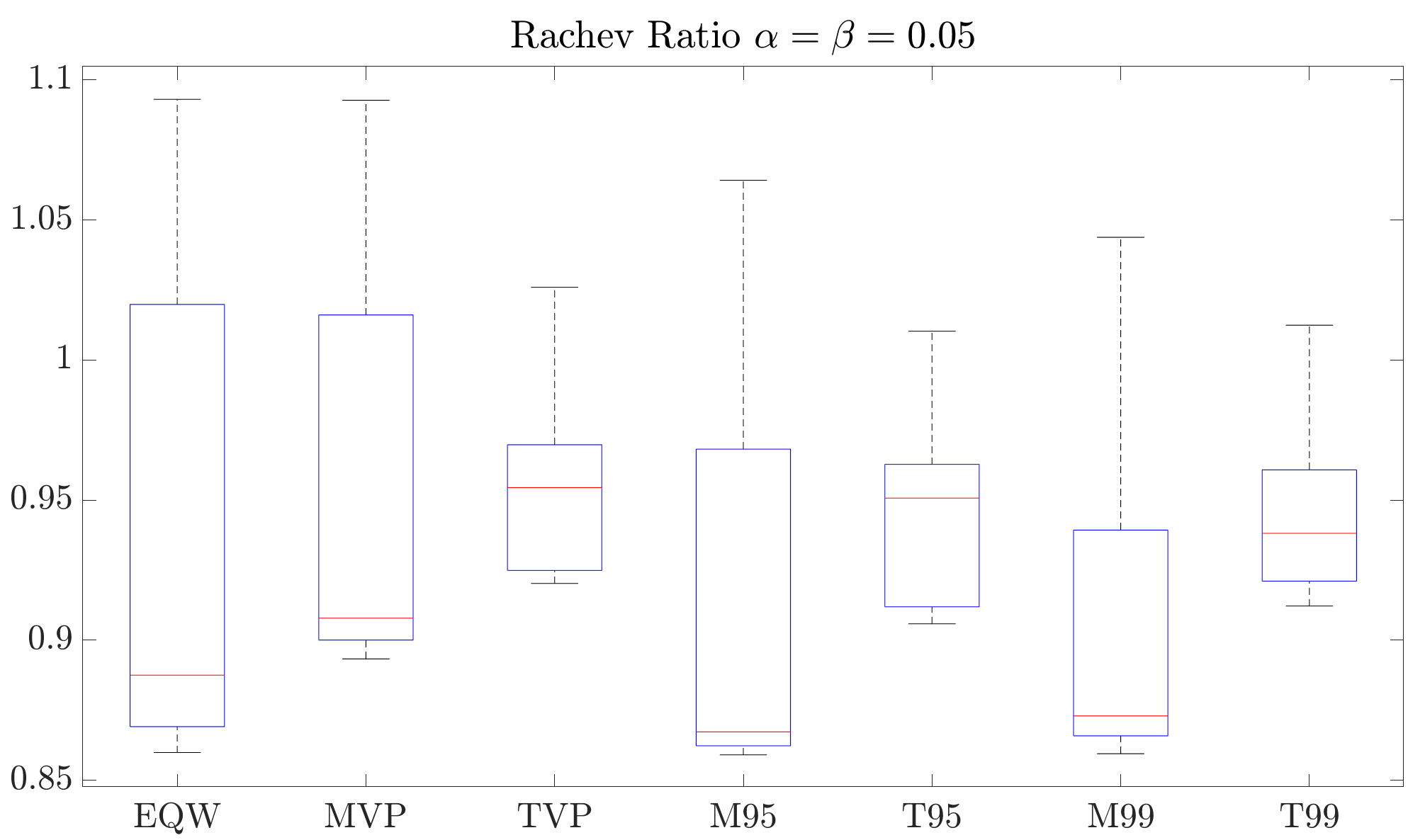}
\end{subfigure}
\hspace{0.5cm}
\begin{subfigure}{0.45\textwidth}
    \includegraphics[width=\linewidth]{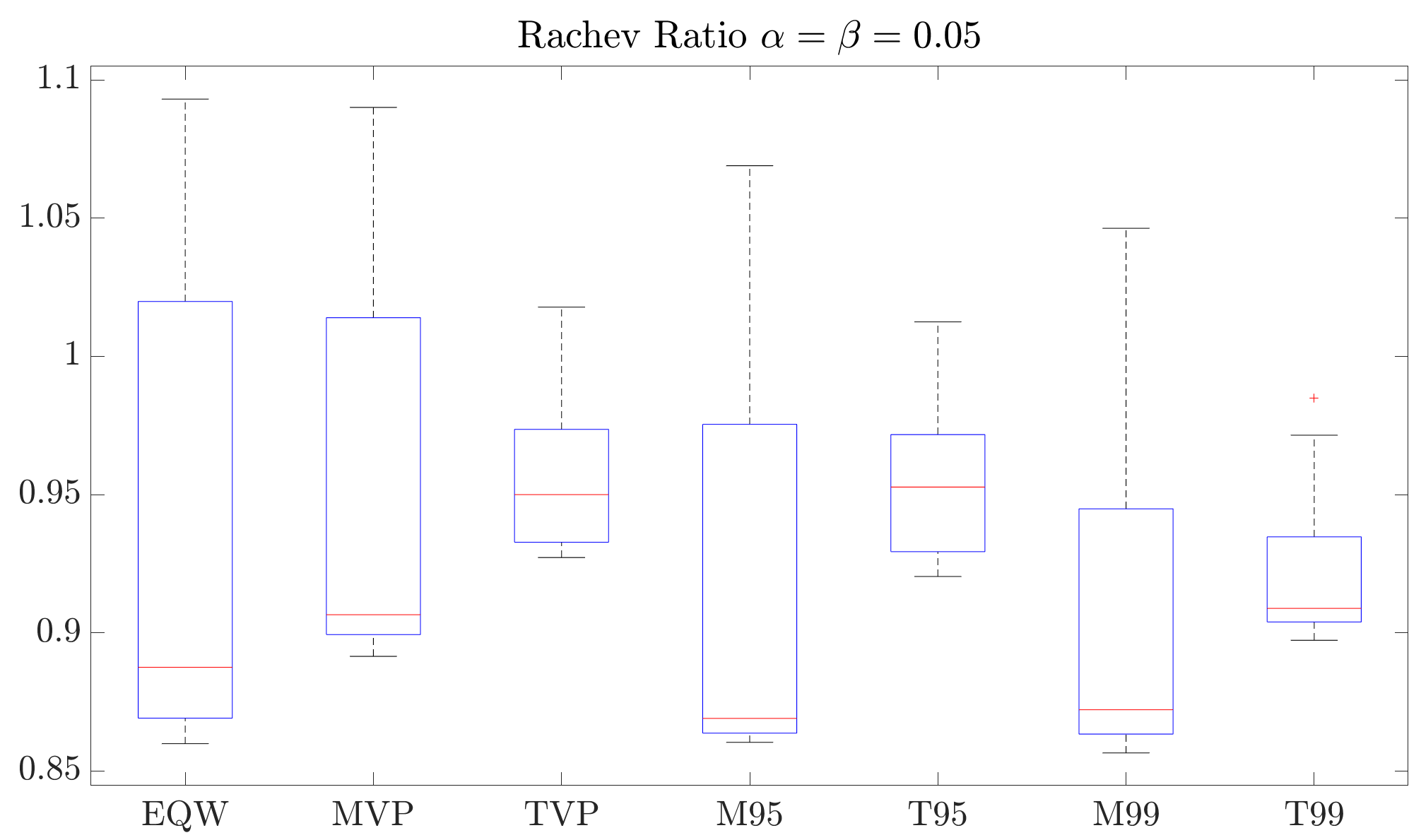}
\end{subfigure}

\caption{RR at the 95\% confidence level for the optimized portfolios and EQW under long-only and long–short strategies. Each figure displays results under four configurations: Long–Only (top-left), Long–Short 10\% (top-right), Long–Short 20\% (bottom-left), and Long–Short 30\% (bottom-right).}
\label{fig:rr95_lo_ls}
\end{figure}

\begin{figure}[H]
\centering

\begin{subfigure}{0.45\textwidth}
    \includegraphics[width=\linewidth]{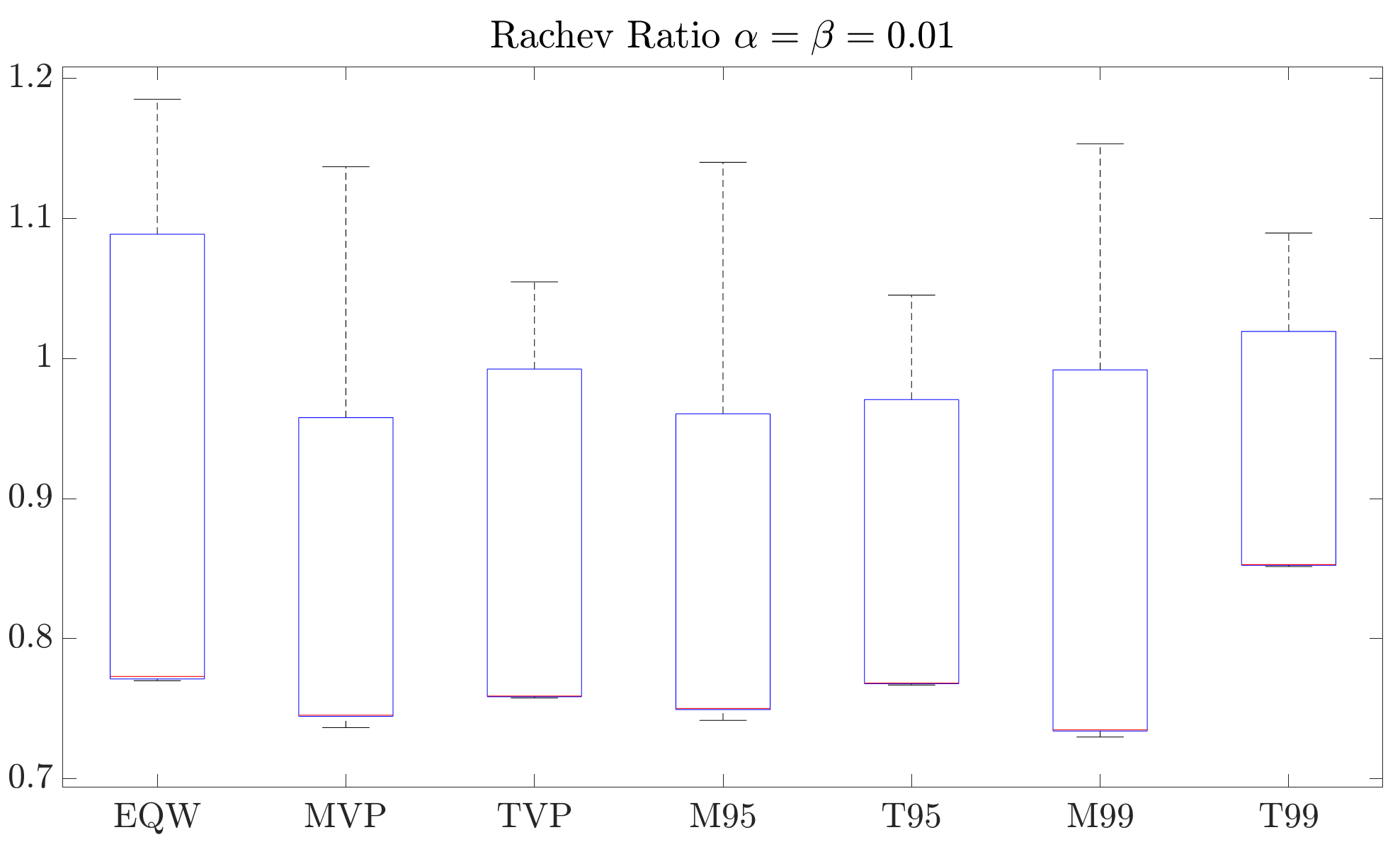}
\end{subfigure}
\hspace{0.5cm}
\begin{subfigure}{0.45\textwidth}
    \includegraphics[width=\linewidth]{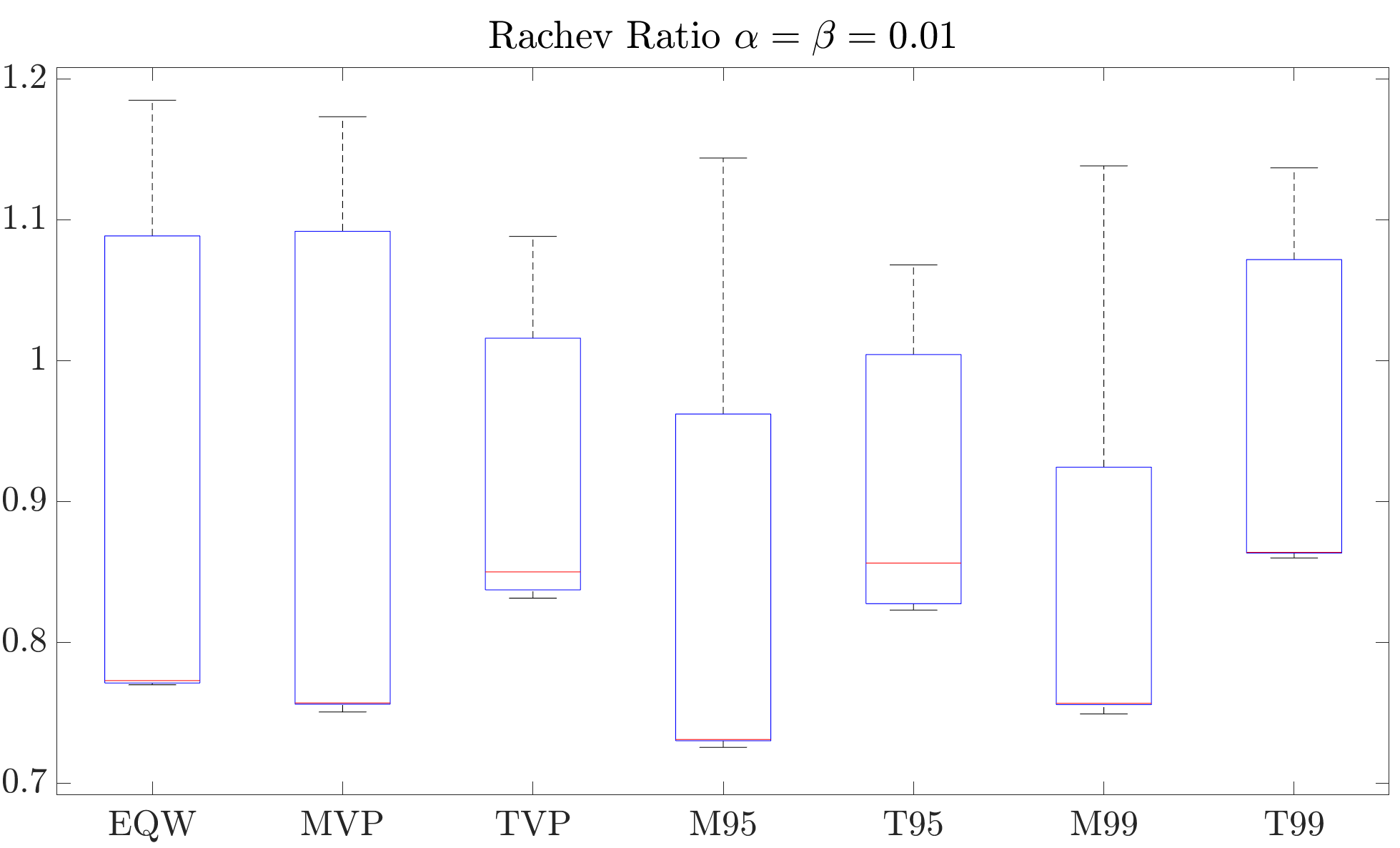}
\end{subfigure}

\vspace{0.1cm}

\begin{subfigure}{0.45\textwidth}
    \includegraphics[width=\linewidth]{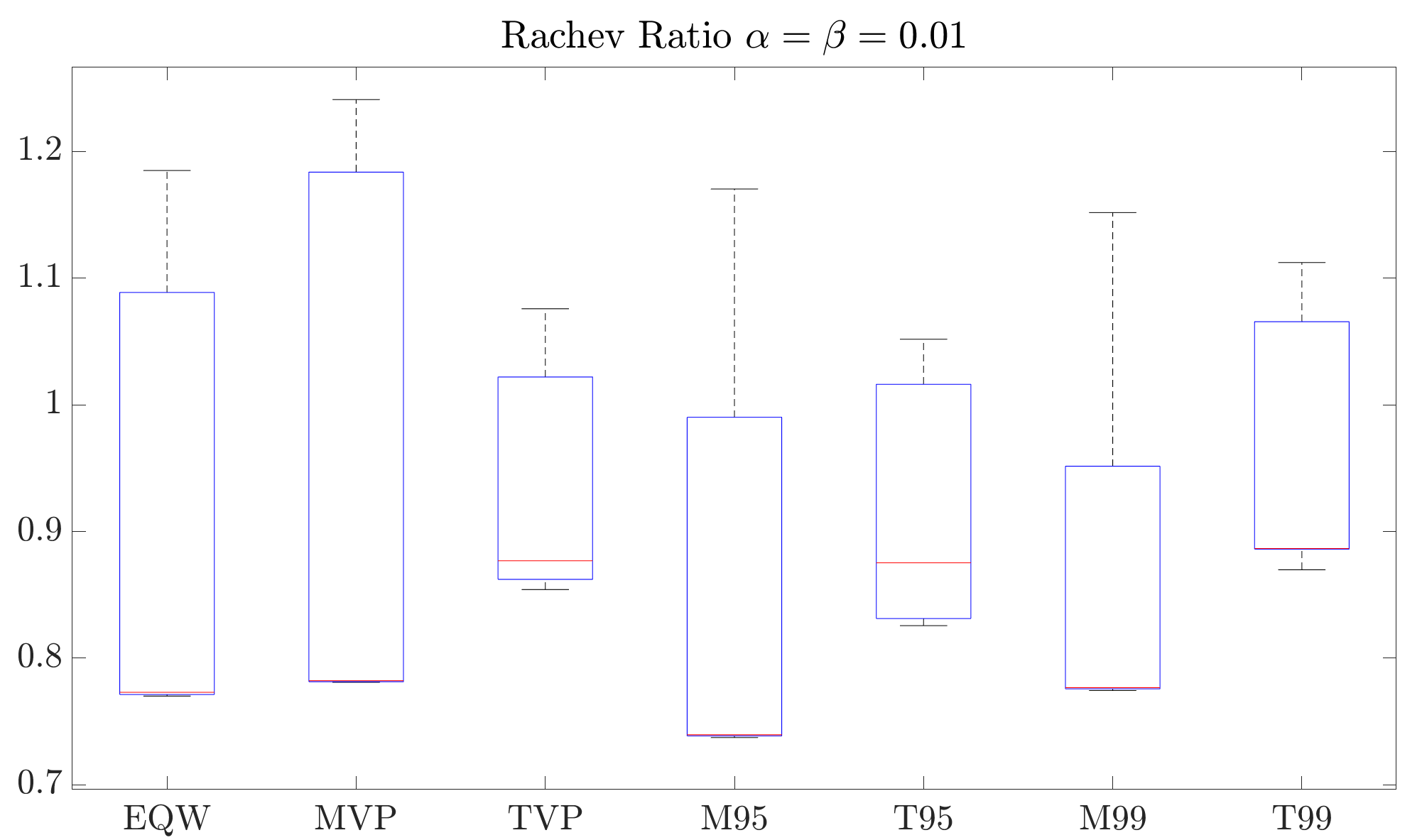}
\end{subfigure}
\hspace{0.5cm}
\begin{subfigure}{0.45\textwidth}
    \includegraphics[width=\linewidth]{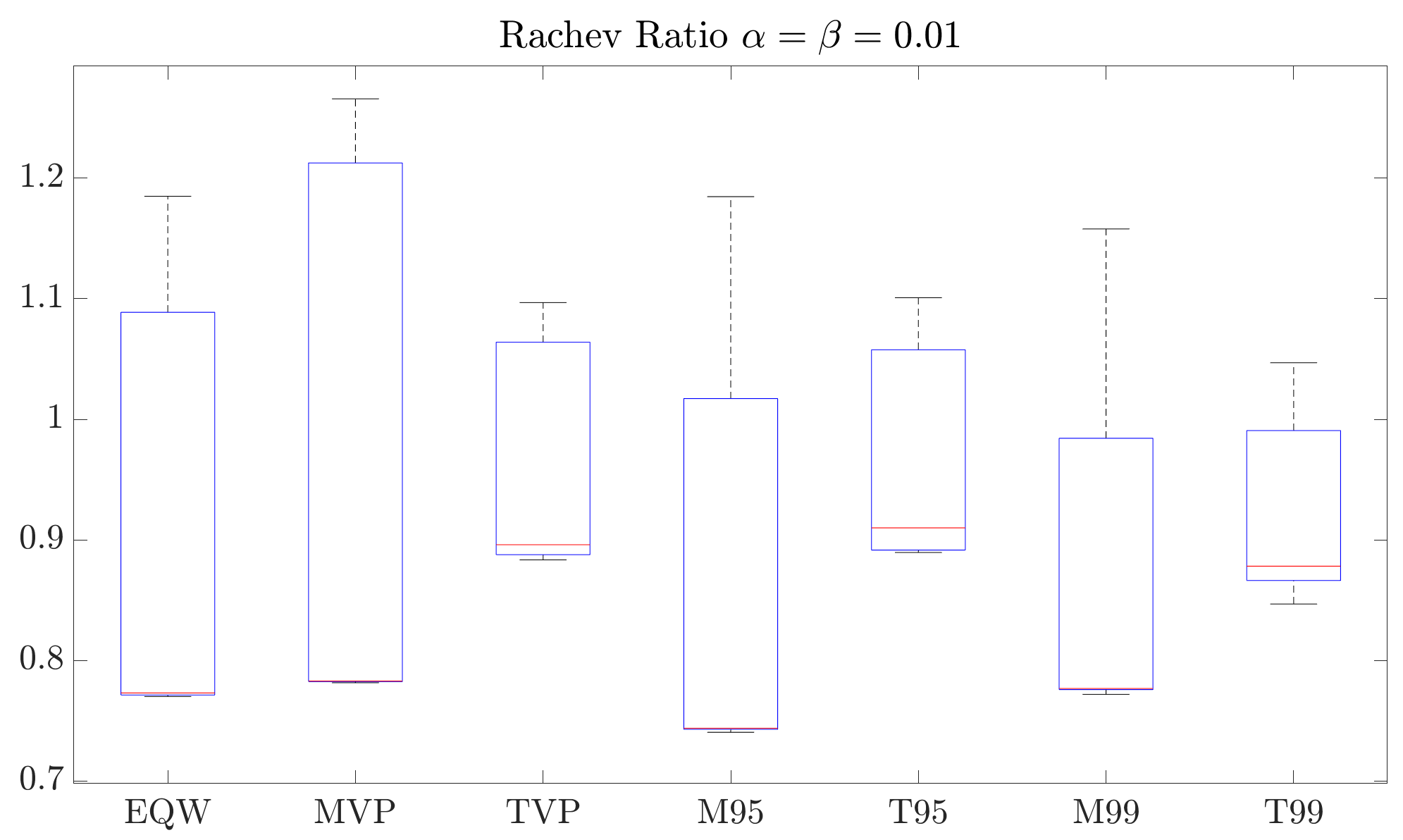}
\end{subfigure}

\caption{RR at the 99\% confidence level for the optimized portfolios and EQW under long-only and long–short strategies. 
Each figure displays results under four configurations: 
Long–Only (top-left), Long–Short 10\% (top-right), Long–Short 20\% (bottom-left), and Long–Short 30\% (bottom-right).}
\label{fig:rr99_lo_ls}
\end{figure}

While the Sharpe ratio (Figure \ref{fig:sharpe_lo_ls}) captures mean–variance efficiency, 
the RR in Figures \ref{fig:rr95_lo_ls} and \ref{fig:rr99_lo_ls} accounts for asymmetry between extreme gains and losses by measuring the ratio of expected returns in the right tail (favorable outcomes) to expected losses in the left tail (adverse outcomes). 
The box-plots display the distribution of these ratio values. 
Values of RR near one indicate roughly symmetric tail behavior, 
whereas higher values (greater than one) reflect stronger upside potential relative to downside risk.
The RR plots show greater dispersion than the Sharpe ratio, 
reflecting the variability in upside versus downside tail performance. 
At the 95\% confidence level (Figure \ref{fig:rr95_lo_ls}), 
long-only portfolios generally display higher and more stable RR values than their long-short counterparts, 
with MVP and TVP performing relatively well, 
suggesting that moderate diversification helps manage tail losses effectively. 
As leverage increases in the long-short portfolios, 
the RR distribution widens, indicating heightened exposure to tail risk. 
Moving to the 99\% confidence level (Figure \ref{fig:rr99_lo_ls}), 
the dispersion becomes even more pronounced, 
as the stricter tail threshold captures rarer and more extreme events. 
Here, the median RR values decline slightly across all strategies, 
showing that under extreme conditions, 
downside risk dominates more strongly. 
The difference between the 95\% and 99\% levels illustrates how sensitive the RR is to extreme losses, a
s the confidence level increases, 
tail performance becomes less stable and less favorable. 
Overall, while Sharpe ratios emphasize average risk-adjusted efficiency, 
the RR provides a more nuanced view of tail vulnerability, 
revealing that higher leverage and stricter confidence levels magnify downside asymmetry in returns.

\subsection{STARR Ratio}
The stable tail adjusted return ratio (STARR) builds on the Rachev ratio and the Sharpe ratio by focusing on reward per unit of expected tail loss. 
By replacing volatility with CVaR in its denominator, 
it captures downside sensitivity more effectively and offers a conservative measure of performance \citep{rachev2008advanced}.

\citet{martin2003portfolios} first introduced the STAR ratio, 
which later evolved into the STARR measure, 
and research such as \citet{sehgal2021robust} and \citet{zhang2025risk} shows that incorporating CVaR improves the robustness of portfolio optimization, particularly during turbulent markets. 
\citet{zhang2025risk} further demonstrates that STARR based momentum strategies achieve stronger downside risk control and sector-level consistency, confirming its practical relevance in heterogeneous market conditions.
The STARR ratio is mathematically expressed as

\begin{equation}
STARR(T) = \frac{\mathbb{E}[r_p(t) - r_f(t)]}{\mathrm{CVaR}_{\alpha}[r_p(t) - r_f(t)]}
\end{equation}

By focusing on CVaR, the STARR ratio captures how efficiently a portfolio compensates investors for 
exposure to extreme downside risk, offering a more conservative assessment of performance.

\begin{figure}[H]
\centering

\begin{subfigure}{0.45\textwidth}
    \includegraphics[width=\linewidth]{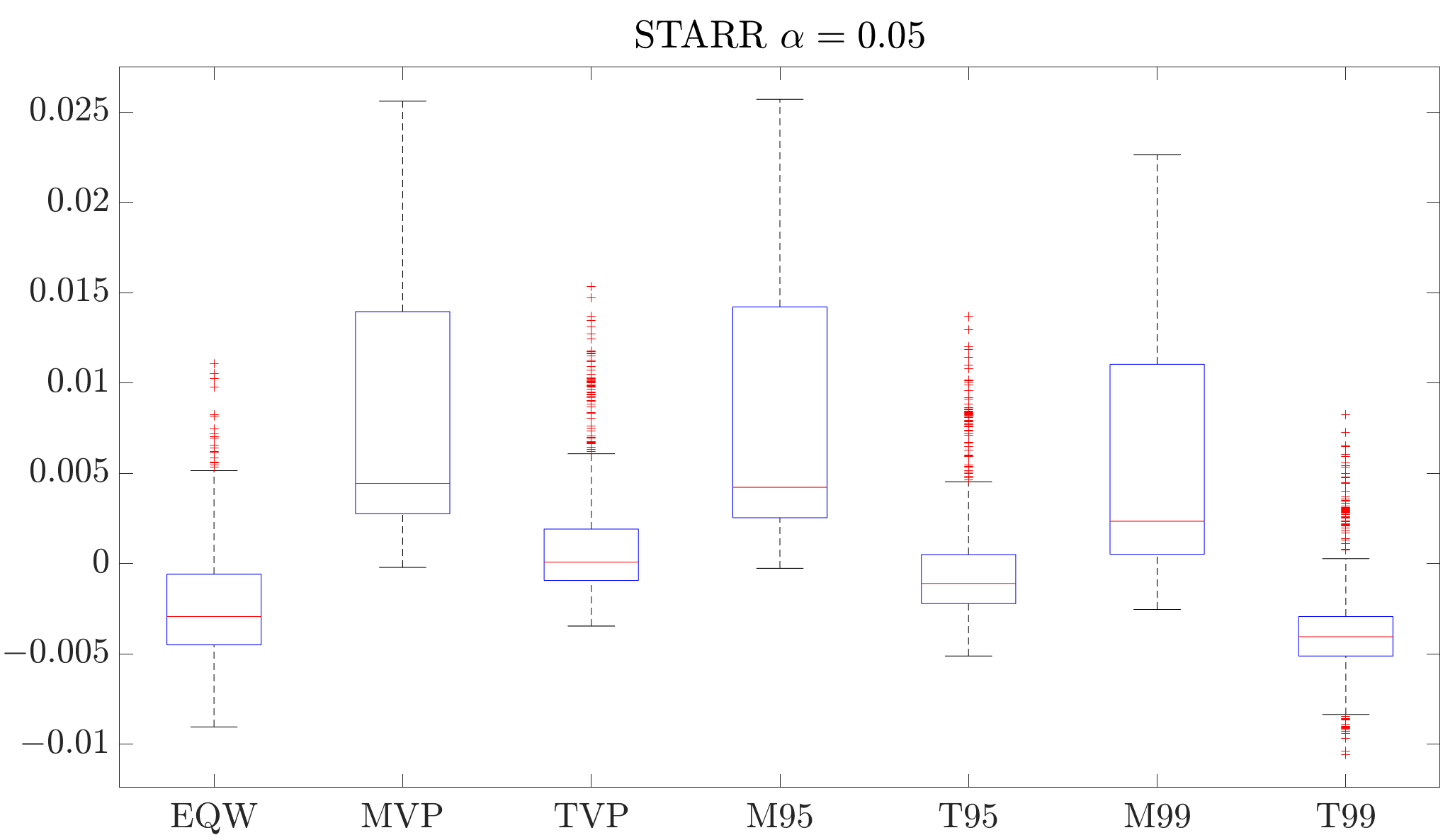}
\end{subfigure}
\hspace{0.5cm}
\begin{subfigure}{0.45\textwidth}
    \includegraphics[width=\linewidth]{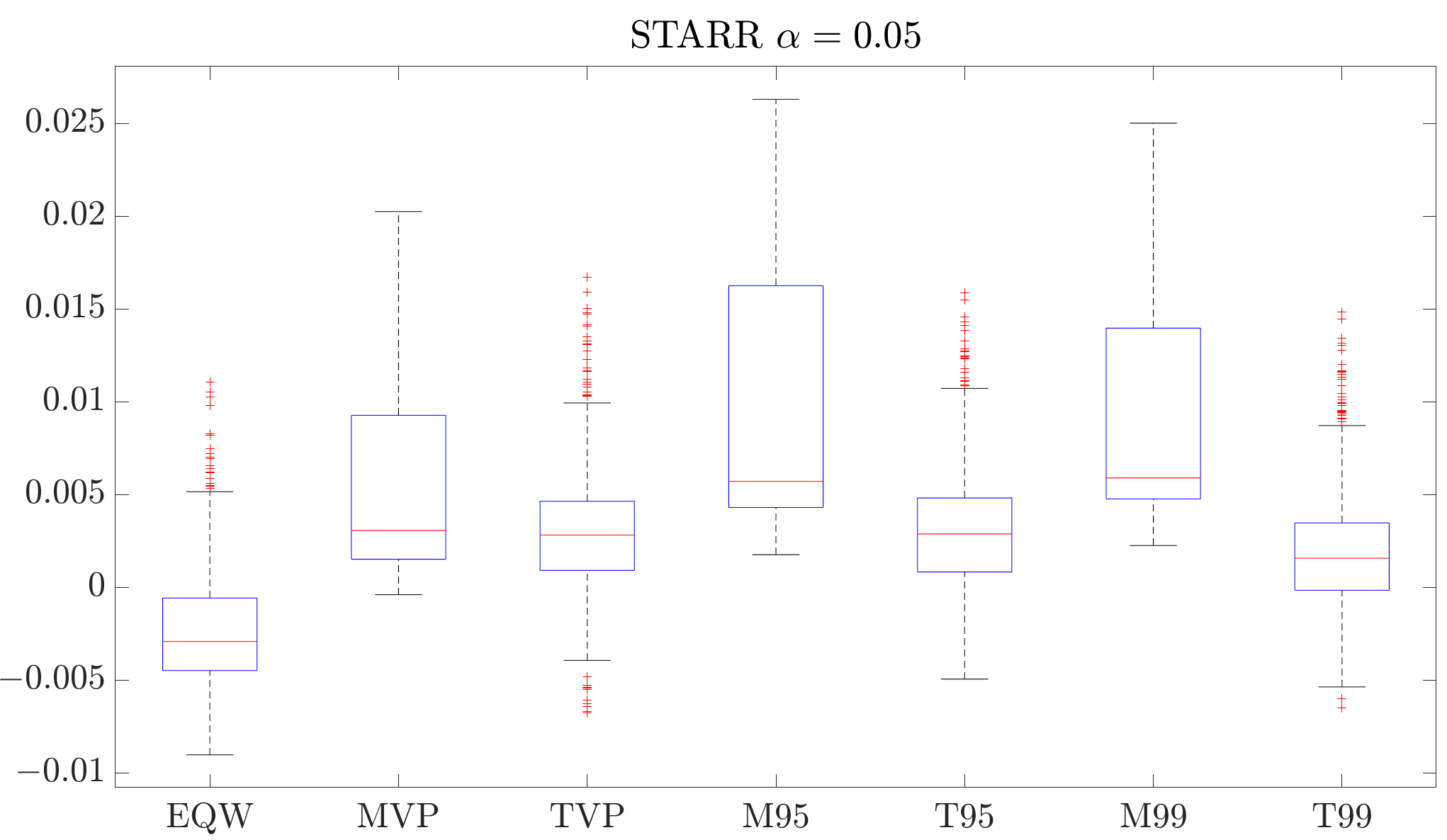}
\end{subfigure}

\vspace{0.1cm}

\begin{subfigure}{0.45\textwidth}
    \includegraphics[width=\linewidth]{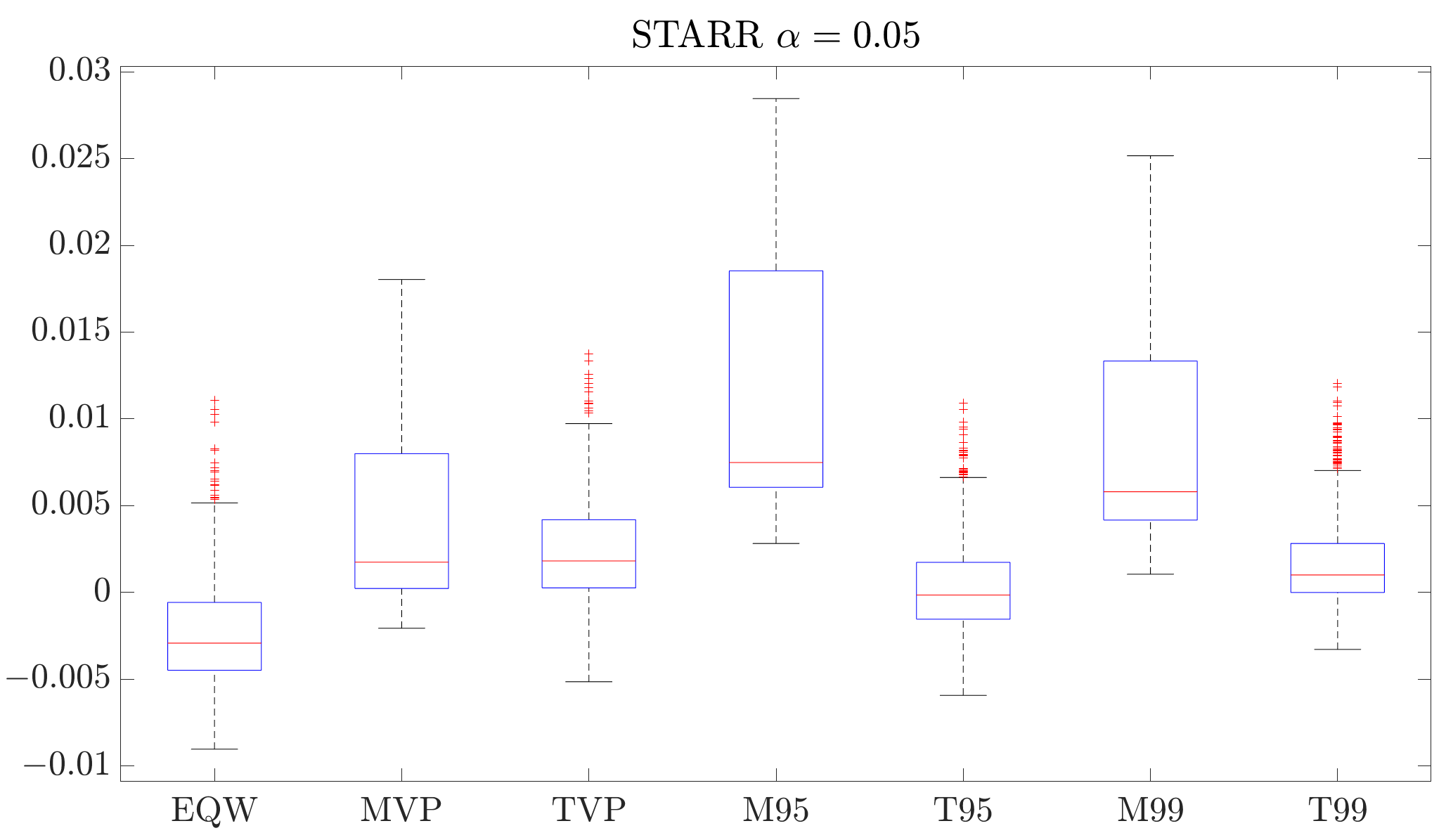}
\end{subfigure}
\hspace{0.5cm}
\begin{subfigure}{0.45\textwidth}
    \includegraphics[width=\linewidth]{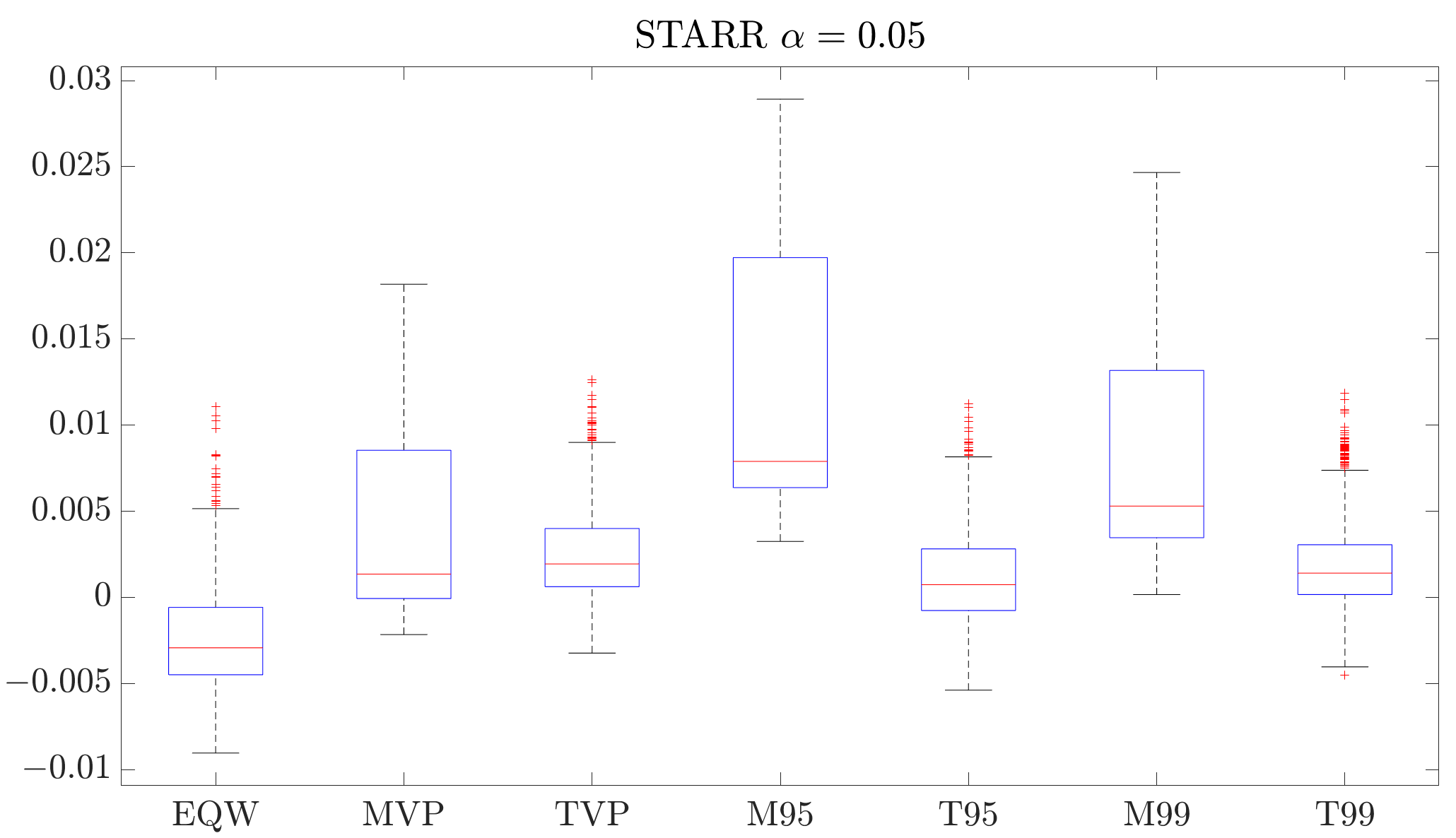}
\end{subfigure}

\caption{STARR at the 95\% confidence level for the optimized portfolios and EQW under long-only and long–short strategies. 
Each figure displays results under four configurations: 
Long–Only (top-left), Long–Short 10\% (top-right), Long–Short 20\% (bottom-left), and Long–Short 30\% (bottom-right).}
\label{fig:starr95_lo_ls}
\end{figure}

\begin{figure}[H]
\centering

\begin{subfigure}{0.45\textwidth}
    \includegraphics[width=\linewidth]{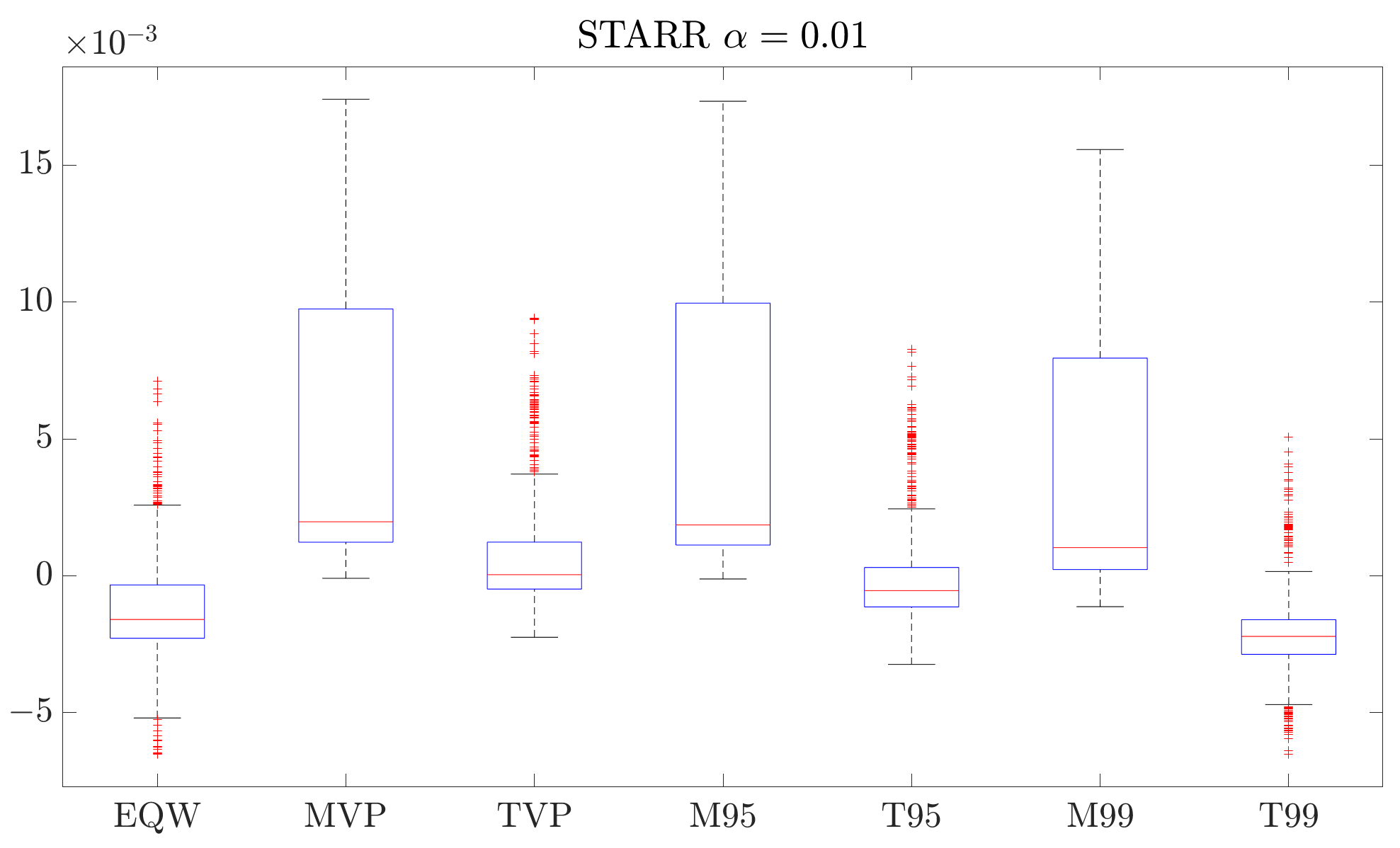}
\end{subfigure}
\hspace{0.5cm}
\begin{subfigure}{0.45\textwidth}
    \includegraphics[width=\linewidth]{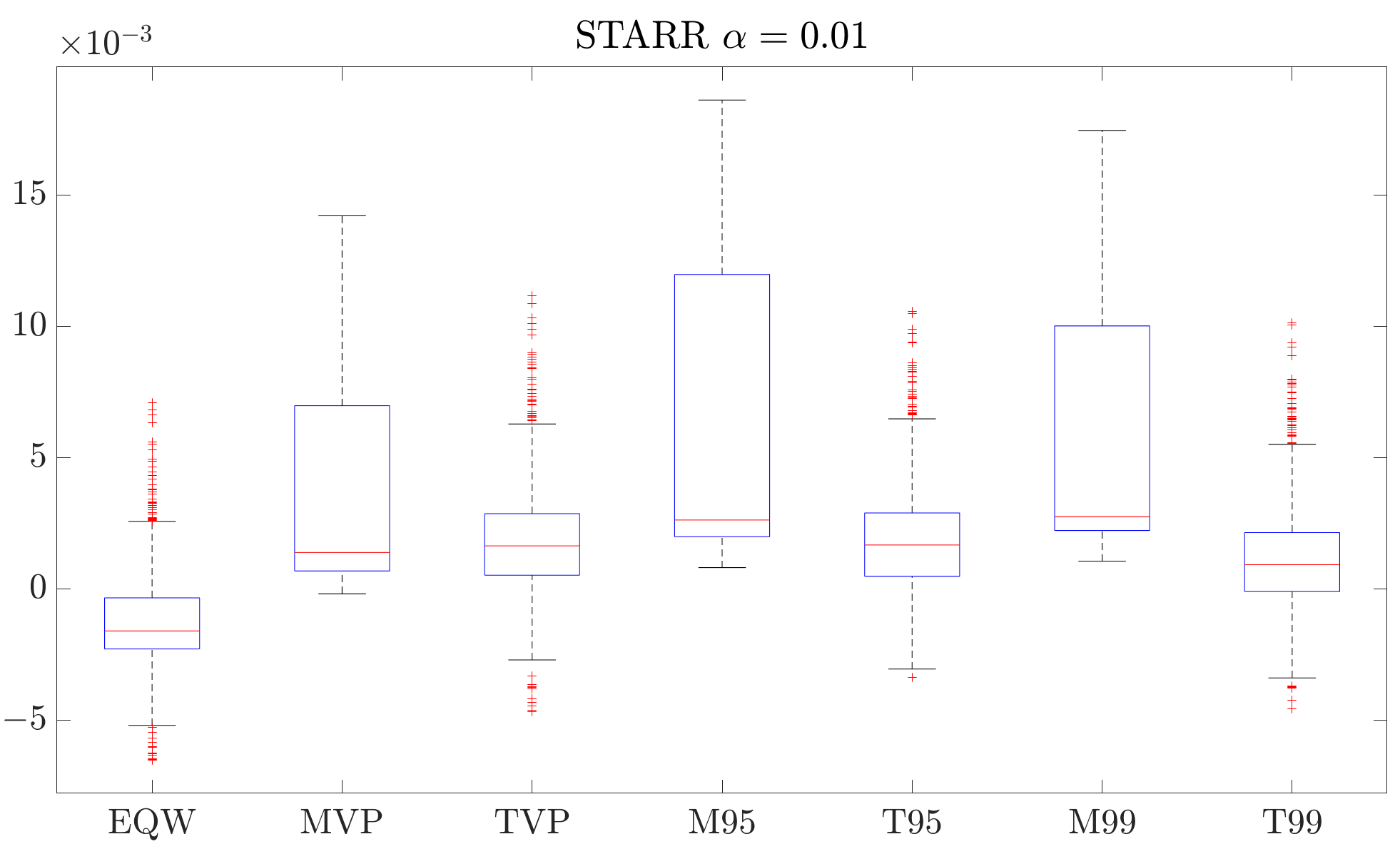}
\end{subfigure}

\vspace{0.1cm}

\begin{subfigure}{0.45\textwidth}
    \includegraphics[width=\linewidth]{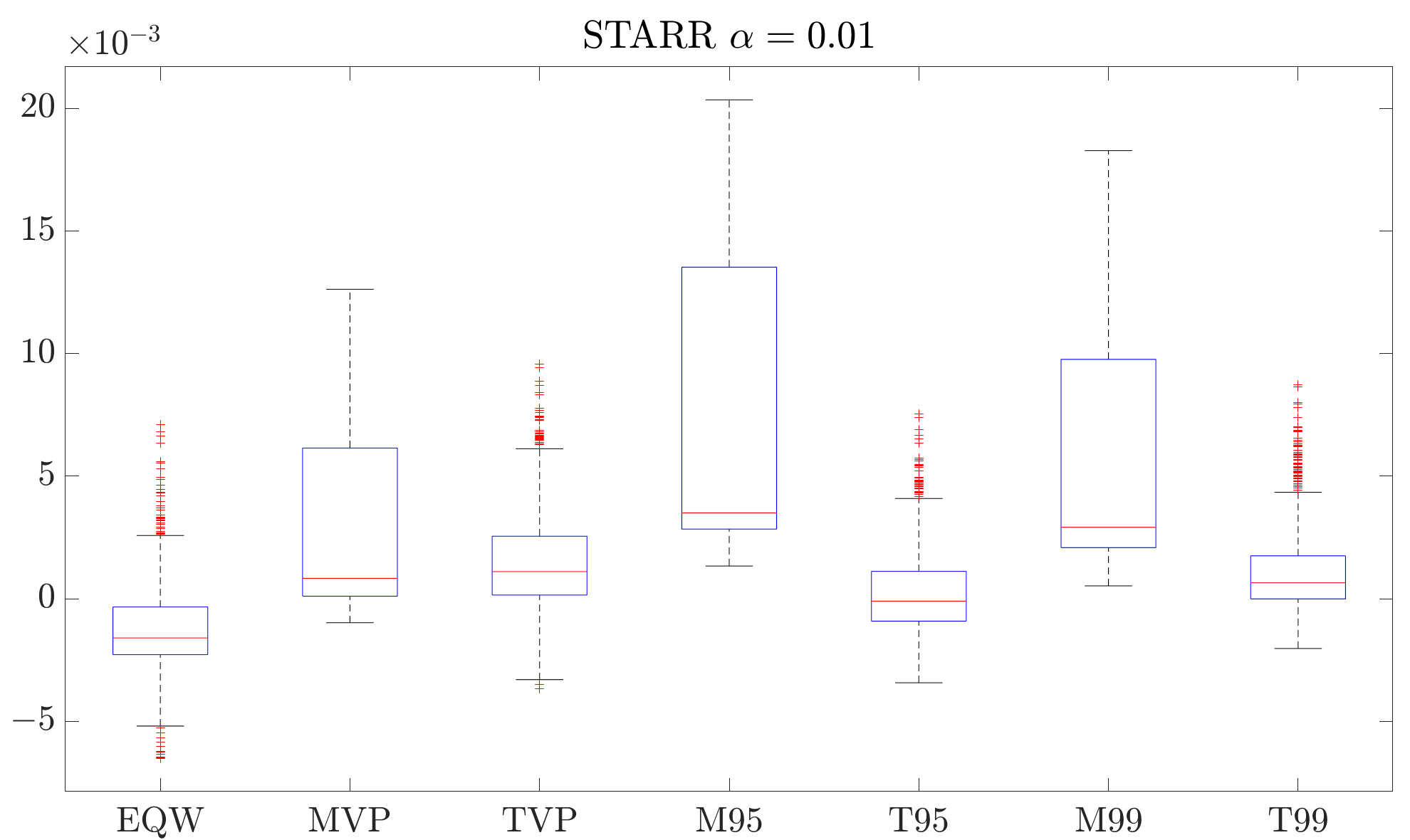}
\end{subfigure}
\hspace{0.5cm}
\begin{subfigure}{0.45\textwidth}
    \includegraphics[width=\linewidth]{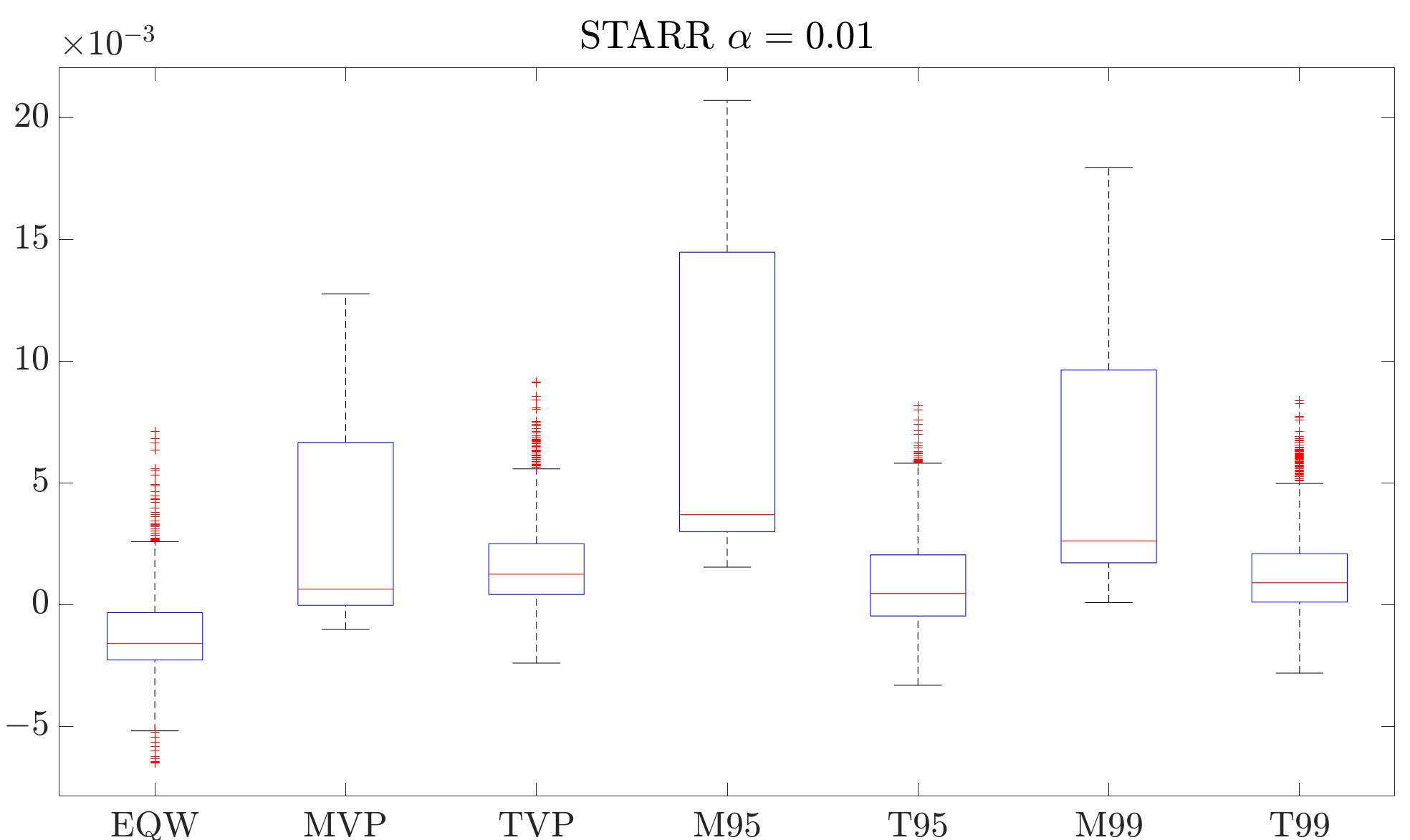}
\end{subfigure}

\caption{STARR at the 99\% confidence level for the optimized portfolios and EQW under long-only and long–short strategies. 
Each figure displays results under four configurations: 
Long–Only (top-left), Long–Short 10\% (top-right), Long–Short 20\% (bottom-left), and Long–Short 30\% (bottom-right). 
}
\label{fig:starr99_lo_ls}
\end{figure}

The STARR ratio provides a tail sensitive measure of performance by considering returns relative to downside tail risk, 
making it more robust in identifying portfolios that balance return stability with tail protection. 
In Figures \ref{fig:starr95_lo_ls} and \ref{fig:starr99_lo_ls}, 
the results at the 95\% and 99\% confidence levels respectively show similar patterns to the Sharpe ratio, 
but with greater sensitivity to downside risk. 
At the 95\% confidence level, the long only portfolios, 
particularly the MVP and M95, 
demonstrate higher median STARR values, 
indicating more consistent risk adjusted returns when accounting for moderate tail risk. 
In contrast, 
the long short portfolios display wider dispersions, 
suggesting that the introduction of leverage amplifies both gains and losses in the tails. 
When moving to the 99\% confidence level, 
the distributions become more dispersed and median STARR values decline slightly across all strategies, 
emphasizing the effect of extreme tail events. 
This pattern reflects that as the confidence level tightens, 
tail losses weigh more heavily on performance metrics. 
Overall, 
the STARR analysis confirms that while long only portfolios tend to deliver more stable tail adjusted returns, 
increasing leverage under long short strategies exposes portfolios to higher downside volatility,
particularly under extreme market conditions captured by the 99\% level.

\medskip
Together, 
the Sharpe, Rachev, and STARR ratios form a complementary set of reward–risk measures that capture different aspects of portfolio performance. 
The Sharpe ratio highlights average efficiency, 
the Rachev ratio emphasizes asymmetry in tail performance, 
and the STARR ratio balances risk and reward under extreme loss scenarios. 
According to \citet{cogneau2009ways}, 
these measures collectively provide a robust framework for assessing portfolio efficiency across both normal and turbulent markets.

\section{Tail Risk Analysis using the Hill Estimator}

Understanding tail behavior is essential in financial risk management, 
as it captures the probability and severity of extreme market movements that traditional measures such as volatility often underestimate. 
To better assess this extreme behavior, 
this study employs the Hill estimator, 
a classical and widely used method for estimating the tail index of heavy-tailed distributions.

Introduced by \citet{hill1975simple}, 
the Hill estimator provides a nonparametric approach to characterizing the heaviness of distributional tails. 
By evaluating the conditional likelihood of the extreme order statistics, 
it offers valuable insight into how rapidly the tails decay—an important factor in understanding the frequency of large losses or gains.

Suppose we have a sample of independent and identically distributed random variables whose distribution satisfies
\[
P(X > x) = x^{-\alpha}L(x), \quad \text{for } x \ge 0,
\]
where \(L(x)\) is a slowly varying function and \(\alpha > 0\) is the tail index, which measures how heavy the tail is. 
A smaller value of \(\alpha\) indicates heavier tails and a higher likelihood of extreme outcomes.

The Hill estimator is commonly expressed as
\[
\hat{\alpha}^{-1} = \frac{1}{k} \sum_{i=1}^{k} \ln\left(\frac{X_{(i)}}{X_{(k+1)}}\right),
\]
where \(X_{(1)} \ge X_{(2)} \ge \dots \ge X_{(k+1)}\) are the largest order statistics of the sample. 
In simple terms, 
the Hill estimator measures the average spacing between the largest observations relative to a threshold, 
providing a quantitative measure of how slowly the tail of the distribution decays. 
The smaller the estimated \(\alpha\), 
the heavier the tail of the return distribution, 
meaning that extreme losses or gains occur more frequently.

Over time, 
several studies have refined and extended Hill’s methodology to address its limitations and enhance robustness in financial applications. 
For example, 
\citet{wagner2000adaptive} proposed adaptive estimators to handle dependence structures such as ARCH-type processes, 
while \citet{brilhante2013generalisation} developed a generalized version to correct bias and incorporate higher-order moments. 
Similarly, 
\citet{aban2001shifted} introduced a shifted version of the estimator to achieve scale and shift invariance, 
and \citet{hill2010tail} examined its asymptotic properties under dependence and heterogeneity in financial time series.

The Hill estimator remains widely used because of its flexibility and diagnostic usefulness in empirical finance. 
\citet{drees2000hillplot} provided guidance on constructing Hill plots to visualize tail index estimates, 
while \citet{davletov2022estimating} demonstrated its practical application in analyzing the heavy-tailed behavior of equity returns, 
highlighting structural breaks during periods of financial distress such as the 2008 Global Financial Crisis.

In this study, 
the Hill estimator is applied to compare tail risk across the DJIA, EQW, and 29 individual Asian ETFs. 
Plots for EQW and DJIA are discussed in the main text to highlight differences in tail behavior, 
while the tail index estimates for the individual ETFs are presented in Appendix C to provide a broader perspective on the heterogeneity of tail behavior across the Asian ETF universe. 
In interpreting the Hill plots, 
a flatter slope corresponds to heavier tails, 
indicating a greater probability of extreme losses, 
whereas a steeper slope suggests thinner tails and more stable return distributions.

\begin{figure}[H]
\centering

\begin{subfigure}{0.45\textwidth}
    \includegraphics[width=\linewidth]{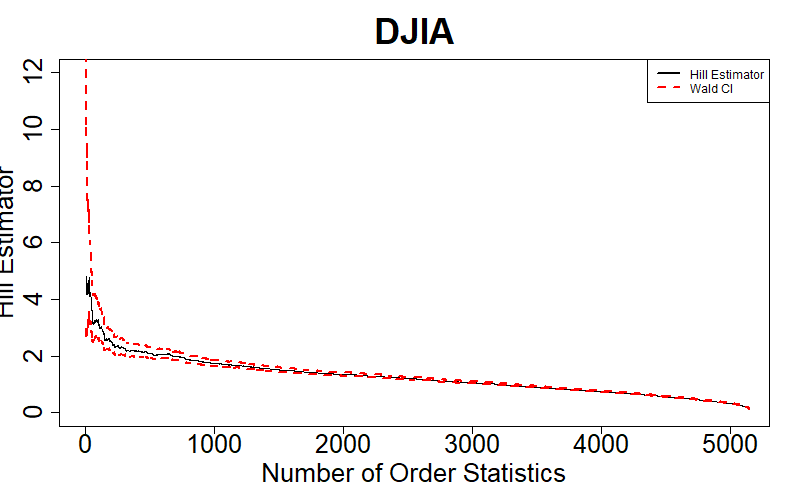}
\end{subfigure}
\hspace{0.5cm}
\begin{subfigure}{0.45\textwidth}
    \includegraphics[width=\linewidth]{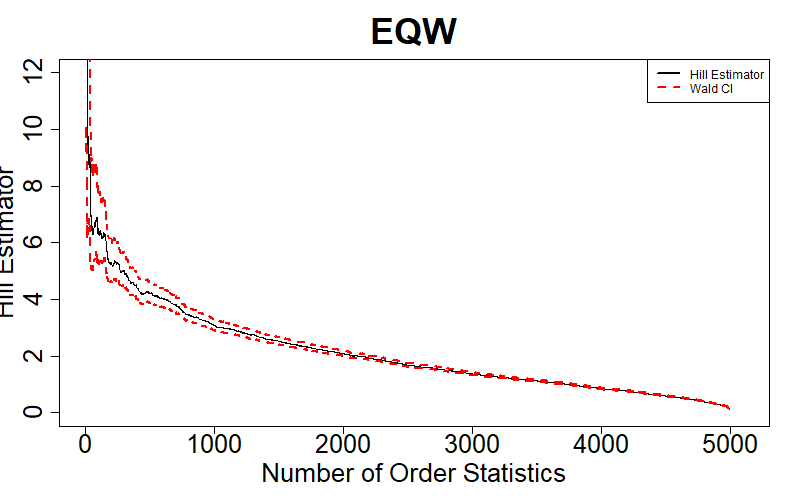}
\end{subfigure}

\caption{DJIA and EQW estimated tail index, along with the Wald confidence interval (CI)}
\label{fig:djia_eqw_hill}
\end{figure}

Figure \ref{fig:djia_eqw_hill} compares the Hill estimators for the DJIA index and the EQW portfolio, 
revealing distinct differences in tail behavior. 
Both series exhibit high initial fluctuations due to sensitivity to extreme observations, 
after which the estimates gradually stabilize. 
The DJIA converges to a higher tail index, 
indicating thinner tails and a lower probability of extreme returns. 
In contrast, 
the EQW portfolio settles at a lower tail index with wider confidence intervals, 
signifying heavier tails and greater exposure to extreme market movements.

The interpretation of these patterns relies on how the Hill estimator evolves as the number of order statistics increases. 
The Hill curve for the DJIA stabilizes around a higher value of the tail index, 
suggesting that extreme events decay more rapidly and the likelihood of very large losses is relatively low. 
Conversely, 
the EQW curve flattens at a lower level, 
indicating a slower rate of tail decay and a heavier tailed distribution. 
This means that extreme returns occur more frequently for EQW, 
reflecting its higher sensitivity to outliers and volatility shocks. 
The wider confidence bands observed for EQW further illustrate this instability, 
implying greater uncertainty in tail behavior and a larger potential range for extreme loss events.

This heavier tailed behavior in the EQW portfolio suggests a higher susceptibility to large downside losses, 
consistent with the findings from the CVaR analysis. As shown in the CVaR efficient frontiers, 
the EQW portfolio exhibits larger expected losses in the worst 1\% and 5\% of outcomes. 
The Hill plot reinforces this result by confirming that EQW’s return distribution is more prone to extreme events. 
Together, 
these findings highlight the structural distinction between the two: 
while the DJIA index’s price weighted composition benefits from the relative stability of large cap firms, 
the EQW portfolio’s equal weighting increases exposure to smaller, 
more volatile assets, 
amplifying tail risk and downside vulnerability under stressed market conditions.

\section{Conclusion}
Using 29 Asian ETFs, 
this study examines portfolio performance across the spectrum of average risk and return trade-offs and extreme tail behavior. 
We construct Markowitz and CVaR efficient frontiers and implement 
six optimized portfolios (MVP, TVP, M95, T95, M99, and T99) 
under both long only and long short configurations with leverage levels of 10\%, 20\%, and 30\%. 
The equally weighted portfolio serves as a neutral benchmark. 
Portfolio performance is evaluated using the Sharpe, Rachev, 
and STARR ratios, and extreme loss exposure is assessed through Hill tail index estimates.
The results reveal three key insights. 
First, optimization consistently enhances efficiency relative to equal weighting, 
as most individual ETFs fall below the efficient frontiers, 
while EDIV aligns closely with the tangent portfolios in both mean variance and CVaR frameworks, 
indicating strong standalone efficiency. 
Second, leverage acts as a double-edged tool: 
moderate leverage (10 to 20\%) improves returns but increases dispersion, 
whereas higher leverage (30\%) amplifies volatility and draw-downs. 
Sharpe oriented portfolios perform best when leverage is introduced, 
while minimum risk and CVaR based strategies exhibit smoother, more resilient trajectories. 
Third, tail risk analysis reinforces these findings. 
Hill plots show that the EQW portfolio settles at a lower tail index than the DJIA, 
indicating heavier tails and greater susceptibility to extreme losses. 
This pattern mirrors the CVaR results, 
where EQW demonstrates larger expected losses in the worst 1\% and 5\% of outcomes, 
confirming that systematic optimization offers superior control over downside risk and extreme market events.

These results hold clear practical implications for investors seeking Asian exposure. 
Equal weighting provides diversification but lacks efficiency compared with systematic optimization. 
In stable or trending markets, 
tangent portfolios with restrained leverage can enhance returns, 
while in volatile or stress prone periods, 
CVaR focused or minimum variance portfolios offer stronger downside protection. 
The proximity of EDIV to the tangent portfolios suggests dividend oriented ETFs can serve as efficient anchors in balanced allocations. 

Several limitations remain.
This study assumes frictionless trading and excludes borrow, financing, and shorting costs.
It also relies on rolling-window estimates that are susceptible to sampling error and model instability.
Results may vary with alternative window lengths, 
re-balancing frequencies, or estimation techniques.
In addition, 
tail index estimates depend on threshold selection and sample size, 
which may affect robustness.
Future research can address these limitations by incorporating transaction and financing costs and applying Bayesian estimators to enhance input stability.
Further extensions could explore alternative tail estimators, 
automated threshold selection, 
and the effects of currency hedging, as well as assess liquidity and implementation costs for long–short portfolios to strengthen external validity.
Moreover, 
integrating option pricing and implied volatility analysis could provide deeper insights into risk dynamics.
This remains a promising avenue for future research, 
particularly within dynamic portfolio optimization frameworks.
Despite these limitations, 
the findings presented here offers a clear and actionable map of the risk–return and tail-risk landscape for Asian ETFs.

\appendix
\appendixpage

\section{ETF Dataset}\label{app:Data} 
\begin{table}[H]
\centering
\caption{List of ETFs included in the study}
\label{tab:etfs}
\begin{tabular}{lll}
\hline
\textbf{Ticker} & \textbf{Bloomberg Code} & \textbf{ETF Name} \\
\hline
AAXJ & AAXJ US Equity & iShares MSCI All Country Asia ex Japan ETF \\
ACWX & ACWX US Equity & iShares MSCI ACWI ex U.S. ETF \\
AIA  & AIA US Equity  & iShares Asia 50 ETF \\
ASEA & ASEA US Equity & Global X FTSE Southeast Asia ETF \\
CHIQ & CHIQ US Equity & MSCI China Consumer Discretionary ETF \\
DVYA & DVYA US Equity & iShares Asia/Pacific Dividend ETF \\
DVYE & DVYE US Equity & iShares Emerging Markets Dividend ETF \\
EDIV & EDIV US Equity & SPDR S\&P Emerging Markets Dividend ETF \\
EEM  & EEM US Equity  & iShares MSCI Emerging Markets ETF \\
EEMA & EEMA US Equity & iShares MSCI Emerging Markets Asia ETF \\
EEMS & EEMS US Equity & iShares MSCI Emerging Markets Small-Cap ETF \\
EIDO & EIDO US Equity & iShares MSCI Indonesia ETF \\
EPHE & EPHE US Equity & iShares MSCI Philippines ETF \\
EWM  & EWM US Equity  & iShares MSCI Malaysia ETF \\
EWT  & EWT US Equity  & iShares MSCI Taiwan ETF \\
EWX  & EWX US Equity  & SPDR S\&P Emerging Markets Small Cap ETF \\
EWY  & EWY US Equity  & iShares MSCI South Korea ETF \\
FM   & FM US Equity   & iShares Frontier and Select EM ETF \\
GMF  & GMF US Equity  & SPDR S\&P Emerging Asia Pacific ETF \\
GXC  & GXC US Equity  & SPDR S\&P China ETF \\
HAUZ & HAUZ US Equity & Xtrackers International Real Estate ETF \\
HYEM & HYEM US Equity & VanEck Emerging Markets High Yield Bond ETF \\
INDA & INDA US Equity & iShares MSCI India ETF \\
KBWB & KBWB US Equity & Invesco KBW Bank ETF \\
KWEB & KWEB US Equity & KraneShares CSI China Internet ETF \\
MCHI & MCHI US Equity & iShares MSCI China ETF \\
THD  & THD US Equity  & iShares MSCI Thailand ETF \\
VNM  & VNM US Equity  & VanEck Vietnam ETF \\
VWO  & VWO US Equity  & Vanguard FTSE Emerging Markets ETF \\
\hline
\end{tabular}
\end{table}

\section{Risk-Adjusted Ratios for the 29 Individual Asian ETFs}

\begin{figure}[H]
\centering
\includegraphics[width=0.47\linewidth]{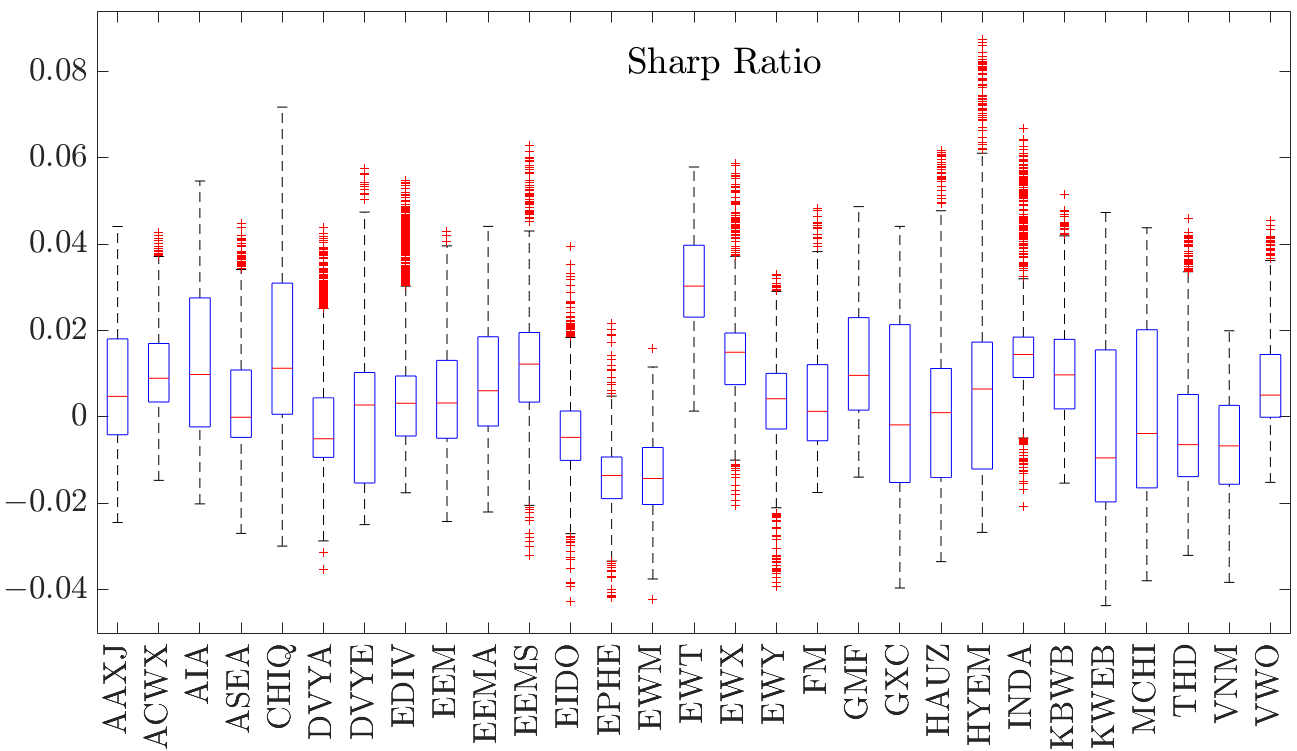}
\caption*{Figure B1: Boxplot of Sharpe ratios for the 29 individual Asian ETFs.}
\label{fig:B1_sharpe}
\end{figure}

\vspace{-0.5cm}

\begin{figure}[H]
\centering
\begin{subfigure}{0.47\textwidth}
    \includegraphics[width=\linewidth]{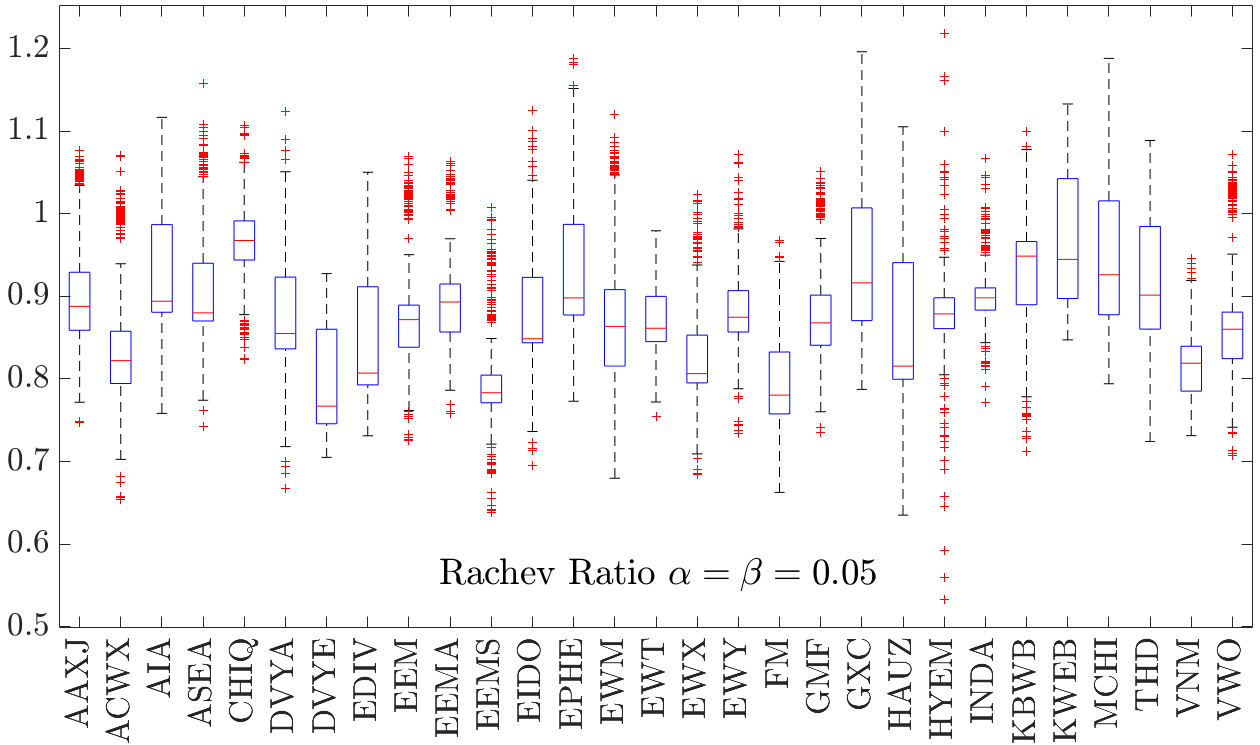}
\end{subfigure}
\hspace{0.5cm}
\begin{subfigure}{0.47\textwidth}
    \includegraphics[width=\linewidth]{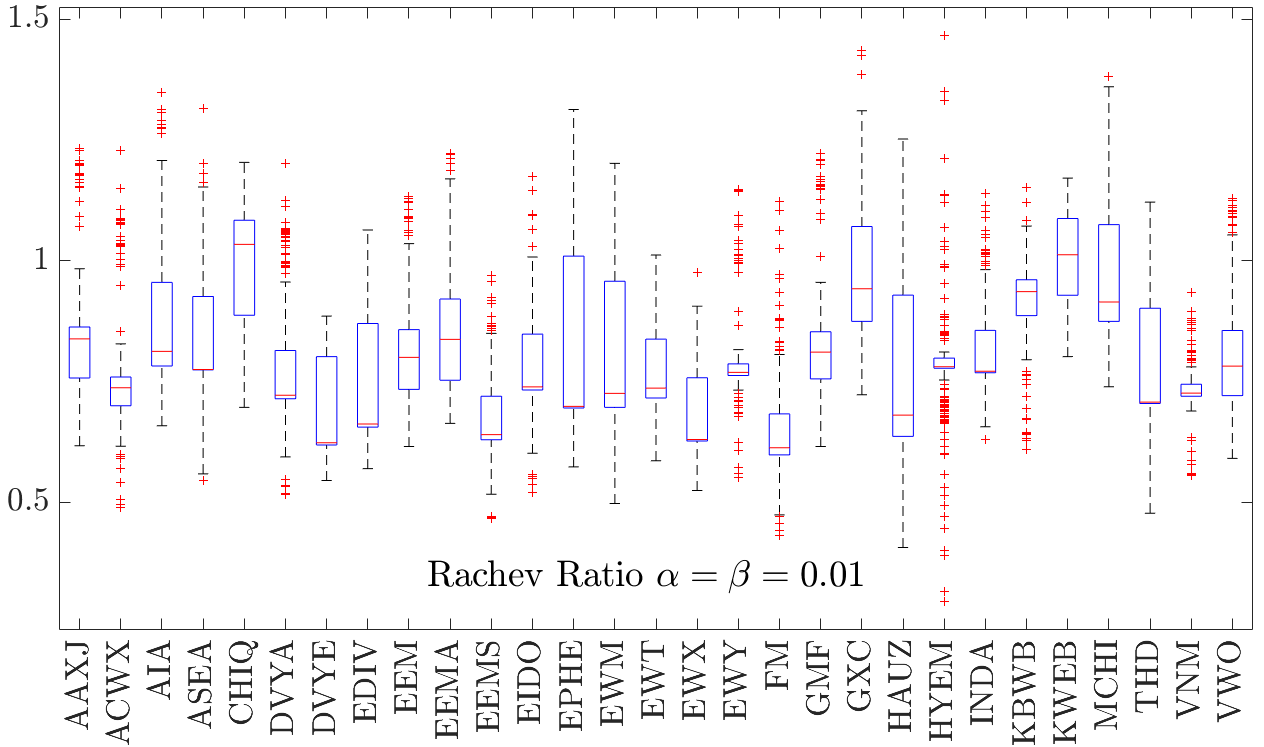}
\end{subfigure}
\caption*{Figure B2: Boxplots of RR (95\% and 99\%) for the 29 individual Asian ETFs.}
\label{fig:B2_rachev}
\end{figure}

\begin{figure}[H]
\centering
\begin{subfigure}{0.47\textwidth}
    \includegraphics[width=\linewidth]{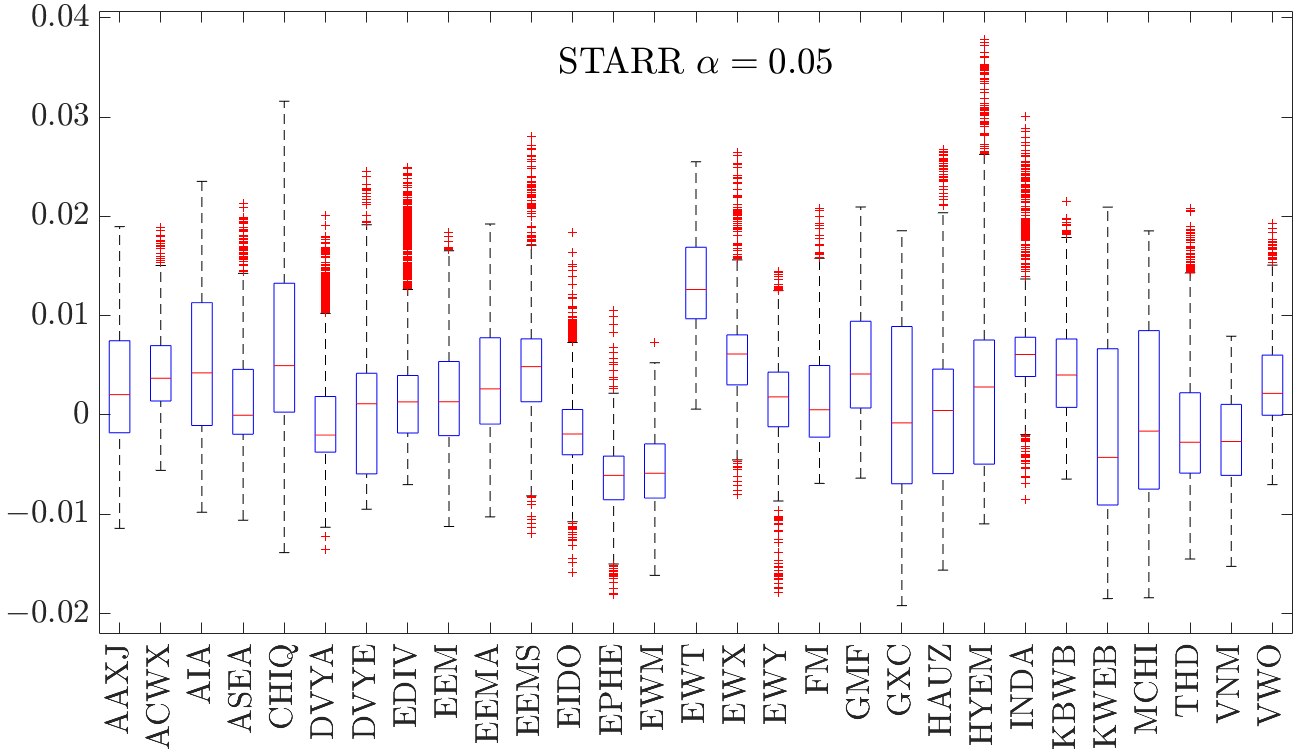}
\end{subfigure}
\hspace{0.25cm}
\begin{subfigure}{0.49\textwidth}
    \includegraphics[width=\linewidth]{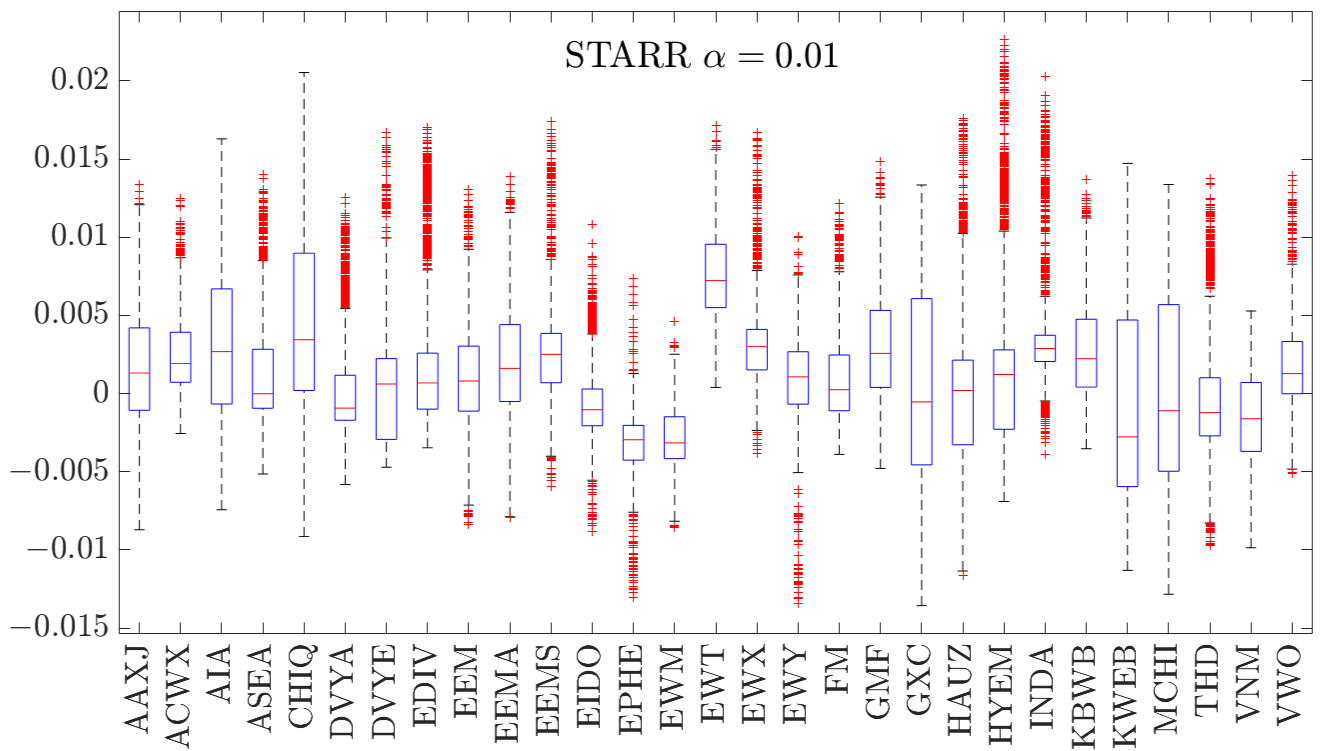}
\end{subfigure}
\caption*{Figure B3: Boxplots of STARR ratios (95\% and 99\%) for the 29 individual Asian ETFs.}
\label{fig:B3_starr}
\end{figure}

\section{Hill Estimator Plots for Individual Asian ETFs}
This appendix presents the Hill estimator results for the 29 individual Asian ETFs included in this study. Each plot illustrates the estimated tail index across the number of order statistics, along with the corresponding Wald confidence intervals. These visualizations provide deeper insight into the tail behavior and extreme risk characteristics of each ETF.

\begin{figure}[H]
\centering
\begin{subfigure}{0.3\textwidth}
    \includegraphics[width=\linewidth]{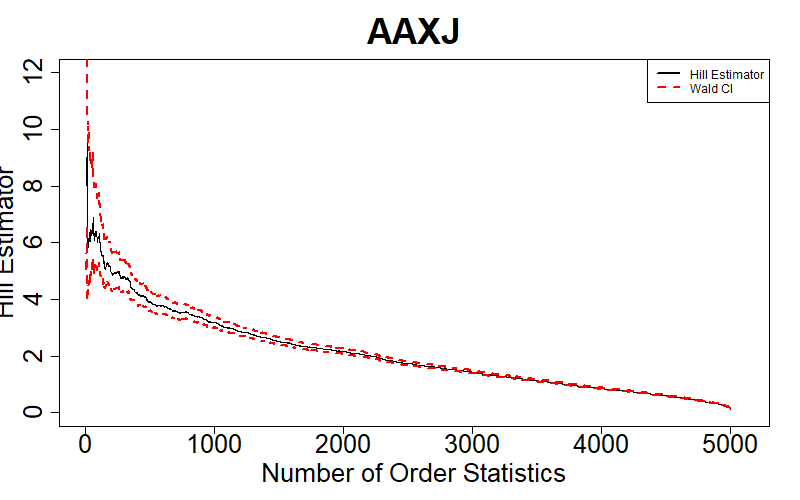}
\end{subfigure}
\hspace{0.3cm}
\begin{subfigure}{0.3\textwidth}
    \includegraphics[width=\linewidth]{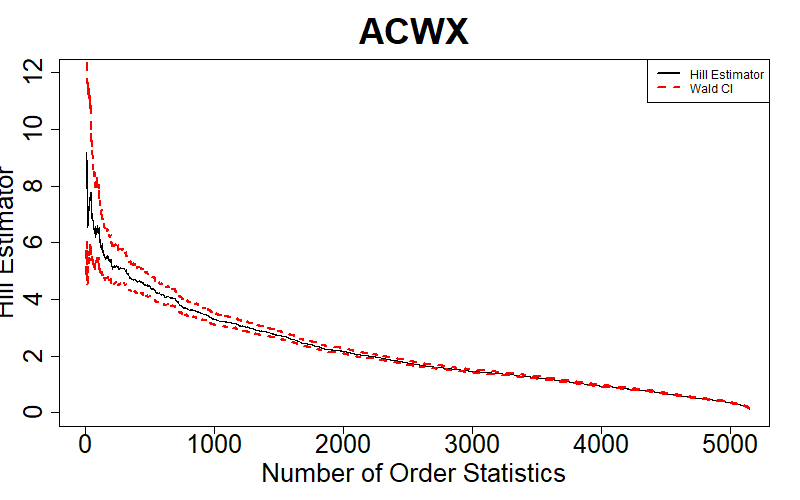}
\end{subfigure}
\hspace{0.3cm}
\begin{subfigure}{0.3\textwidth}
    \includegraphics[width=\linewidth]{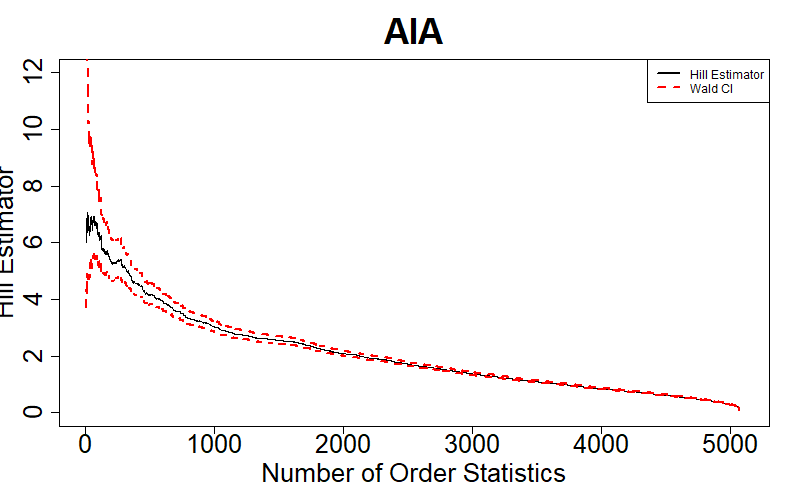}
\end{subfigure}

\vspace{0.2cm}

\begin{subfigure}{0.3\textwidth}
    \includegraphics[width=\linewidth]{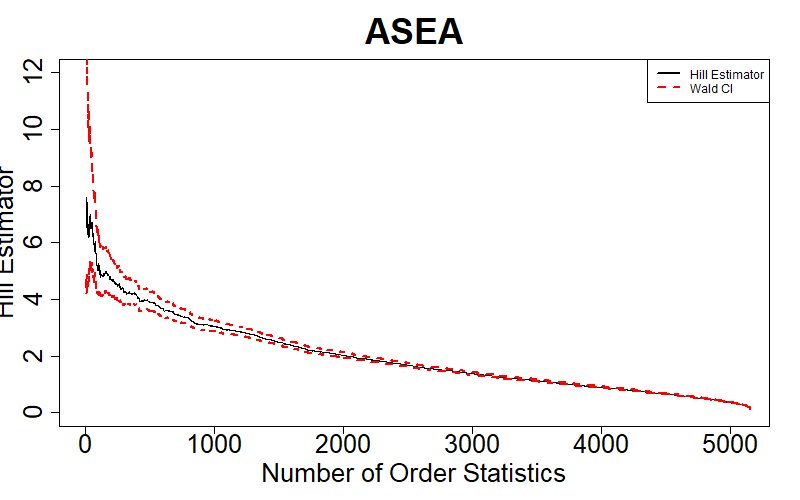}
\end{subfigure}
\hspace{0.3cm}
\begin{subfigure}{0.3\textwidth}
    \includegraphics[width=\linewidth]{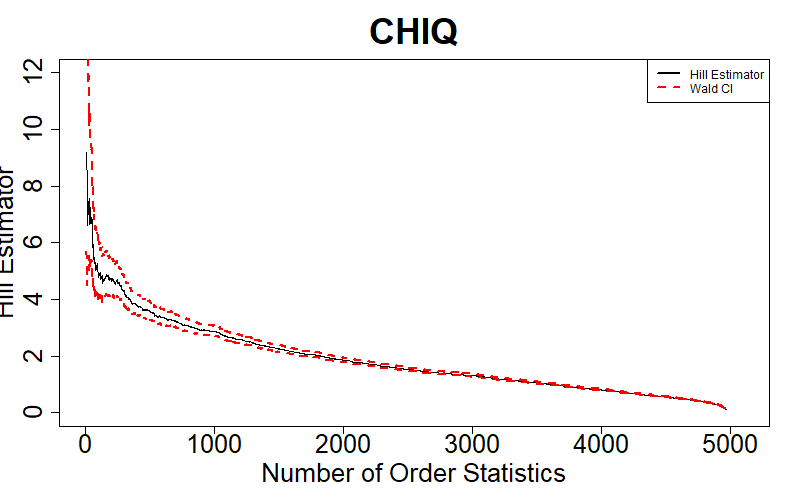}
\end{subfigure}
\hspace{0.3cm}
\begin{subfigure}{0.3\textwidth}
    \includegraphics[width=\linewidth]{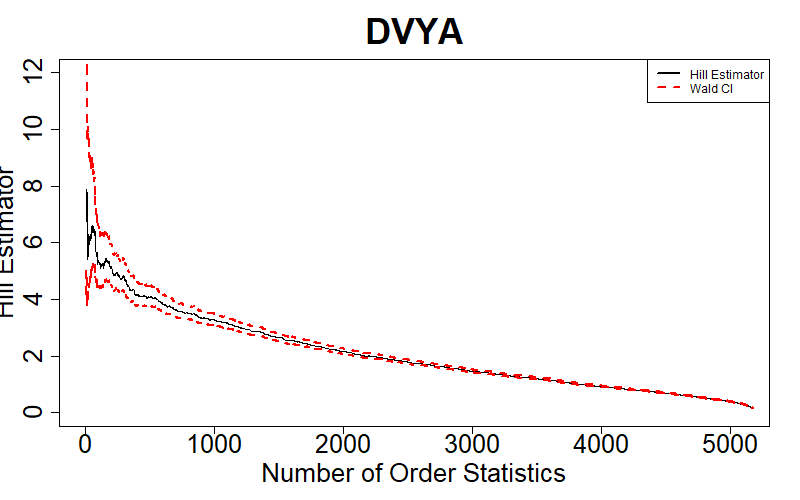}
\end{subfigure}

\vspace{0.2cm}

\begin{subfigure}{0.3\textwidth}
    \includegraphics[width=\linewidth]{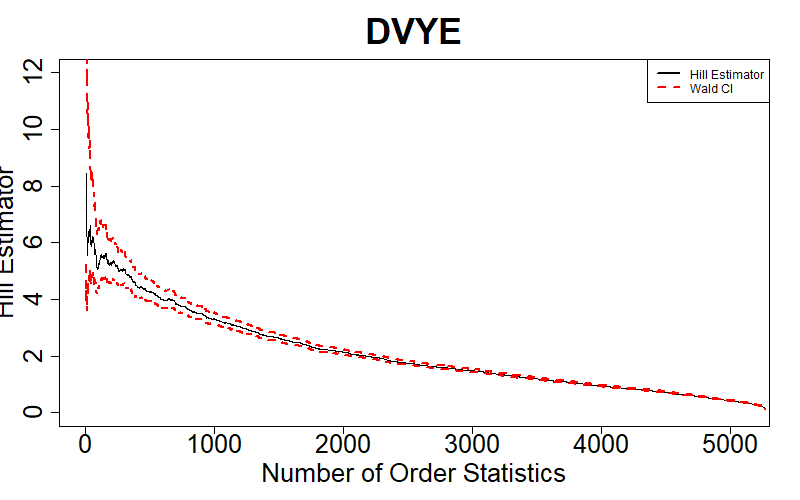}
\end{subfigure}
\hspace{0.3cm}
\begin{subfigure}{0.3\textwidth}
    \includegraphics[width=\linewidth]{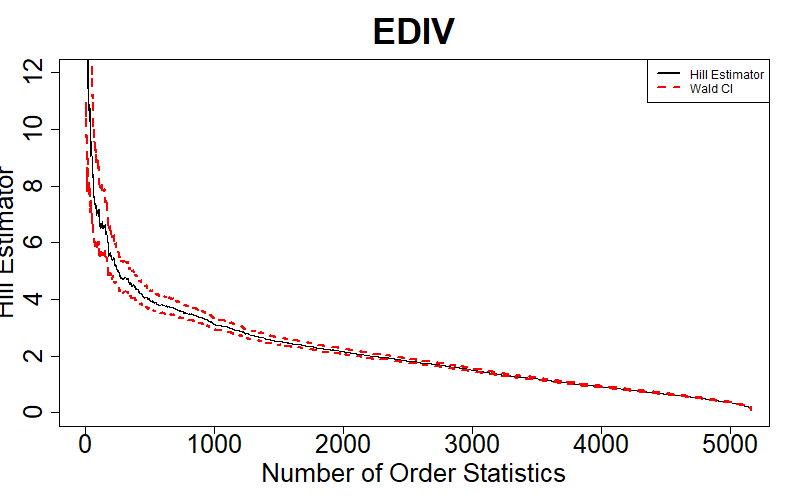}
\end{subfigure}
\hspace{0.3cm}
\begin{subfigure}{0.3\textwidth}
    \includegraphics[width=\linewidth]{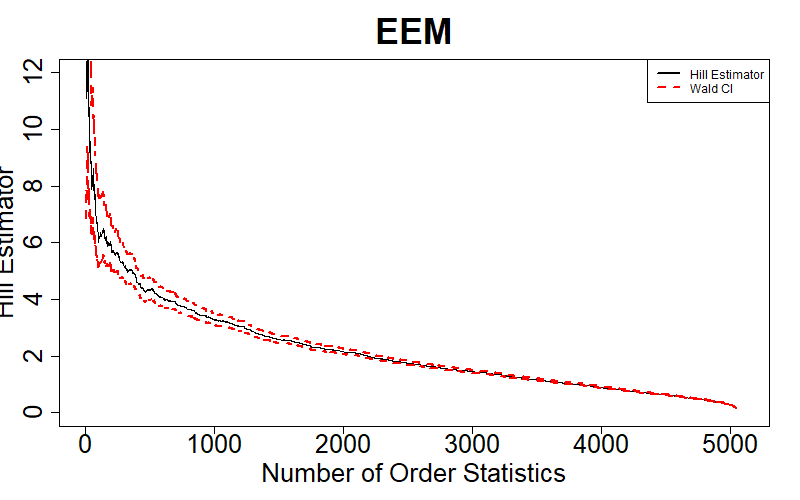}
\end{subfigure}

\vspace{0.2cm}

\begin{subfigure}{0.3\textwidth}
    \includegraphics[width=\linewidth]{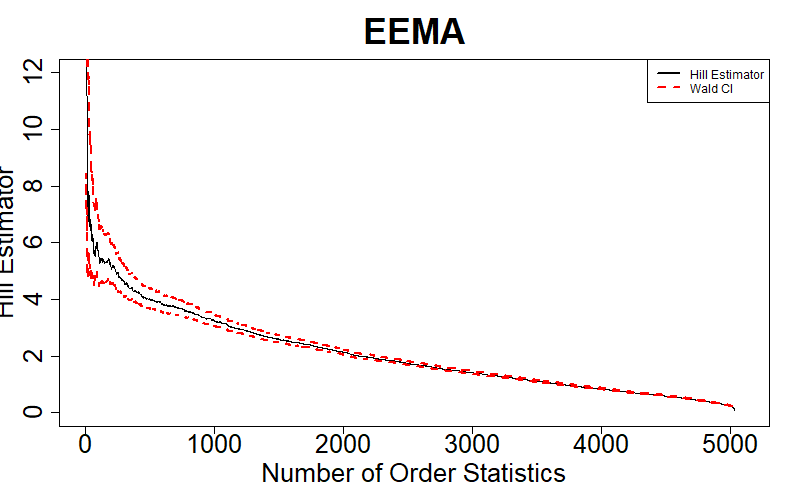}
\end{subfigure}
\hspace{0.3cm}
\begin{subfigure}{0.3\textwidth}
    \includegraphics[width=\linewidth]{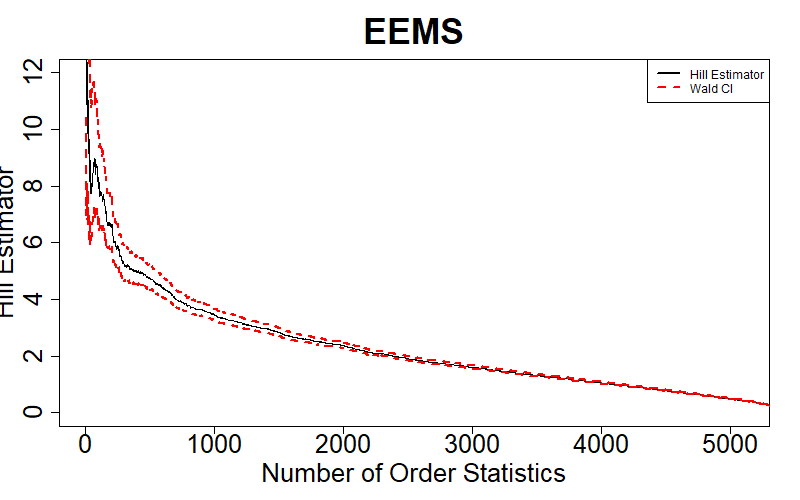}
\end{subfigure}
\hspace{0.3cm}
\begin{subfigure}{0.3\textwidth}
    \includegraphics[width=\linewidth]{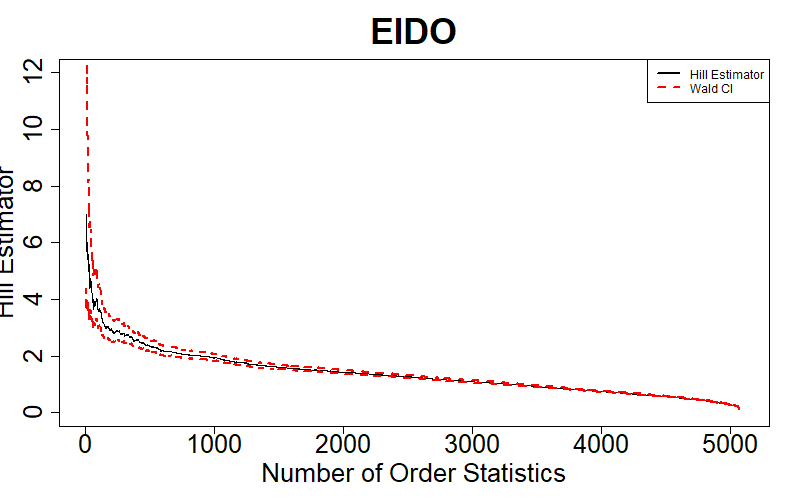}
\end{subfigure}

\vspace{0.2cm}

\begin{subfigure}{0.3\textwidth}
    \includegraphics[width=\linewidth]{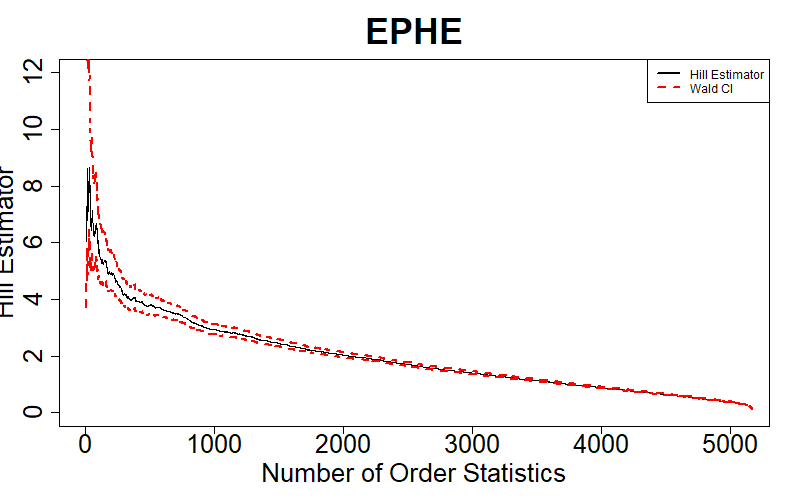}
\end{subfigure}
\hspace{0.3cm}
\begin{subfigure}{0.3\textwidth}
    \includegraphics[width=\linewidth]{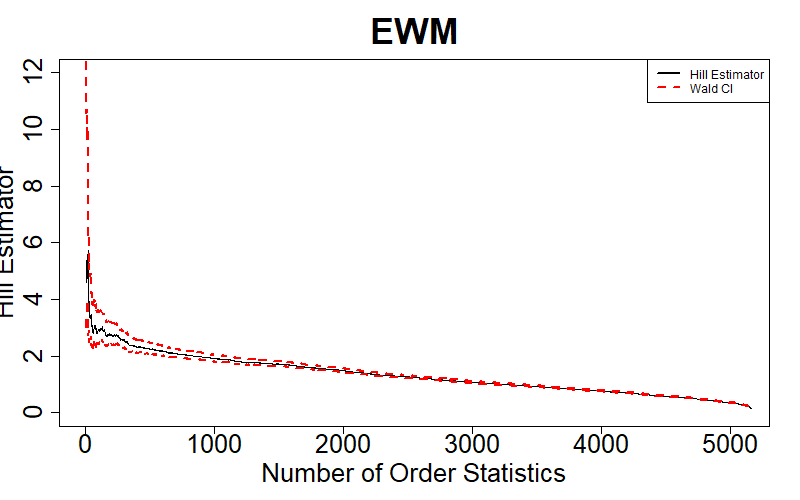}
\end{subfigure}
\hspace{0.3cm}
\begin{subfigure}{0.3\textwidth}
    \includegraphics[width=\linewidth]{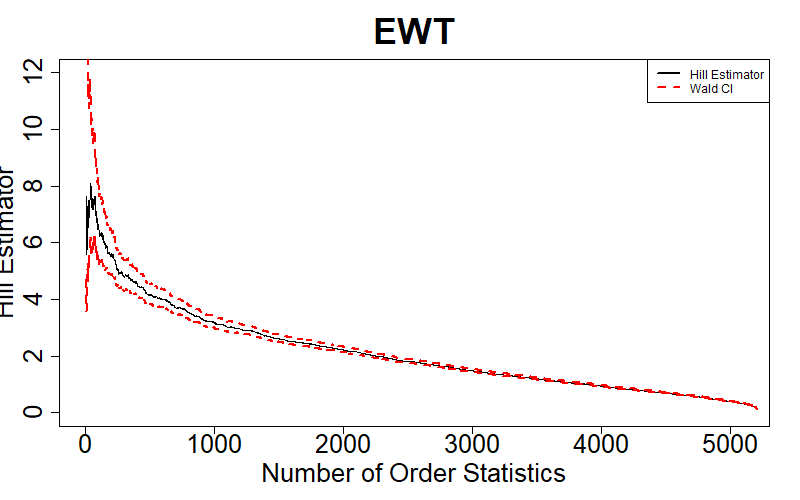}
\end{subfigure}
\vspace{0.2cm}

\begin{subfigure}{0.3\textwidth}
    \includegraphics[width=\linewidth]{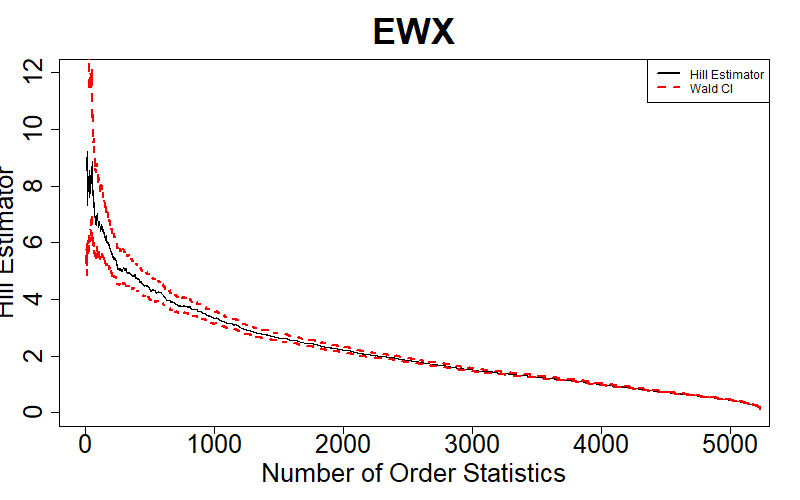}
\end{subfigure}
\hspace{0.3cm}
\begin{subfigure}{0.3\textwidth}
    \includegraphics[width=\linewidth]{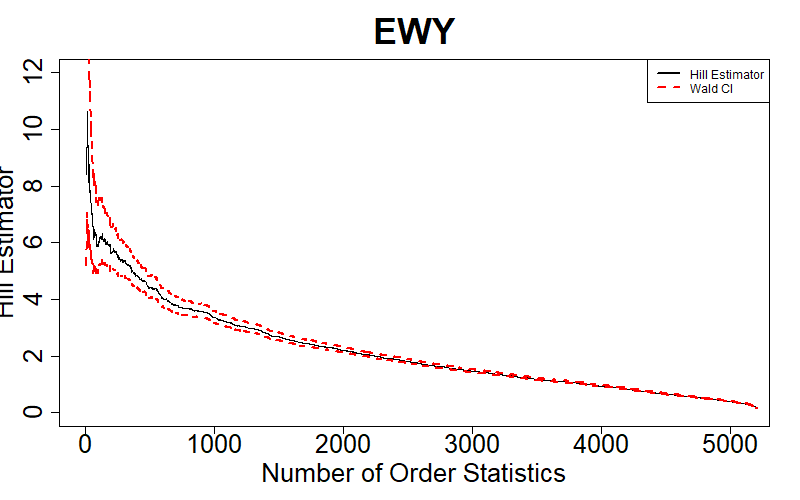}
\end{subfigure}
\hspace{0.3cm}
\begin{subfigure}{0.3\textwidth}
    \includegraphics[width=\linewidth]{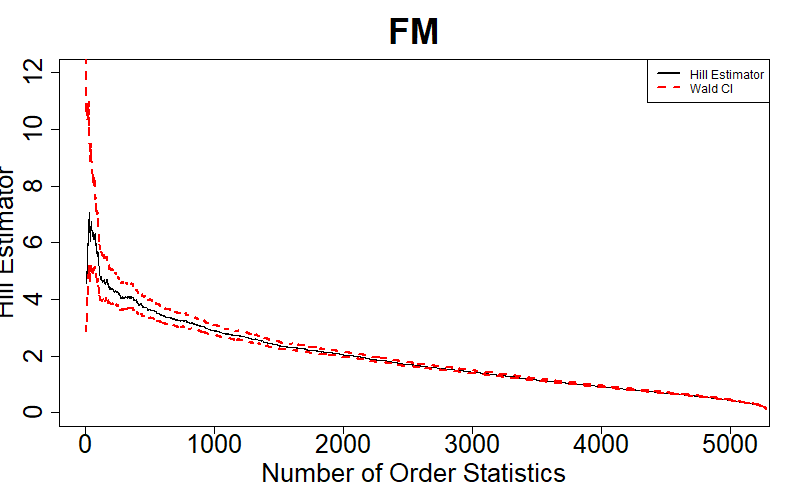}
\end{subfigure}
\label{fig:C1_hill}
\end{figure}

\vspace{0.2cm}
\begin{figure}[H]
\centering
\begin{subfigure}{0.3\textwidth}
    \includegraphics[width=\linewidth]{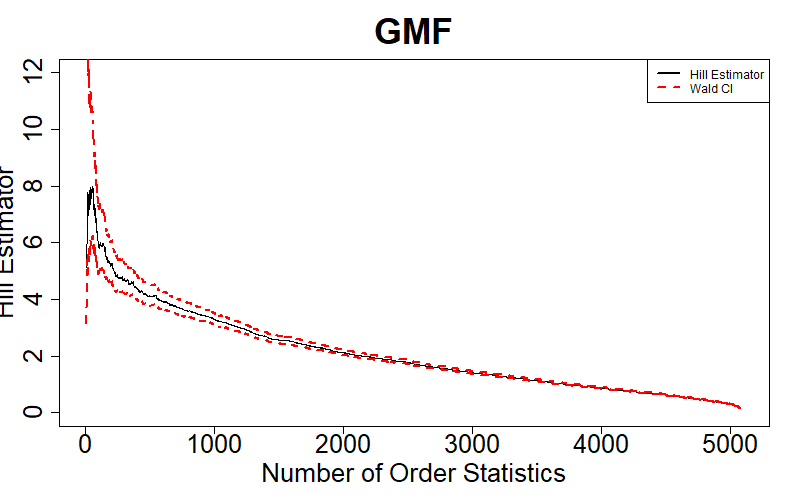}
\end{subfigure}
\hspace{0.3cm}
\begin{subfigure}{0.3\textwidth}
    \includegraphics[width=\linewidth]{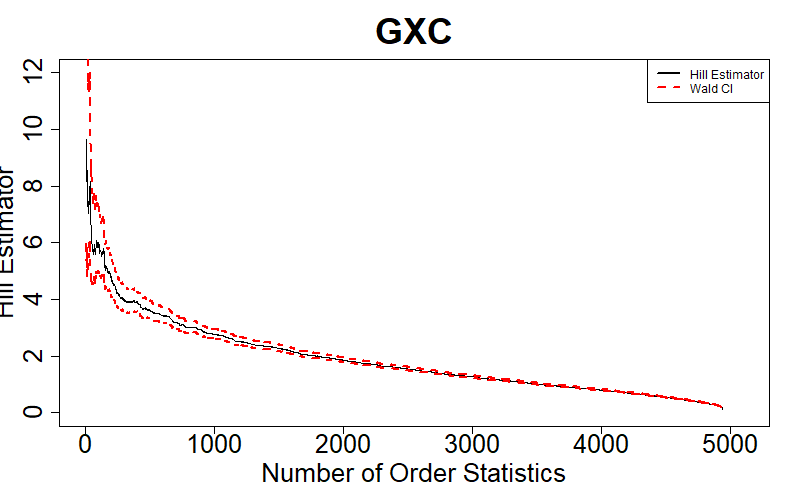}
\end{subfigure}
\hspace{0.3cm}
\begin{subfigure}{0.3\textwidth}
    \includegraphics[width=\linewidth]{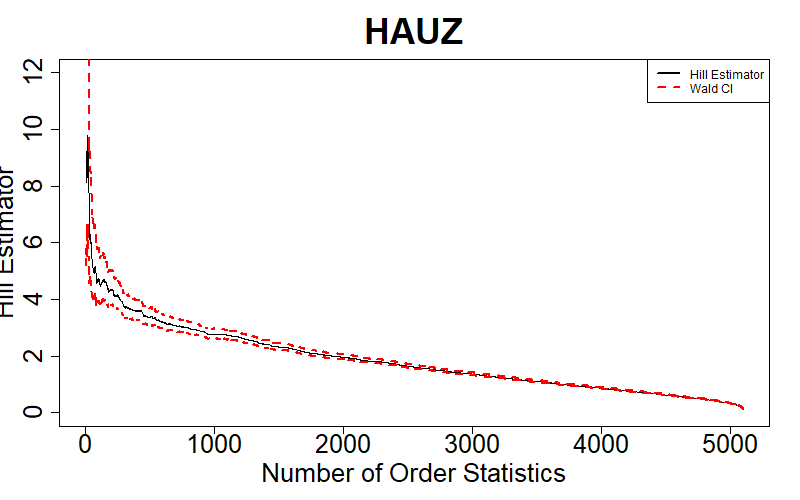}
\end{subfigure}

\vspace{0.2cm}
\begin{subfigure}{0.3\textwidth}
    \includegraphics[width=\linewidth]{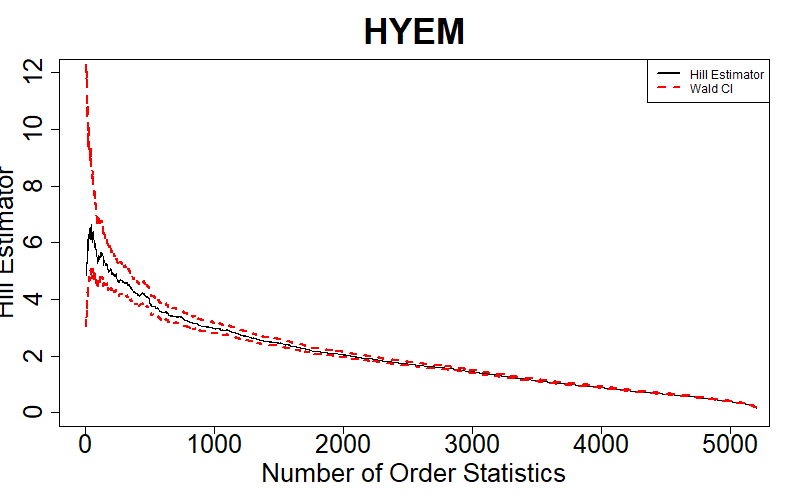}
\end{subfigure}
\hspace{0.3cm}
\begin{subfigure}{0.3\textwidth}
    \includegraphics[width=\linewidth]{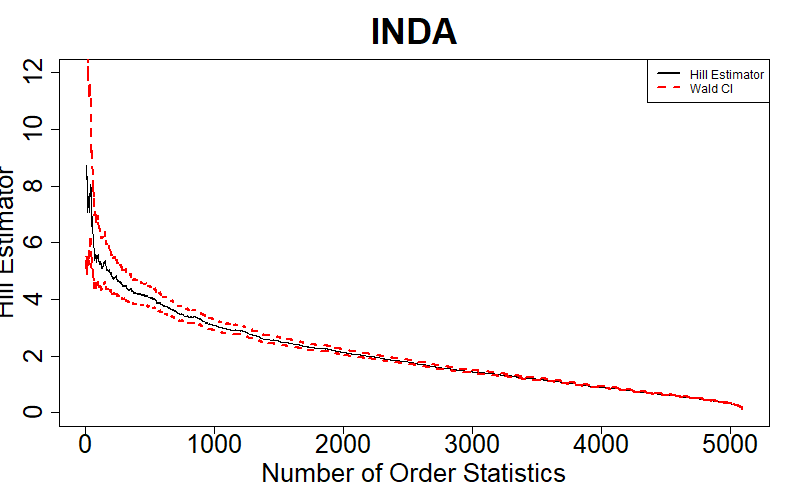}
\end{subfigure}
\hspace{0.3cm}
\begin{subfigure}{0.3\textwidth}
    \includegraphics[width=\linewidth]{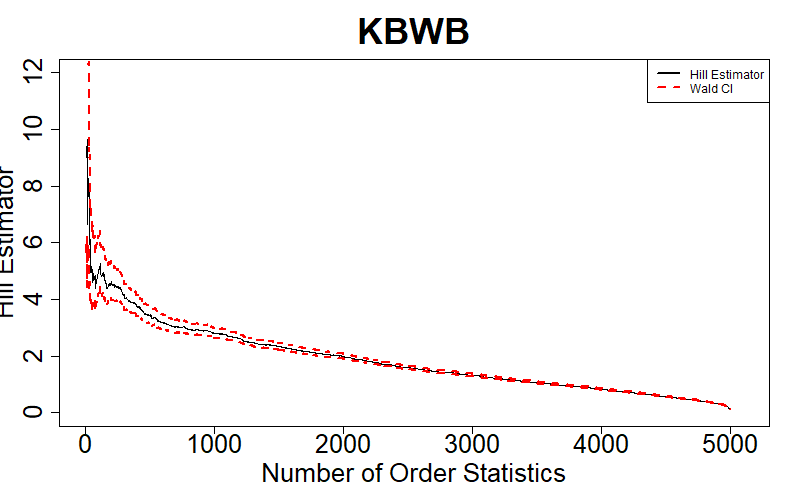}
\end{subfigure}

\vspace{0.2cm}
\begin{subfigure}{0.3\textwidth}
    \includegraphics[width=\linewidth]{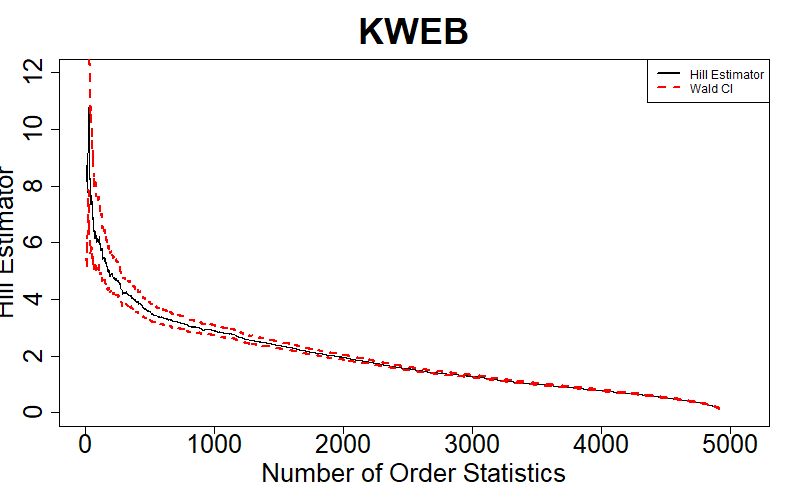}
\end{subfigure}
\hspace{0.3cm}
\begin{subfigure}{0.3\textwidth}
    \includegraphics[width=\linewidth]{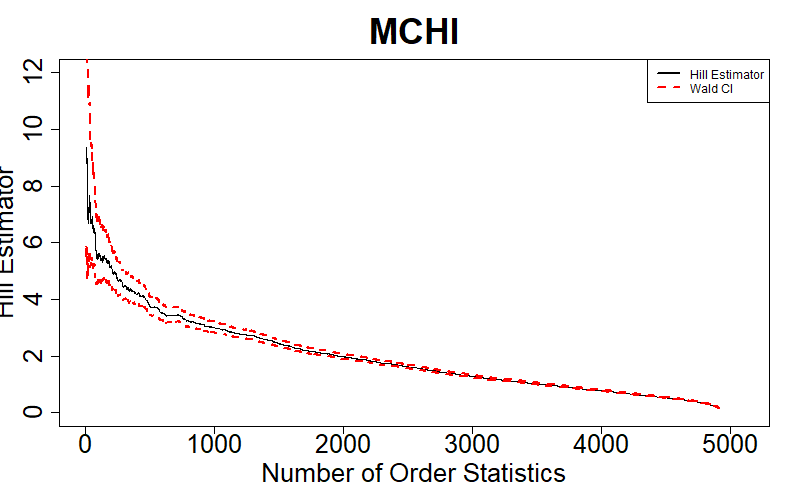}
\end{subfigure}
\hspace{0.3cm}
\begin{subfigure}{0.3\textwidth}
    \includegraphics[width=\linewidth]{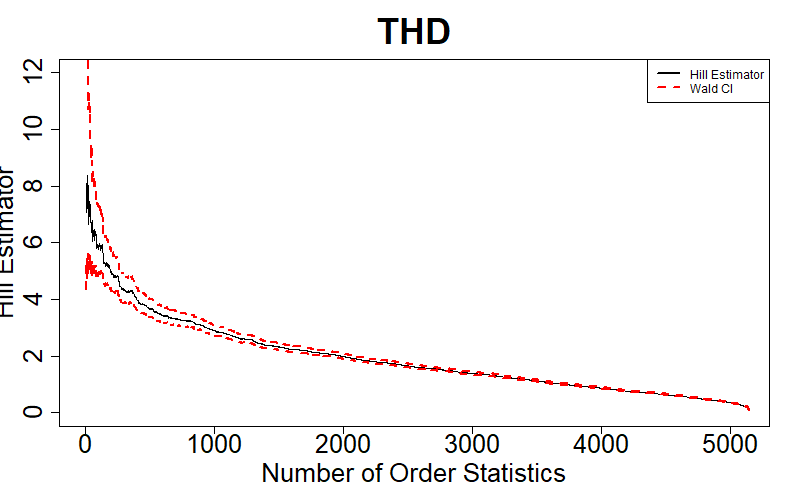}
\end{subfigure}

\vspace{0.2cm}
\begin{subfigure}{0.3\textwidth}
    \includegraphics[width=\linewidth]{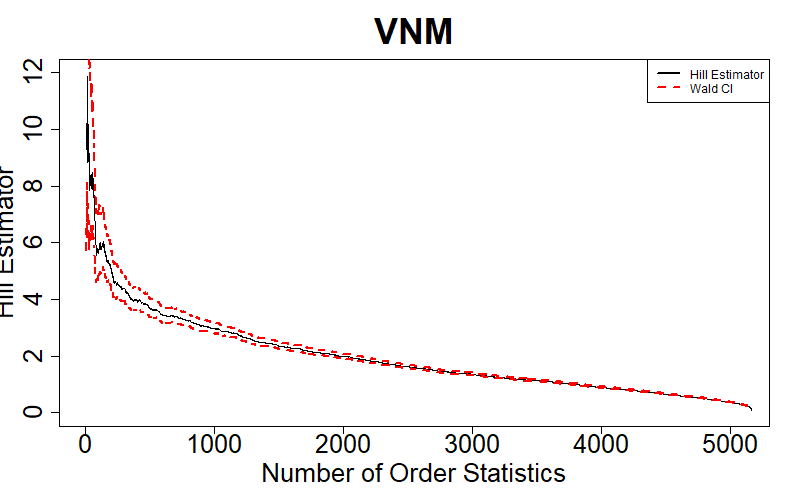}
\end{subfigure}
\hspace{0.3cm}
\begin{subfigure}{0.3\textwidth}
    \includegraphics[width=\linewidth]{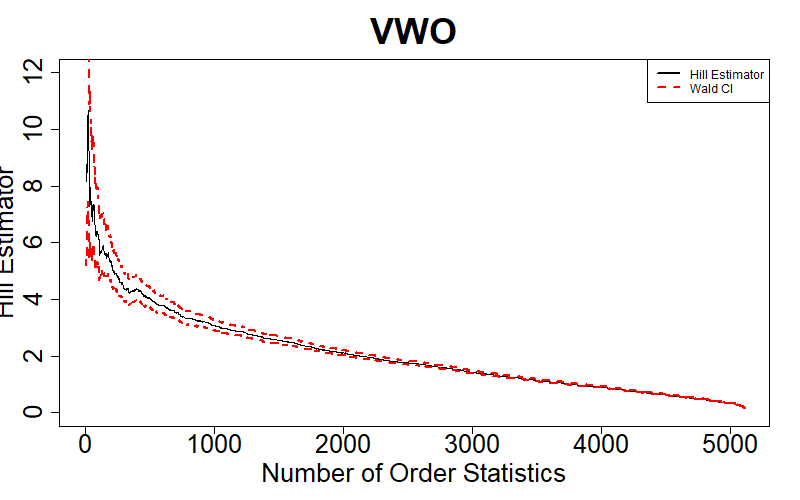}
\end{subfigure}

\vspace{0.3cm}
\caption*{Figure C1: Hill estimator plots for Asian ETFs, along with the Wald confidence interval (CI).}
\label{fig:C2_hill}
\end{figure}

\singlespacing
\normalem		
\bibliographystyle{chicago}
\bibliography{AsianOpt} 

\begin{thebibliography}{}

\bibitem[\protect\citeauthoryear{Aban and Meerschaert}{Aban and Meerschaert}{2001}]{aban2001shifted}
Aban, I.~B. and M.~M. Meerschaert (2001).
\newblock Shifted {H}ill’s {E}stimator for {H}eavy {T}ails.
\newblock {\em Communications in Statistics – Simulation and Computation\/}~{\em 30\/}(4), 949--962.

\bibitem[\protect\citeauthoryear{Bahadar, Gan, and Nguyen}{Bahadar et~al.}{2020}]{bahadar2020etf}
Bahadar, S., C.~Gan, and C.~Nguyen (2020).
\newblock Performance {D}ynamics of {I}nternational {E}xchange-{T}raded {F}unds.
\newblock {\em Journal of Risk and Financial Management\/}~{\em 13\/}(8), 169.

\bibitem[\protect\citeauthoryear{Ben-David, Franzoni, and Moussawi}{Ben-David et~al.}{2016}]{Ben-David2016}
Ben-David, I., F.~Franzoni, and R.~Moussawi (2016).
\newblock Exchange {T}raded {F}unds ({ETF}s).
\newblock Working Paper 22829, National Bureau of Economic Research.

\bibitem[\protect\citeauthoryear{Biglova, Ortobelli, Rachev, and Stoyanov}{Biglova et~al.}{2004}]{biglova2004different}
Biglova, A., S.~Ortobelli, S.~Rachev, and S.~Stoyanov (2004).
\newblock Different approaches to risk estimation in portfolio theory.
\newblock {\em The Journal of Portfolio Management\/}~{\em 31}, 103--112.

\bibitem[\protect\citeauthoryear{Blitz and Vidojevic}{Blitz and Vidojevic}{2021}]{blitz2021etf}
Blitz, D. and M.~Vidojevic (2021).
\newblock The {P}erformance of {E}xchange-{T}raded {F}unds.
\newblock {\em The Journal of Alternative Investments\/}~{\em 23\/}(3).

\bibitem[\protect\citeauthoryear{Brilhante, Gomes, and Pestana}{Brilhante et~al.}{2013}]{brilhante2013generalisation}
Brilhante, M.~F., M.~I. Gomes, and D.~Pestana (2013, January).
\newblock A simple generalisation of the {H}ill estimator.
\newblock {\em Computational Statistics \& Data Analysis\/}~{\em 57\/}(1), 518--535.

\bibitem[\protect\citeauthoryear{Cheng, Fung, and Tse}{Cheng et~al.}{2008}]{cheng2008}
Cheng, L. T.~W., H.-G. Fung, and Y.~Tse (2008).
\newblock China’s {E}xchange {T}raded {F}und: {I}s {T}here a {T}rading {P}lace {B}ias?
\newblock {\em Review of Pacific Basin Financial Markets and Policies\/}~{\em 11\/}(1), 61--74.

\bibitem[\protect\citeauthoryear{Choi, Kim, and Mitov}{Choi et~al.}{2015}]{choi2015reward}
Choi, J., Y.~S. Kim, and I.~Mitov (2015).
\newblock Reward-risk momentum strategies using classical tempered stable distribution.
\newblock {\em Journal of Banking \& Finance\/}~{\em 58}, 194--213.

\bibitem[\protect\citeauthoryear{Cogneau and H{\"u}bner}{Cogneau and H{\"u}bner}{2009}]{cogneau2009ways}
Cogneau, P. and G.~H{\"u}bner (2009).
\newblock The 101 {W}ays to {M}easure {P}ortfolio {P}erformance.
\newblock {\em SSRN Electronic Journal\/}.

\bibitem[\protect\citeauthoryear{da~Fonseca}{da~Fonseca}{2020}]{dafonseca2020performance}
da~Fonseca, J.~S. (2020).
\newblock Performance {R}atios for {S}electing {I}nternational {P}ortfolios: {A} {C}omparative {A}nalysis {U}sing {S}tock {M}arket {I}ndices in the {E}uro {A}rea.
\newblock Technical report, SSRN.
\newblock Available at SSRN: \url{https://ssrn.com/abstract=3567140} or \url{http://dx.doi.org/10.2139/ssrn.3567140}.

\bibitem[\protect\citeauthoryear{Davletov}{Davletov}{2022}]{davletov2022estimating}
Davletov, F. (2022).
\newblock Estimating the {T}ail {I}ndex of {C}onditional {D}istribution of {A}sset {R}eturns.
\newblock {\em International Journal of Financial Research\/}~{\em 13\/}(2), 21615.

\bibitem[\protect\citeauthoryear{Drees, de~Haan, and Resnick}{Drees et~al.}{2000}]{drees2000hillplot}
Drees, H., L.~de~Haan, and S.~Resnick (2000).
\newblock How to make a {H}ill plot.
\newblock {\em Annals of Statistics\/}~{\em 28}, 254--274.

\bibitem[\protect\citeauthoryear{Gatfaoui}{Gatfaoui}{2009}]{gatfaoui2009sharpe}
Gatfaoui, H. (2009).
\newblock Sharpe {R}atios and {T}heir {F}undamental {C}omponents: {A}n {E}mpirical {S}tudy.

\bibitem[\protect\citeauthoryear{Goetzmann, Ingersoll, Spiegel, and Welch}{Goetzmann et~al.}{2002}]{goetzmann2002sharpening}
Goetzmann, W., J.~Ingersoll, M.~Spiegel, and I.~Welch (2002).
\newblock Sharpening {S}harpe {R}atios.
\newblock NBER Working Paper 9116, National Bureau of Economic Research.

\bibitem[\protect\citeauthoryear{Hill}{Hill}{1975}]{hill1975simple}
Hill, B.~M. (1975).
\newblock A {S}imple {G}eneral {A}pproach to {I}nference about the {T}ail of a {D}istribution.
\newblock {\em The Annals of Statistics\/}~{\em 3}, 1163--1174.

\bibitem[\protect\citeauthoryear{Hill}{Hill}{2010}]{hill2010tail}
Hill, J.~B. (2010).
\newblock On tail index estimation for dependent, heterogeneous data.
\newblock {\em Econometric Theory\/}~{\em 26\/}(5), 1398--1436.

\bibitem[\protect\citeauthoryear{Hilliard and Le}{Hilliard and Le}{2022}]{hilliard2022}
Hilliard, J. and T.~D. Le (2022).
\newblock Exchange-{T}raded {F}unds {I}nvesting in the {E}uropean {E}merging {M}arkets.
\newblock {\em Journal of Eastern European and Central Asian Research / Manuscripts\/}~{\em 9\/}(2).

\bibitem[\protect\citeauthoryear{Jacobs, Levy, and Starer}{Jacobs et~al.}{1999}]{jacobs1999long}
Jacobs, B.~I., K.~N. Levy, and D.~Starer (1999).
\newblock Long-short portfolio management.
\newblock {\em Journal of Portfolio Management\/}~{\em 25\/}(2), 23--32.

\bibitem[\protect\citeauthoryear{Jaffri, Shirvani, Jha, Rachev, and Fabozzi}{Jaffri et~al.}{2025}]{jaffri2025}
Jaffri, A., A.~Shirvani, A.~Jha, S.~T. Rachev, and F.~J. Fabozzi (2025).
\newblock Optimizing {P}ortfolios with {P}akistan-{E}xposed {E}xchange-{T}raded {F}unds: {R}isk and {P}erformance {I}nsight.
\newblock {\em Journal of Risk and Financial Management\/}~{\em 18\/}(3), 158.

\bibitem[\protect\citeauthoryear{Jares and Lavin}{Jares and Lavin}{2004}]{jares2004}
Jares, T.~E. and A.~M. Lavin (2004).
\newblock Japan and {H}ong {K}ong {E}xchange-{T}raded {F}unds ({ETF}s): {D}iscounts, {R}eturns, and {T}rading {S}trategies.
\newblock {\em Journal of Financial Services Research\/}~{\em 25}, 57--69.

\bibitem[\protect\citeauthoryear{Joshi and Dash}{Joshi and Dash}{2024}]{joshi2024etf}
Joshi, G. and R.~K. Dash (2024).
\newblock Exchange-{T}raded {F}unds and the {F}uture of {P}assive {I}nvestments: {A} {B}ibliometric {R}eview and {F}uture {R}esearch {A}genda.
\newblock {\em Future Business Journal\/}~{\em 10}, 17.

\bibitem[\protect\citeauthoryear{Krokhmal, Palmquist, and Uryasev}{Krokhmal et~al.}{2002}]{krokhmal2002portfolio}
Krokhmal, P., J.~Palmquist, and S.~Uryasev (2002).
\newblock Portfolio optimization with conditional value-at-risk objective and constraints.
\newblock {\em Journal of Risk\/}~{\em 4\/}(2), 11--27.

\bibitem[\protect\citeauthoryear{Liebi}{Liebi}{2020}]{liebi2020}
Liebi, L. (2020).
\newblock The effect of {ETF}s on financial markets: a literature review.
\newblock {\em Financial Markets and Portfolio Management\/}.

\bibitem[\protect\citeauthoryear{Lindquist, Rachev, Hu, and Shirvani}{Lindquist et~al.}{2021}]{lindquist2021advanced}
Lindquist, W.~B., S.~T. Rachev, Y.~Hu, and A.~Shirvani (2021).
\newblock {\em Advanced REIT Portfolio Optimization: Innovative Tools for Risk Management}.
\newblock Springer.

\bibitem[\protect\citeauthoryear{Liu, Zhou, and Niu}{Liu et~al.}{2023}]{liu2023portfolio}
Liu, Y., Y.~Zhou, and J.~Niu (2023).
\newblock Portfolio optimization: A multi-period model with dynamic risk preference and minimum lots of transaction.
\newblock {\em Finance Research Letters\/}~{\em 55\/}(Part B), 103964.

\bibitem[\protect\citeauthoryear{Lo}{Lo}{2002}]{lo2002statistics}
Lo, A.~W. (2002).
\newblock The {S}tatistics of {S}harpe {R}atios.
\newblock {\em Financial Analysts Journal\/}~{\em 58\/}(4), 36--52.

\bibitem[\protect\citeauthoryear{Lo and Patel}{Lo and Patel}{2008}]{lo2008new}
Lo, A.~W. and P.~N. Patel (2008).
\newblock 130/30: The {N}ew {L}ong-{O}nly.
\newblock {\em Journal of Portfolio Management\/}~{\em 34\/}(2), 12--38.

\bibitem[\protect\citeauthoryear{Malhotra and Sinha}{Malhotra and Sinha}{2023}]{malhotra2023}
Malhotra, P. and P.~Sinha (2023).
\newblock Exchange-{T}raded {F}unds in {I}ndia {A}mid {COVID}-19 {C}risis: {A}n {E}mpirical {A}nalysis of the {P}erformance.
\newblock {\em Metamorphosis\/}~{\em 22\/}(1), 38--54.

\bibitem[\protect\citeauthoryear{Malladi and Fabozzi}{Malladi and Fabozzi}{2017}]{malladi2017}
Malladi, R. and F.~J. Fabozzi (2017).
\newblock Equal-weighted strategy: Why it outperforms value-weighted strategies? {T}heory and evidence.
\newblock {\em Journal of Asset Management\/}~{\em 18\/}(3), 188--208.

\bibitem[\protect\citeauthoryear{Markowitz}{Markowitz}{1952}]{markowitz1952portfolio}
Markowitz, H. (1952).
\newblock Portfolio selection.
\newblock {\em The Journal of Finance\/}~{\em 7\/}(1), 77--91.

\bibitem[\protect\citeauthoryear{Marszk, Lechman, and Kato}{Marszk et~al.}{2019}]{marszk2019}
Marszk, A., E.~Lechman, and Y.~Kato (2019).
\newblock {\em Exchange-{T}raded {F}unds {M}arket {D}evelopment in {A}sia-{P}acific {R}egion}, pp.\  83--142.
\newblock Cham: Springer International Publishing.

\bibitem[\protect\citeauthoryear{Martin, Rachev, and Siboulet}{Martin et~al.}{2003}]{martin2003portfolios}
Martin, R.~D., S.~Rachev, and F.~Siboulet (2003, November).
\newblock Phi-alpha optimal portfolios and extreme risk management.
\newblock {\em Willmot Magazine of Finance\/}, 70--83.

\bibitem[\protect\citeauthoryear{Mittal and Richu}{Mittal and Richu}{2017}]{mittal2017}
Mittal, S.~K. and Richu (2017).
\newblock Exchange {T}raded {F}und: A {H}istorical {R}eview.
\newblock {\em Asian Journal of Management\/}~{\em 8\/}(2), 349--352.

\bibitem[\protect\citeauthoryear{Ning, Zhao, and Jiang}{Ning et~al.}{2024}]{ning2024}
Ning, W., J.~Zhao, and F.~Jiang (2024).
\newblock {ETF}s and {T}ail {D}ependence: {E}vidence from the {C}hinese {S}tock {M}arket.
\newblock {\em Journal of International Money and Finance\/}~{\em 149}, 103194.

\bibitem[\protect\citeauthoryear{Ou}{Ou}{2023}]{ou2023portfolio}
Ou, S. (2023, October).
\newblock Portfolio {O}ptimization and {A}nalysis using {M}odern {P}ortfolio {T}heory.
\newblock {\em Dean \& Francis Academic Publishing\/}~{\em 1\/}(3).

\bibitem[\protect\citeauthoryear{Petajisto}{Petajisto}{2017}]{petajisto2017}
Petajisto, A. (2017).
\newblock Inefficiencies in the pricing of exchange-traded funds.
\newblock {\em Financial Analysts Journal\/}~{\em 73\/}(1), 24--54.

\bibitem[\protect\citeauthoryear{Rachev, Stoyanov, and Fabozzi}{Rachev et~al.}{2008}]{rachev2008advanced}
Rachev, S.~T., S.~V. Stoyanov, and F.~J. Fabozzi (2008).
\newblock {\em Advanced Stochastic Models, Risk Assessment, and Portfolio Optimization}.
\newblock Hoboken, NJ: Wiley.

\bibitem[\protect\citeauthoryear{Rockafellar and Uryasev}{Rockafellar and Uryasev}{2000}]{rockafellar2000optimization}
Rockafellar, R.~T. and S.~Uryasev (2000).
\newblock Optimization of conditional value-at-risk.
\newblock {\em Journal of Risk\/}~{\em 2\/}(3), 21--41.

\bibitem[\protect\citeauthoryear{Rompotis}{Rompotis}{2024a}]{rompotis2024a}
Rompotis, G.~G. (2024a).
\newblock Performance and {R}isk of the {US-L}isted {A}sia {P}acific {E}xchange {T}raded {F}unds.
\newblock {\em Asian Journal of Economics and Business\/}~{\em 5\/}(2), 119--157.

\bibitem[\protect\citeauthoryear{Rompotis}{Rompotis}{2024b}]{rompotis2024b}
Rompotis, G.~G. (2024b).
\newblock A {S}tudy on the {P}erformance of {J}apanese {ETF}s.
\newblock {\em Economic Analysis Letters\/}~{\em 3\/}(3), 64.

\bibitem[\protect\citeauthoryear{Salo, Doumpos, Liesi{\"o}, and Zopounidis}{Salo et~al.}{2024}]{salo2024fifty}
Salo, A., M.~Doumpos, J.~Liesi{\"o}, and C.~Zopounidis (2024).
\newblock Fifty years of portfolio optimization.
\newblock {\em European Journal of Operational Research\/}~{\em 318\/}(1), 1--18.

\bibitem[\protect\citeauthoryear{Sehgal and Mehra}{Sehgal and Mehra}{2021}]{sehgal2021robust}
Sehgal, R. and A.~Mehra (2021).
\newblock Robust reward–risk ratio portfolio optimization.
\newblock {\em International Transactions in Operational Research\/}~{\em 28\/}(4), 2169--2190.

\bibitem[\protect\citeauthoryear{Sen and Dasgupta}{Sen and Dasgupta}{2023}]{sen2023portfolio}
Sen, J. and S.~Dasgupta (2023).
\newblock Portfolio {O}ptimization: {A} {C}omparative {S}tudy.
\newblock Preprint accepted for publication in *Deep Learning - Recent Finding and Researches*, IntechOpen, London, UK, January 2024.

\bibitem[\protect\citeauthoryear{Sharpe, Alexander, and Bailey}{Sharpe et~al.}{1999}]{sharpe1999investments}
Sharpe, W.~F., G.~J. Alexander, and J.~V. Bailey (1999).
\newblock {\em Investments\/} (6th ed.).
\newblock Prentice Hall.

\bibitem[\protect\citeauthoryear{Smetters and Zhang}{Smetters and Zhang}{2013}]{smetters2013sharper}
Smetters, K. and X.~Zhang (2013).
\newblock A {S}harper {R}atio: {A} {G}eneral {M}easure for {C}orrectly {R}anking {N}on-{N}ormal {I}nvestment {R}isks.
\newblock NBER Working Paper 19500, National Bureau of Economic Research.
\newblock JEL No. G11.

\bibitem[\protect\citeauthoryear{Steiner}{Steiner}{2011}]{steiner2011sharpe}
Steiner, A. (2011).
\newblock Sharpe {R}atio {C}ontribution and {A}ttribution {A}nalysis.
\newblock Technical report, SSRN.
\newblock Available at SSRN: \url{https://ssrn.com/abstract=1839166} or \url{http://dx.doi.org/10.2139/ssrn.1839166}.

\bibitem[\protect\citeauthoryear{Taljaard and Maré}{Taljaard and Maré}{2021}]{Taljaard2021}
Taljaard, B.~H. and E.~Maré (2021).
\newblock Why has the equal weight portfolio underperformed and what can we do about it?
\newblock {\em Quantitative Finance\/}~{\em 21\/}(11), 1855--1868.

\bibitem[\protect\citeauthoryear{Tarassov}{Tarassov}{2016}]{tarassov2016}
Tarassov, E. (2016).
\newblock Exchange {T}raded {F}unds ({ETF}): {H}istory, {M}echanism, {A}cademic {L}iterature {R}eview and {R}esearch {P}erspectives.
\newblock {\em Journal of Corporate Finance Research\/}~{\em 38\/}(2), 89--108.

\bibitem[\protect\citeauthoryear{Vu and Tskhoidze}{Vu and Tskhoidze}{2021}]{vu2021etfperformance}
Vu, T. K.~L. and S.~Tskhoidze (2021).
\newblock Analysis of the {P}erformance of {ETF}s: {A} {S}tudy on the {US} {M}arket.
\newblock Student Paper.

\bibitem[\protect\citeauthoryear{Wagner and Marsh}{Wagner and Marsh}{2000}]{wagner2000adaptive}
Wagner, N. and T.~Marsh (2000).
\newblock On {A}daptive {T}ail {I}ndex {E}stimation for {F}inancial {R}eturn {M}odels.
\newblock Research Program in Finance Working Paper RPF295, University of California, Berkeley.

\bibitem[\protect\citeauthoryear{Wang, Hussain, and Adnan}{Wang et~al.}{2010}]{wang2010gold}
Wang, L., I.~Hussain, and A.~Adnan (2010).
\newblock Gold {E}xchange {T}raded {F}unds: {C}urrent {D}evelopments and {F}uture {P}rospects in {C}hina.
\newblock {\em Asian Social Science\/}~{\em 6\/}(7), 119--119.

\bibitem[\protect\citeauthoryear{Wu, Xiong, and Gao}{Wu et~al.}{2021}]{wu2021}
Wu, C., X.~Xiong, and Y.~Gao (2021).
\newblock Performance {C}omparisons between {ETF}s and {T}raditional {I}ndex {F}unds: {E}vidence from {C}hina.
\newblock {\em Finance Research Letters\/}~{\em 40}, 101740.

\bibitem[\protect\citeauthoryear{Yeniley}{Yeniley}{2025}]{Yenileg2025}
Yeniley, U. (2025).
\newblock The {R}ole of {G}old {ETF}s in {M}arket {S}tability, {P}rice {D}iscovery, and {E}conomic {D}ynamics during {C}rises.
\newblock {\em Asian Journal of Economics, Business and Accounting\/}~{\em 25\/}(1), 60--74.

\bibitem[\protect\citeauthoryear{Young and Chuahay}{Young and Chuahay}{2019}]{young2019meanvariance}
Young, M.~N. and T.~J. T.~N. Chuahay (2019).
\newblock Mean-{V}ariance {P}ortfolio {S}election {U}tilizing {E}xchange {T}raded {F}unds in {A}sia.
\newblock In {\em 2019 IEEE 6th International Conference on Engineering Technologies and Applied Sciences (ICETAS)}, pp.\  1--5.

\bibitem[\protect\citeauthoryear{Zawadzki}{Zawadzki}{2020}]{zawadzki2020etf}
Zawadzki, K. (2020).
\newblock The {P}erformance of {ETF}s on {D}eveloped and {E}merging {M}arkets with {C}onsideration of {R}egional {D}iversity.
\newblock {\em Quantitative Finance and Economics\/}~{\em 4\/}(3), 515--525.

\bibitem[\protect\citeauthoryear{Zhang}{Zhang}{2025}]{zhang2025risk}
Zhang, X. (2025).
\newblock Risk-adjusted {M}omentum {S}trategy {C}onstruction and {I}ndustry {H}eterogeneity {A}nalysis based on {STARR} {I}ndicator.
\newblock {\em Advances in Economics, Management and Political Sciences\/}~{\em 207\/}(1), 43--53.

\end{thebibliography}

\end{document}